\documentclass[12pt]{article}
\usepackage{graphicx} 
\usepackage{authblk}
\usepackage[T1]{fontenc}
\usepackage[utf8]{inputenc}
\usepackage{lmodern}
\usepackage{amsmath,amssymb}
\usepackage{geometry}
\usepackage[ruled,vlined,linesnumbered]{algorithm2e}
\usepackage[title]{appendix}
\SetKwInput{KwInput}{Input}
\SetKwInput{KwOutput}{Output}

\title{Optimization of Material Selection for Modular Buildings}

\author[1]{Marcin Anholcer\thanks{The research was supported by the SUMOC project co-financed by the National Science Centre (NCN) through grant 2024/55/I/HS4/01297, the Slovenian Research and Innovation Agency (ARIS) through grant N2-0434 and the Austrian Science Fund (FWF) through grant 10.55776/PIN4864224.}}
\author[2]{Monika Siewczyńska\thanks{This research is co-funded by the European Union under Horizon Europe FP4-Rail4EARTH project “Sustainable and green rail systems” (GA 101101917) and co-funded by the Polish Ministry of Science and Higher Education (former Ministry of Education and Science) under the “International Co-Financed Projects” programme (Agreement No. 5400/HE/2023/2).}}

\affil[1]{\scriptsize{}Pozna\'n University of Economics and Business, Institute of Informatics and Quantitative Economy}
\affil[ ]{Al.Niepodleg{\l}o\'sci 10, 61-875 Pozna\'n, Poland, \textit{marcin.anholcer@ue.poznan.pl}}
\affil[ ]{}
\affil[2]{Poznań University of Technology, Faculty of Civil and Transport Engineering, Institute of Building Engineering}
\affil[ ]{ul. Piotrowo 5, 61-138 Pozna\'n, Poland, \textit{monika.siewczynska@put.poznan.pl}}
\affil[ ]{}

\date{\today}

\usepackage[table]{xcolor}

\begin{document}

\maketitle

\begin{abstract}

We study the selection of construction materials for modular buildings with two competing objectives: the total cost and the carbon footprint (Global Warming Potential) of material production (life-cycle stages A1--A3) and of transport to the prefabrication plant. The problem is formulated as a bi-objective combinatorial problem in which one product variant is selected for each of 23 component groups. Two solution procedures are analysed: TOPSIS with a sweep of criterion weights and a greedy Pareto-sequence procedure, which generates the complete Pareto frontier of the convex hull of the feasible set of the problem. The procedures are applied to nine modular commercial and railway station buildings developed within the Rail4Earth project, using emission data from EPDs, and cost data from the Sekocenbud database, manufacturers and calibrated LLM estimates. A geometric rule applied to the Pareto sequence identifies a balanced solution for each building and a set of recommended materials common to all buildings. Across the analysed buildings, the attainable carbon footprint range is 20--26\% of the balanced value, while the corresponding cost range is below 8\%. Optimising without the transport stage changes the recommended set of materials. The results show that simple optimisation tools can support early design decisions and reduce the carbon footprint of modular buildings at an acceptable cost.
\end{abstract}

 \medskip

 \medskip


\section{Introduction}

Cost rationalisation is usually a priority for investors. However, due to the growing need to reduce CO\textsubscript{2} emissions, there is also a need to optimise aspects of the environmental impact of buildings. This paper focuses on the problem of selecting construction materials, using the example of modular buildings for small railway stations, whose concept is being developed within the Rail4Earth project.
The literature presents various approaches to building material optimisation, including multi-criteria decision-making methods, evolutionary algorithms, BIM models integrated with LCA analysis, and artificial intelligence techniques. Existing studies, however, most often focus on selected groups of materials or individual building envelope components, while comprehensive optimisation of material sets for entire buildings is addressed less frequently. An additional limitation is the lack of complete environmental data for locally available materials, which makes reliable comparative analyses difficult.
This paper presents a method for optimising the selection of construction materials for modular buildings, aimed at simultaneously minimising construction costs and the carbon footprint associated with material production and transportation. The objective of the study is to demonstrate that the use of optimisation tools can support the design decision-making process and reduce CO\textsubscript{2} emissions while maintaining acceptable investment costs.

\section{Literature review}

\subsection{LCA Analysis of Modular Buildings}

Studies reported in the literature on building design optimisation vary in scope. In many cases, optimisation involves modifying building performance parameters, material configurations of individual envelope components \cite{aljalal_integrated_2023, cheng_carbon-footprint-based_2025}, the thickness of thermal insulation in external partitions, as well as heating systems \cite{pal_life_2017} and ventilation systems \cite{aljalal_integrated_2023}, to minimise operational energy consumption. These analyses are typically integrated with environmental impact indicators derived from Life Cycle Assessment (LCA).

Material optimisation is generally conducted only to a limited extent, most often by substituting materials that have the greatest influence on LCA outcomes. This typically concerns primary structural materials \cite{santos_bim-based_2020}, such as concrete \cite{mlybari_smart_2025, gan_parametric_2023} and steel \cite{gan_parametric_2023}, and sometimes envelope elements, for instance, glass façades \cite{santos_bim-based_2020}.

An alternative approach to reducing the environmental impact of buildings is to select the most advantageous solution from a predefined set of design variants \cite{ploszaj-mazurek_artificial_2024, cheng_carbon-footprint-based_2025}. In some studies, optimisation focuses more on construction processes than on material selection \cite{zhang_multiobjective_2025}, while also demonstrating methods for incorporating transport-related CO\textsubscript{2} emissions. The scope of LCA analyses varies across studies, which complicates direct comparison of results. Many analyses adopt a cradle-to-grave perspective (A1--C4)  \cite{mlybari_smart_2025, gan_parametric_2023, figueiredo_sustainable_2021, santos_bim-based_2020}, often concluding that the operational phase dominates overall environmental performance. This may, however, obscure the contribution of the material production stage. Only a limited number of studies focus explicitly on stages A1--A3. The most commonly assessed indicator is Global Warming Potential (GWP) \cite{mlybari_smart_2025, gan_parametric_2023, santos_bim-based_2020}.

The GWP parameter is most often analysed. In addition, some studies consider other impact categories, such as operational energy demand, acidification potential, eutrophication potential, and tropospheric ozone formation. When transport-related emissions are included \cite{gan_parametric_2023, zhang_multiobjective_2025}, they are typically calculated as a function of material mass, transport distance, and emission factors associated with specific transport modes (e.g., road, rail, maritime), expressed in ton-kilometres. LCA is often combined with Life Cycle Cost Analysis (LCCA) \cite{mlybari_smart_2025, pal_life_2017, santos_bim-based_2020}, and in some cases extended to include social dimensions \cite{figueiredo_sustainable_2021}, such as working conditions, health and safety, and impacts on local communities.

A wide range of tools is employed to support environmental and cost analyses. Commonly used solutions include LCA platforms such as One Click LCA \cite{yu_data-driven_2025, aljalal_integrated_2023}, frequently integrated with Building Information Modelling (BIM) environments, particularly Autodesk Revit \cite{gan_parametric_2023, cheng_carbon-footprint-based_2025, santos_bim-based_2020}, and, in studies involving energy optimisation, DesignBuilder. Some research also utilises parametric design tools such as Grasshopper, artificial intelligence algorithms \cite{ploszaj-mazurek_artificial_2024}, and proprietary decision-support systems \cite{soust-verdaguer_utilization_2023, soust-verdaguer_supporting_2024}, especially in complex multi-criteria optimisation contexts.

Input data aggregation methods vary; however, quantities of construction materials are most often extracted directly from BIM models \cite{gan_parametric_2023, cheng_carbon-footprint-based_2025, santos_bim-based_2020}. Certain software tools, including Revit and DesignBuilder, support plugins that enable automatic transfer of quantity data to platforms such as One Click LCA, facilitating process automation. In early design stages, conceptual BIM models may also support the optimisation of LEED certification scores \cite{bank_decision-making_2011}, where the transport distance of construction materials is a relevant criterion. In earlier versions of the U.S. Green Building Council's LEED system, transport distance was evaluated directly in kilometres, promoting the use of regional materials within an 800 km radius. Currently, materials transported within 160 km are classified as local \cite{chen_integrate_2015}.

Environmental parameters of materials are primarily derived from Environmental Product Declarations (EPDs). In their absence, they are supplemented with literature data, expert judgement \cite{mlybari_smart_2025}, or averaged values \cite{gan_parametric_2023}, such as emission factors for standard materials and processes, including transport. Cost data remain highly heterogeneous and difficult to compare, as they depend on both temporal and location factors. They are obtained from LCA tool databases (e.g., One Click LCA) \cite{yu_data-driven_2025}, market and cost estimation data \cite{santos_bim-based_2020, cheng_carbon-footprint-based_2025}, literature sources \cite{pal_life_2017}, or expert assessments \cite{mlybari_smart_2025}, and are often based on averaged assumptions reflecting regional conditions.

The literature review reveals several research gaps. The scope of optimisation in existing studies is highly variable and often insufficiently focused on the selection of structural materials, particularly for stages A1--A3, i.e., the production of construction materials and their transport to the construction site (A4). Structural components and building envelope elements are characterised by high durability and are only partially replaced during the building life cycle, in contrast to HVAC and MEP systems, which can be more easily upgraded to more energy-efficient or lower-carbon alternatives. Consequently, decisions regarding material selection are crucial and should also incorporate economic considerations.

Material optimisation is typically limited in scope, usually addressing individual materials or a small number of configurations. Moreover, the reliability of LCA results depends on the quality of input data. Currently, in many countries -- including Poland -- there are no regulatory requirements forcing EPDs for construction materials. As a result, analyses often rely on averaged data, which introduces more uncertainty.

One of the objectives of the Rail4Earth project, `Sustainable and green rail systems', is to develop a catalogue of modular (volumetric) railway station buildings for small stations and to demonstrate that the proposed solution enables at least a 30\% reduction in carbon footprint relative to the reference object. Achieving this goal requires using construction materials with lower GWP values. Material selection is based on EPD data. However, for many material categories, such declarations are not available for products manufactured in Poland, necessitating the use of materials sourced from other European countries. In such cases, it is essential to account for both emissions and transport costs associated with delivering materials to the construction site or prefabrication facility in Poland.

Due to the specific procurement model for railway station construction by the State Treasury company (PKP S.A.), which involves separate tenders for individual station buildings, it is not possible to optimise transport load sizes for construction materials. This constraint must also be taken into account in the analyses.

\subsection{Optimisation}

Various optimisation techniques have been developed to support the selection of optimal building materials, particularly in sustainable construction. These approaches typically combine multi-criteria decision-making (MCDM), life-cycle assessment, and advanced optimisation algorithms to address complex and uncertain design environments.

A large group of studies relies on multi-criteria decision-making frameworks, often combined with fuzzy logic to capture uncertainty in expert judgments. For instance, Aghazadeh et al. \cite{aghazadeh_hybrid_2022} propose a hybrid fuzzy MCDM framework (Delphi-SWARA-ARAS) for optimising structural systems compatible with sustainable materials in mass housing. The method reduces an initial set of 42 indicators to 14 key criteria using the fuzzy Delphi method, then applies fuzzy SWARA for weighting and fuzzy ARAS for ranking. This allows effective handling of uncertainty while integrating technical, economic, environmental, and socio-cultural dimensions, finally identifying LSF, ICF, and PRC systems as the best options.

Similarly, Mathiyazhagan et al. \cite{mathiyazhagan_sustainable_2019} present a hybrid framework integrating the Best-Worst Method (BWM) and Fuzzy TOPSIS. Optimisation is achieved by first deriving consistent weights for 23 sustainability criteria using BWM and then ranking alternatives under uncertainty through Fuzzy TOPSIS. This structured approach reduces decision complexity while maintaining robustness and consistency, with fly ash brick identified as the optimal solution.

Fazeli et al. \cite{fazeli_bim-integrated_2022} also employ fuzzy TOPSIS within a BIM-integrated optimisation framework at the conceptual design stage. Their model transforms subjective user preferences into weighted criteria using fuzzy logic, then ranks alternatives based on proximity to an ideal solution. The integration of BIM data, fuzzy weighting, and TOPSIS enables automated, uncertainty-aware multi-criteria optimisation, with results further validated through life-cycle cost (LCC) analysis.

A related but different approach is proposed by the authors of \cite{figueiredo_sustainable_2021}, who combine BIM, life-cycle sustainability assessment (LCSA), and fuzzy AHP. The method derives criteria weights under uncertainty and leverages BIM data to improve the accuracy of life-cycle impact calculations. This integration enables systematic optimisation that joins environmental, economic, and social objectives. However, results, as usual, depend on data completeness.

Other studies extend MCDM frameworks by combining multiple methods or incorporating dynamic system behaviour. Mlybari \cite{mlybari_smart_2025} develops an integrated LCA-MCDM framework combining AHP-TOPSIS and VIKOR, supported by Monte Carlo uncertainty analysis for robustness. Osama et al. \cite{osama_analyzing_2025} integrate system dynamics (SD) with BWM-TOPSIS, enabling optimisation by dynamically simulating interactions among sustainability criteria before ranking alternatives. Ruslan et al. \cite{ruslan_value-based_2021} propose a value-based optimisation framework that combines AHP-derived weights with function-over-cost evaluation, explicitly balancing performance and lifecycle costs while accounting for stakeholder priorities.

Additionally, Gurupatham et al. \cite{gurupatham_building_2024} employ AHP to derive weights for ecological, economic, social, and cultural criteria, which are then used in TOPSIS to rank material alternatives. Their results emphasise the dominance of ecological factors (46\%), followed by economic (30\%), social (19\%), and cultural (5\%), enabling systematic optimisation across multiple sustainability dimensions.

Another approach focuses on integrating life-cycle cost (LCC), environmental impacts, and carbon footprint into unified optimisation frameworks. In \cite{abouhamad_framework_2019}, a life-cycle-based framework is developed that integrates BIM, energy simulation, and end-of-life modelling to minimise total LCC. The approach includes Monte Carlo sensitivity analysis to identify key cost factors and uses sustainability metrics such as recycling and waste management, enabling robust early-stage decision-making under uncertainty.

Cheng et al. \cite{cheng_carbon-footprint-based_2025} extend this idea by introducing a Carbon Footprint-Integrated LCC model that combines construction costs, carbon footprint, and waste disposal into a single objective function. Optimisation is performed using the Symbiotic Organic Search (SOS) 2.0 algorithm, which identifies the configuration that minimises total life-cycle cost across different carbon-pricing scenarios. Sun et al. \cite{sun_selection_2020} further develop life-cycle optimisation by explicitly introducing uncertainty. Their framework combines LCA and economic analysis, first using deterministic optimisation and then extending it to a robust formulation that accounts for uncertainties in material properties and operational conditions. This enables optimisation of insulation material type and thickness while guaranteeing performance under worst-case scenarios.

Many studies apply advanced multi-objective optimisation algorithms, particularly evolutionary methods, to explore trade-offs between competing objectives. Ansah et al. \cite{ansah_integrated_2022} propose a BIM-based lifecycle optimisation framework using NSGA-II to minimise both environmental (LCA) and economic (LCC) impacts. The approach integrates embodied and operational effects into a single process, generating Pareto-optimal solutions and avoiding trade-offs between lifecycle stages.

Similarly, Pal et al. \cite{pal_life_2017} use NSGA‑II to simultaneously minimise life-cycle cost and carbon footprint, demonstrating how full lifecycle integration influences cost-carbon trade-offs. Zhang et al. \cite{zhang_multiobjective_2025} apply NSGA‑II within a 6D BIM environment to optimise construction time, cost, and carbon emissions, enabling dynamic visualisation of trade-offs in project execution strategies. Shi et al. \cite{shi_multi-objective_2025} propose a two-phase optimisation framework for early-stage design that joins structural selection with parametric optimisation. Using NSGA‑II, the approach captures trade-offs between embodied carbon, energy use, and daylighting performance, highlighting the importance of early design decisions.

More advanced hybrid approaches integrate machine learning with optimisation algorithms. Bai and Xue \cite{bai_energy_2024} develop an ANN-NSGA-II framework where a neural network model replaces computationally expensive simulations, enabling fast exploration of design alternatives. Asghari et al. \cite{asghari_machine_2026} further develop this concept by combining GWO-based geometry optimisation, EnsembleXGBoost-SHAP sensitivity analysis, and NSGA‑II. Their approach uses explainable machine learning to identify key variables and reduce the search space, improving efficiency and facilitating interpretation.

In contrast to heuristic algorithms, some studies employ rigorous mathematical programming techniques. Kilis et al. \cite{kilis_multi-criteria_2023} develop a multi-criteria optimisation framework using mixed-integer nonlinear programming (MINLP) implemented in GAMS and solved with the BARON optimiser. The model simultaneously evaluates thermal performance, cost, and environmental impact, treating material and design parameters as decision variables. By applying multi-objective formulations and goal programming, the framework systematically explores trade-offs and identifies balanced solutions through weighting and sensitivity analysis.

The reviewed approaches demonstrate a clear trend towards integrating multiple sustainability dimensions -- environmental, economic, and social -- into unified optimisation frameworks. Methods combine MCDM techniques, life-cycle assessment, and advanced optimisation algorithms, often enhanced by BIM and machine learning, to support robust, uncertainty-aware decision-making for sustainable material selection.

More papers touching on the problem of optimal selection of building materials can be found, e.g., in the following two studies.

Shi et al. \cite{shi_assessment_2025} offer a comprehensive review of optimisation approaches for building material selection, with particular emphasis on embodied carbon emissions (ECE). Their analysis highlights a clear transition from traditional life-cycle assessment (LCA)-based methods toward hybrid, data-driven, and multi-objective optimisation frameworks. In this context, optimisation is increasingly achieved through machine learning (ML) techniques that enable dynamic performance prediction, uncertainty quantification, and automated exploration of design alternatives. At the same time, multi-criteria decision-making (MCDM) approaches remain essential for balancing competing objectives such as carbon emissions, cost, and performance. The review emphasises the increasing use of hybrid ML-MCDM frameworks that integrate predictive modelling with structured decision-making while accounting for regional variability and full lifecycle impacts. This change reflects a move toward more adaptive, data-rich, and uncertainty-aware optimisation paradigms.

A complementary perspective is provided by Soliman et al. \cite{soliman_sustainable_2026}, who conduct a meta-analysis of the foundations of optimisation for sustainable material selection, focusing in particular on the role of criterion weighting within MCDM frameworks. In their analysis, optimisation is described as a weighting problem in which the relative importance assigned to environmental, technical, economic, and social criteria directly determines the resulting material choices. Reviewing 15 optimisation frameworks, the authors identify a consistent prioritisation pattern in which environmental criteria dominate (36.4\%), followed by technical, economic, and social factors. This shows a systematic bias in current optimisation practices and explains variability in results across studies. The findings demonstrate that key drivers such as environmental impact minimisation, constructability, and cost significantly influence the feasible region, while differences in weightings have a greater influence on outcomes than the choice of optimisation method. Consequently, the study underscores that achieving reliable, robust optimisation requires more balanced, data-driven weighting strategies and stronger integration of stakeholder preferences. Such improvements are critical for enhancing the realism and effectiveness of multi-objective optimisation in sustainable building design.

\section{Case Study: Small Modular Railway Station Buildings}

\begin{figure}[htbp]
\centering
\includegraphics[width = 1\textwidth]{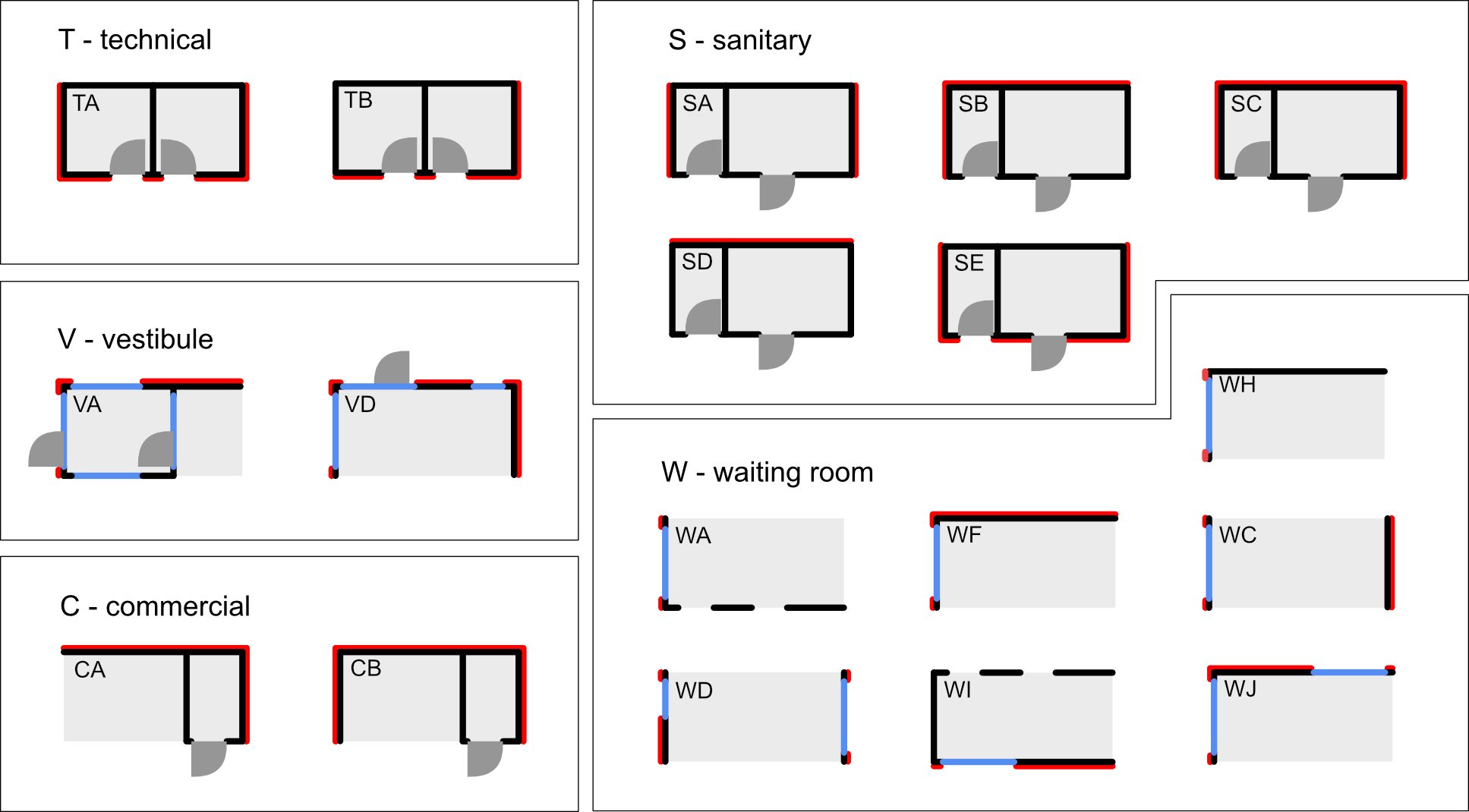}
\caption{Overview of module variants: black -- interior wall, red -- exterior cladding, blue -- glass facade, grey -- door.}
Source: own.\\
\label{fig:Fig 1}
\end{figure}

A modular railway station system was developed based on five primary groups of modules, differing in their functional purpose and building-envelope configurations. The system included technical (T), sanitary (S), vestibule (V), waiting room (W), and commercial modules (C). The modules used in the analysis are shown in Figure \ref{fig:Fig 1}.

\begin{figure}[h!t]
\centering
\includegraphics[width = 0.99\textwidth]{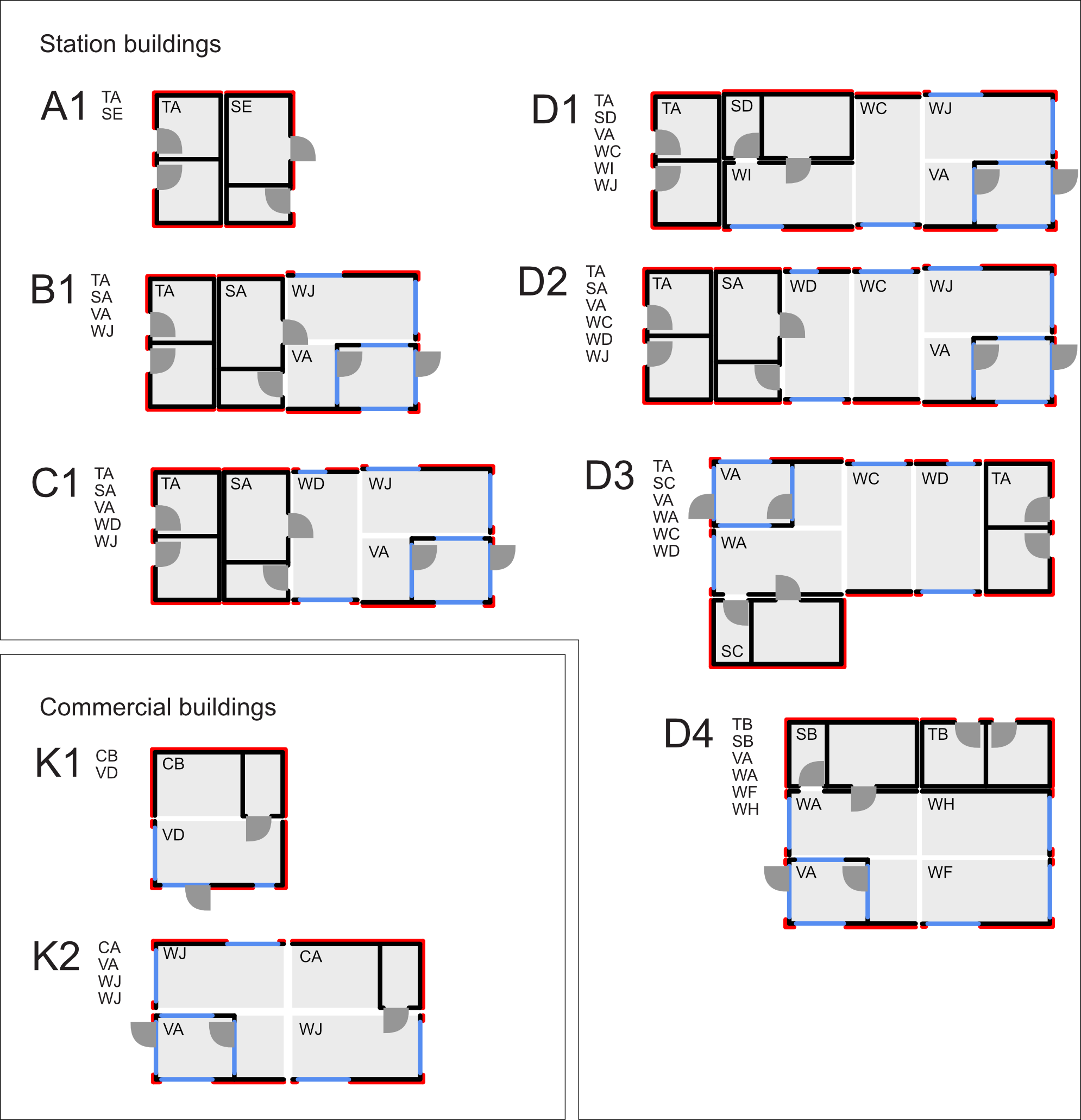}
\caption{Example module configurations forming railway station buildings of different sizes: T -- technical, S -- sanitary, V -- vestibule, W -- waiting room, C -- commercial.}
Source: own.\\
\label{fig:Fig 2}
\end{figure}

For each module type and variant, a detailed geometric analysis was performed, including calculations of external and internal wall areas, floor and roof areas, and the number and surface area of windows and doors. Based on the developed modules, several representative railway station configurations differing in size and functional programme, together with two commercial buildings (K1 and K2) composed of the same module types, were assembled (Figure \ref{fig:Fig 2}):
\begin{itemize}
\item[]    
A1 -- two-module building: T + S,
\item[]
B1 -- four-module building: T + S + V + W,
\item[]
C1 -- five-module building: T + S + V + 2W,
\item[]
D1--D4 -- six-module buildings: T + S + V + 3W,
\item[]
K1 -- two-module commercial building: C + V,
\item[]
K2 -- four-module commercial building: C + V + 2W.
\end{itemize}

In the next stage, detailed material quantity take-offs were prepared for each building variant. Material quantities were expressed in cubic metres (m\textsuperscript{3}), kilograms (kg), and square metres (m\textsuperscript{2}) for planar elements. Based on these quantities, an environmental assessment was conducted by calculating the total Global Warming Potential (GWP) for life-cycle stages A1--A3, covering raw material extraction, transport to manufacturing facilities, and material production processes.
In addition, the environmental impact associated with transporting materials from production plants to the assumed modular factory location in Oborniki Wielkopolskie (Poland) was assessed. An economic analysis was also conducted, including both material procurement costs and transportation costs.

\section{Data Preparation}

The data collection process was carried out in several stages and included the preparation of material, environmental, cost, and transportation datasets for the analysed modular railway station building variants. Based on the building model (BIM), detailed material quantity take-offs were prepared. Each material was assigned a place parameter indicating its location within the building (e.g., foundations, walls, or roof) and a type parameter defining the material category (e.g., concrete, ceramic tiles, tile adhesive). Additional information regarding density, volume, mass, and surface area was collected for each material.

In the next stage, Global Warming Potential (GWP) values were determined for all construction materials. For each material category, three to four specific products were selected based on the availability of environmental data within the OneClickLCA database. Preference was given to products manufactured in Poland. However, when suitable datasets were unavailable, products from other European countries were also considered. Where possible, the percentage of recycled content used in the manufacturing process was additionally recorded.

The environmental assessment adopted GWP values corresponding to life-cycle stages A1--A3, in accordance with Environmental Product Declarations (EPDs). These stages include raw material extraction, transport to the manufacturing plant, and the production process itself. Depending on the material type, GWP values were expressed in kg CO\textsubscript{2}e/m\textsuperscript{3}, kg CO\textsubscript{2}e/kg, or kg CO\textsubscript{2}e/m\textsuperscript{2}. While some materials were supported by full EPD documentation, others relied on generic datasets generated by the OneClickLCA platform. These generic values were derived from average environmental data for comparable products manufactured in other countries and adjusted to account for differences in national energy mixes relative to Poland.

Simultaneously, material cost data were collected. Several sources were investigated. Although cost information available in OneClickLCA was initially analysed, it was excluded from further consideration due to its approximate nature and the lack of differentiation between products from various manufacturers and countries. The most reliable sources were identified as the Sekocenbud cost database (4th quarter of 2024) and information obtained directly from manufacturers. However, such data were available only for commonly used construction materials on the Polish market.

For materials lacking reliable cost information, price estimates generated by large language models (ChatGPT and Grok) were used. Since the values generated by the language models were generally higher than those reported in the Sekocenbud database, a correction coefficient ($C_{cor}$) was introduced to calibrate the estimated unit prices. Unit costs expressed in PLN/m\textsuperscript{3}, PLN/kg, or PLN/m\textsuperscript{2} were multiplied by material quantities to obtain total material costs. The following notation was adopted:

\begin{itemize}
\item[]
$K_{s,i}$ - material cost based on Sekocenbud prices,
\item[]
$K_{c,i}$ - material cost based on ChatGPT-generated prices,
\item[]
$K_{g,i}$ - material cost based on Grok-generated prices,
\item[]
$i$ = 1 - reference material with available Sekocenbud pricing,
\item[]
$i$ = 2 or 3 - alternative materials characterized by lower GWP values than the reference material.
\end{itemize}

\noindent{}For each material category, the average LLM-generated cost was calculated as:

$$K_{av,i} = \frac{K_{c,i} + K_{g,i}}{2}.$$

\noindent{}Subsequently, correction coefficients were determined according to:

\begin{align*}
C_{cor,2} &= \frac{K_{av,2}}{K_{av,1}},\\
C_{cor,3} &= \frac{K_{av,3}}{K_{av,1}}.
\end{align*}

For each material group, the reference material cost was calculated by multiplying the available reference unit price by the corresponding material quantity. For alternatives without direct pricing data, the estimated total material cost was obtained by applying the dimensionless correction coefficient to the total cost of the reference material within the same group:

\begin{equation*}
K_i = C_{\mathrm{cor},i}K_{s,1},
\end{equation*}
where
\begin{itemize}
\item[]
$K_i$ -- the estimated total cost of alternative $i$,
\item[]
$K_{s,1}$ -- the total cost of the reference material for the analysed building.
\end{itemize}

The correction coefficient expresses the relative cost of the alternative relative to the reference material, derived from averaged LLM estimates. Thus, the LLM estimates determine relative cost differences, while the absolute cost level is anchored to the reference pricing data.

The next stage involved determining transportation costs. For each material, the location of the manufacturing plant, including city and country of origin, was identified. Transportation distances to the assumed modular factory location in Oborniki Wielkopolskie were obtained using Google Maps, excluding ferry crossings. Distances travelled through individual transit countries were also recorded.

Transportation costs were estimated using vehicle-specific tonne-kilometre rates adopted as calculation assumptions for the comparative analysis. For concrete mixer trucks with a capacity of 8 m\textsuperscript{3} and an assumed load of 19.2 t, a rate of 12.50 PLN/km was selected within the assumed range of 10--15 PLN/km, giving approximately 0.65 PLN/tkm. For trucks with a payload capacity of 12 t, the assumed rate of 4.90 PLN/km, within the range of 3--5 PLN/km, gave approximately 0.41 PLN/tkm. For trucks with a payload capacity of 24 t, the assumed rate of 6.50 PLN/km, within the range of 5--7 PLN/km, gave approximately 0.27 PLN/tkm. The latter two rates were selected from the upper part of the respective ranges as a conservative costing assumption. These values were used consistently across the analysed alternatives. They represent neither market-average prices nor offers from individual carriers.

Additional costs included a fixed charge of 200 PLN for the use of a truck-mounted crane (HDS) for deliveries performed by 20--24 tonne trucks. Minimum transport charges for partial loads were also considered, amounting to 100 PLN, 300 PLN, or 400 PLN depending on vehicle type. This adjustment was necessary because some materials were transported in quantities insufficient to fully utilise vehicle capacity.

Road tolls applicable in individual transit countries were also incorporated into the analysis. Country-specific toll rates were adopted based on publicly available data and included Germany, Czech Republic, Austria, France, Belgium, Hungary, Slovakia, Switzerland, Lithuania, Sweden, the United Kingdom, Romania, Italy, Poland, and other relevant countries.

A uniform 10\% reduction was adopted as a scenario assumption representing potential transport cost optimisation, for example through negotiated terms or recurring deliveries. Applied to the adopted base rates of 0.65, 0.41 and 0.27 PLN/tkm, this yielded calculation coefficients of 0.585, 0.369 and 0.243 PLN/tkm, respectively. These coefficients are reported, after rounding to two decimal places, as 0.59 PLN/tkm for concrete mixers, 0.37 PLN/tkm for 12 t trucks and 0.24 PLN/tkm for 24 t trucks. The calculations retained the unrounded coefficients. Road tolls were calculated separately.

The final stage of the methodology involved calculating transportation-related greenhouse gas emissions. Transport-related greenhouse gas emissions were calculated using a conventional diesel emission factor of 2.67 kg CO\textsubscript{2}e/litre, corresponding to the tank-to-wheel (TTW) factor specified in EN 16258:2012 and reproduced in the CLECAT guide, Calculating GHG Emissions for Freight Forwarding and Logistics Services (Table 5, p.26). This factor covers direct emissions from fuel combustion and excludes upstream emissions associated with fuel extraction, production and distribution (well-to-tank, WTT). It was applied consistently to all analysed transport alternatives. Fuel consumption values were adopted for each vehicle type under both loaded and unloaded operating conditions. For concrete mixer trucks, fuel consumption was assumed to be 50 l/100 km when fully loaded and 35 l/100 km when empty. For 7.5--12 tonne trucks, values of 24 l/100 km and 17 l/100 km were adopted, respectively, while 20--24 tonne trucks were assigned values of 32 l/100 km and 23 l/100 km.

For each material, the number of fully loaded delivery trips and the load fraction of the remaining partial-load trip were determined from the required material quantity and the vehicle capacity. Concrete deliveries were assessed using volume and an 8 m\textsuperscript{3} mixer capacity, whereas other deliveries were assessed using material mass and the assigned vehicle payload capacity. Fuel consumption for a partially loaded vehicle was estimated by linear interpolation between the assumed empty-vehicle and fully loaded fuel consumption values:

\begin{align*}
q_a &= q_0 + r(q_f - q_0),\\
E_f &= \frac{2.67q_f}{100},\\
E_a &= \frac{2.67q_a}{100},
\end{align*}

\noindent{}where:
\begin{itemize}
\item[]
$r$ -- the fraction of vehicle capacity used by the partial load,
\item[]
$q_0$ -- fuel consumption in litres per 100 km for an empty vehicle,
\item[]
$q_f$ -- fuel consumption in litres per 100 km for a fully loaded vehicle,
\item[]
$q_a$ -- fuel consumption in litres per 100 km for a partially loaded vehicle.
\end{itemize}

\noindent{}Transport emissions were in turn calculated as:

\begin{equation*}
E_t = D(A_fE_f + A_aE_a),
\end{equation*}

\noindent{}where:
\begin{itemize}
\item[]
$E_t$ -- transport emissions [kg CO\textsubscript{2}e],
\item[]
$D$ -- the one-way delivery distance [km],
\item[]
$A_f$ -- the number of fully loaded delivery trips,
\item[]
$A_a$ -- indicates whether an additional partial-load trip is required (A\textsubscript{a} = 1 if the remaining quantity requires a partially loaded trip, A\textsubscript{a} = 0 otherwise),
\item[]
$E_f$ -- emissions per vehicle-kilometre for fully loaded condition [kg CO\textsubscript{2}e/km]
\item[]
$E_a$ -- emissions per vehicle-kilometre for partially loaded conditions [kg CO\textsubscript{2}e/km].
\end{itemize}

Empty-vehicle fuel consumption was used as the lower bound of the interpolation. Emissions from empty return journeys were not added to the reported transport emissions.

\section{Optimal material selection}

\subsection{Model description}

Consider a fixed building module or complete building assembled from several modules. The module consists of a finite set of required parts
$$\mathcal{I} = \{1,\ldots,n\},$$
such as foundations, reinforcement, windows, insulation, roofing, and other construction components.

For each part $i \in \mathcal{I}$, let
$$\mathcal{T}_i = \{1,\ldots,q_i\}$$
be the set of admissible component types. The number $q_i$ may differ from one part to another, depending on the available suppliers and technical options.

For every part $i \in \mathcal{I}$ and every type $t \in \mathcal{T}_i$, two values are known: $c_{it} \geq 0$, the total cost of using type $t$ for part $i$, and $e_{it} \geq 0$, the total carbon footprint of using type $t$ for part $i$.

Both values are assumed to refer to the full quantity of part $i$ required for the entire module. Therefore, the total cost and the total carbon footprint of a complete module are obtained by summing the selected component-type values.

Define the binary decision variable
$$
x_{it} = 
\begin{cases}
1, & \text{if type } t \text{ is selected for part } i,\\
0, & \text{otherwise.}
\end{cases}
$$
Exactly one type must be selected for each part:
\begin{align*}
&\sum_{t \in \mathcal{T}_i} x_{it} = 1,&i \in \mathcal{I}.
\end{align*}
A complete selection can also be written as
\begin{align*}
&a = (t_1,\ldots,t_n),&t_i \in \mathcal{T}_i,
\end{align*}
where $t_i$ is the selected type for part $i$. Let $\mathcal{A}$ denote the set of feasible complete selections. If all combinations are allowed, then
$$|\mathcal{A}| = \prod_{i=1}^{n} q_i.$$
If technical, supplier, or compatibility restrictions exist, they should be included when defining $\mathcal{A}$. For a feasible selection $a \in \mathcal{A}$, the total cost is
$$
C(a) = \sum_{i \in \mathcal{I}}\sum_{t \in \mathcal{T}_i}c_{it} x_{it}.
$$
The total carbon footprint is
$$
E(a) = \sum_{i \in \mathcal{I}}\sum_{t \in \mathcal{T}_i}e_{it} x_{it}.
$$
The goal is to select component types so that both criteria are as small as possible:
$$
\min_{a \in \mathcal{A}}\bigl(C(a), E(a)\bigr).
$$
Since the least costly selection is not necessarily the lowest-carbon one, the problem is treated as a multi-criteria decision problem.

\subsection{TOPSIS method}

Let
$$
\mathcal{A} = \{a_1,\ldots,a_K\}
$$
be the set of feasible selections considered in the ranking. Since the normalisation used below depends on this set, the TOPSIS ranking is always understood relative to the set of alternatives actually enumerated (see Section~\ref{sec:numerical}). For each alternative $a_k$, compute $y_{k1} = C(a_k)$ and $y_{k2} = E(a_k)$, where $k=1, \ldots, K$. The TOPSIS decision matrix is
$$
Y =
\begin{bmatrix}
C(a_1) & E(a_1)\\
C(a_2) & E(a_2)\\
\vdots & \vdots\\
C(a_K) & E(a_K)
\end{bmatrix}.
$$
The two criteria may have different units, so they are normalized:
\begin{align*}
r_{kj} &= \frac{y_{kj}}{\sqrt{\sum_{\ell=1}^{K} y_{\ell j}^2}}, &k=1,\ldots,K,\quad j=1,2.
\end{align*}

\noindent{}Let
\begin{align*}
&w_1,w_2 \geq 0,&w_1+w_2=1,
\end{align*}
where $w_1$ is the weight of total cost and $w_2$ is the weight of carbon footprint. The weighted normalised values are
\begin{align*}
&v_{kj} = w_j r_{kj}, &k=1,\ldots,K, &j=1,2.
\end{align*}
If both criteria are equally important, then $w_1 = w_2 = \frac{1}{2}$. Both criteria are minimised. Therefore, the positive ideal solution is
\begin{align*}
&A^+ = (v_1^+,v_2^+), &v_j^+ = \min_{k=1,\ldots,K} v_{kj},
\end{align*}
while the negative ideal solution is
\begin{align*}
&A^- = (v_1^-,v_2^-), &v_j^- = \max_{k=1,\ldots,K} v_{kj}.
\end{align*}
Thus, $A^+$ represents the best attainable values among the considered selections, while $A^-$ represents the worst attainable values. For each selection $a_k$, compute the distances from the positive and negative ideals:
\begin{align*}
D_k^+ &= \sqrt{\sum_{j=1}^{2}(v_{kj}-v_j^+)^2},\\
D_k^- &= \sqrt{\sum_{j=1}^{2}(v_{kj}-v_j^-)^2}.
\end{align*}

\noindent{}The TOPSIS closeness coefficient is
\begin{align*}
S_k = \frac{D_k^-}{D_k^+ + D_k^-}, &k=1,\ldots,K.
\end{align*}

A larger value of $S_k$ indicates that the selection is closer to the ideal low-cost, low-carbon solution and farther from the high-cost, high-carbon solution. The preferred selection is the alternative with the largest closeness coefficient:
$$
a^\star = \arg\max_{a_k \in \mathcal{A}} S_k.
$$
The corresponding total cost and carbon footprint are
\begin{align*}
C^\star = C(a^\star), &E^\star = E(a^\star).
\end{align*}

The TOPSIS workflow can be written as Algorithm~\ref{alg:topsis-fixed-module}.

\begin{algorithm}[h]
\caption{TOPSIS for component-type selection}
\label{alg:topsis-fixed-module}
\small
\KwInput{Parts $\mathcal{I}=\{1,\ldots,n\}$; type sets $\mathcal{T}_i$; values $c_{it}$ and $e_{it}$; feasible alternatives $\mathcal{A}=\{a_1,\ldots,a_K\}$; weights $w_1,w_2$, where $w_1+w_2=1$}
\KwOutput{Ranking of feasible alternatives and preferred selection $a^\star$}

\ForEach{alternative $a_k \in \mathcal{A}$}{
    $y_{k1} \leftarrow C(a_k)=\sum_{i \in \mathcal{I}}\sum_{t \in \mathcal{T}_i} c_{it}x_{it}$\;
    $y_{k2} \leftarrow E(a_k)=\sum_{i \in \mathcal{I}}\sum_{t \in \mathcal{T}_i} e_{it}x_{it}$\;
}
Build the decision matrix $Y=[y_{kj}]_{K \times 2}$\;

\For{$j \leftarrow 1$ \KwTo $2$}{
    $d_j \leftarrow \sqrt{\sum_{\ell=1}^{K} y_{\ell j}^2}$\;
    \For{$k \leftarrow 1$ \KwTo $K$}{
        $r_{kj} \leftarrow y_{kj}/d_j$\;
        $v_{kj} \leftarrow w_j r_{kj}$\;
    }
    $v_j^+ \leftarrow \min_{k=1,\ldots,K} v_{kj}$ and
    $v_j^- \leftarrow \max_{k=1,\ldots,K} v_{kj}$ \tcp*[r]{both criteria are minimized}
}

\For{$k \leftarrow 1$ \KwTo $K$}{
    $D_k^+ \leftarrow \sqrt{\sum_{j=1}^{2}(v_{kj}-v_j^+)^2}$\;
    $D_k^- \leftarrow \sqrt{\sum_{j=1}^{2}(v_{kj}-v_j^-)^2}$\;
    $S_k \leftarrow D_k^-/(D_k^+ + D_k^-)$\;
}

Rank alternatives $a_k$ by decreasing $S_k$\;
$a^\star \leftarrow \arg\max_{a_k \in \mathcal{A}} S_k$\;
\Return{$a^\star$ and the ranked list of alternatives}\;
\end{algorithm}

\subsection{A set of alternatives}

Instead of selecting only one pair of weights $(w_1,w_2)$, the TOPSIS procedure may be repeated for a sequence of weight pairs. This makes it possible to obtain a set of preferred selections corresponding to different relative priorities between total cost and carbon footprint. Let $\varepsilon > 0$ be a weight step. For simplicity, assume that $\frac{1}{\varepsilon} = H$ is an integer. The sequence of weight pairs is then defined as
\begin{align*}
&(w_1^{(h)},w_2^{(h)}) = (h\varepsilon, 1-h\varepsilon), &h=0,1,\ldots,H.
\end{align*}

For each value $h=0,1,\ldots,H$, TOPSIS is performed using the corresponding weights $(w_1^{(h)},w_2^{(h)})$. This gives a weighted normalised matrix
\begin{align*}
&v_{kj}^{(h)} = w_j^{(h)} r_{kj}, &k=1,\ldots,K,\quad j=1,2.
\end{align*}
The positive and negative ideal values $v_j^{+,(h)}$, $v_j^{-,(h)}$, the corresponding distances $D_k^{+,(h)}$ and $D_k^{-,(h)}$, the closeness coefficient $S_k^{(h)}$, and the preferred selection $a^{\star,(h)}$ are then computed for that particular weight pair.

Repeating the procedure for all $h=0,1,\ldots,H$ gives the ordered list of TOPSIS selections
$$
\mathcal{L}_{\varepsilon} = \left(a^{\star,(0)}, a^{\star,(1)}, \ldots, a^{\star,(H)}\right).
$$
If only distinct selections are needed, one may define the corresponding set
$$
\mathcal{S}_{\varepsilon} = \left\{a^{\star,(h)} : h=0,1,\ldots,H\right\}.
$$
Different weight pairs may lead to the same selected alternative, so in general
$$
|\mathcal{S}_{\varepsilon}| \leq H+1.
$$
For instance, if $\varepsilon = 0.1$, then $H=10$, and TOPSIS is performed for the $11$ pairs
$$
(0,1), (0.1,0.9),(0.2,0.8), \ldots,  (0.9,0.1), (1,0).
$$
This produces eleven TOPSIS results in the ordered list $\mathcal{L}_{0.1}$, although the number of distinct component-type selections in $\mathcal{S}_{0.1}$ may be smaller if some alternatives are optimal for more than one weight pair.

A more systematic approach may also be used without TOPSIS. Instead of assigning weights and computing a closeness coefficient, one may construct a sequence of solutions moving from the minimum-cost selection toward the minimum-emission selection. The goal is to obtain consecutive nondominated solutions representing the trade-off between total cost and carbon footprint. Before the sequence is generated, component-wise dominated types may be removed from consideration. Let $\widehat{\mathcal{T}}_i \subseteq \mathcal{T}_i$ denote the set of retained types for group $i$. This filtering never renumbers the types, i.e.,  every retained type keeps its original label $t$ from the source data. For example, if $\mathcal{T}_i=\{1,2,3,4\}$ and types $1$ and $3$ are removed, then $\widehat{\mathcal{T}}_i=\{2,4\}$, not $\{1,2\}$. Consequently, all vectors listed below remain directly comparable with the original data.

Because both criteria are additive over parts, the minimum-cost selection can be obtained greedily by choosing, for every part $i$, the type with the smallest cost:
$$
t_i^C \in \arg\min_{t \in \widehat{\mathcal{T}}_i} c_{it}.
$$
The corresponding selection is
$$
a^C = (t_1^C,\ldots,t_n^C).
$$
If several types have the same minimum cost, the one with the smallest carbon footprint may be selected as a tie-breaker.

Analogously, the minimum-emission selection is obtained by choosing, for every part $i$, the type with the smallest carbon footprint:
$$
t_i^E \in \arg\min_{t \in \widehat{\mathcal{T}}_i} e_{it}.
$$
The corresponding selection is
$$
a^E = (t_1^E,\ldots,t_n^E).
$$
If several types have the same minimum carbon footprint, the one with the smallest cost may be selected as a tie-breaker. Let
$$
a^{(r)} = (t_1^{(r)},\ldots,t_n^{(r)})
$$
be the current selection at iteration $r$. A one-group change means replacing the currently selected type in exactly one component group $i$. If type $u \in \widehat{\mathcal{T}}_i$ is chosen instead of $t_i^{(r)}$, the candidate selection is denoted by
$$
a^{(r)}(i \leftarrow u) = (t_1^{(r)},\ldots,t_{i-1}^{(r)},u,t_{i+1}^{(r)},\ldots,t_n^{(r)}).
$$

The change in total cost is
$$
\Delta C_i(u \mid a^{(r)})
= C\bigl(a^{(r)}(i \leftarrow u)\bigr) - C\bigl(a^{(r)}\bigr),
$$
and the reduction in carbon footprint is
$$
\Delta E_i(u \mid a^{(r)})
= E\bigl(a^{(r)}\bigr) - E\bigl(a^{(r)}(i \leftarrow u)\bigr).
$$
The useful changes are those for which
\begin{align*}
&\Delta C_i(u \mid a^{(r)}) > 0,\\
&\Delta E_i(u \mid a^{(r)}) > 0.
\end{align*}
Such a change increases the cost but decreases the carbon footprint.

For every useful one-group change, compute the average emission reduction per additional unit of cost:
$$
\rho_i(u \mid a^{(r)}) =
\frac{\Delta E_i(u \mid a^{(r)})}{\Delta C_i(u \mid a^{(r)})}.
$$
At iteration $r$, choose the change with the largest value of this ratio:
$$
(i_r,u_r) \in \arg\max_{(i,u)}\rho_i(u \mid a^{(r)}),
$$
where the maximum is taken over all useful one-group changes. The next selection is
$$
a^{(r+1)} = a^{(r)}(i_r \leftarrow u_r).
$$
This means that at each step, the procedure selects the component-type replacement that gives the largest carbon-footprint reduction per one additional unit of cost.

If a change satisfies
\begin{align*}
&\Delta C_i(u \mid a^{(r)}) = 0,\\
&\Delta E_i(u \mid a^{(r)}) > 0,
\end{align*}
then it should be performed before positive-cost changes, because it reduces emissions without increasing cost. After the component-wise filtering described above, such zero-cost changes cannot occur, because a type with the same cost and a lower carbon footprint would dominate the current one. We preserve the rule for the sake of generality. The procedure starts from $a^{(0)} = a^C$. Then the Pareto-sequence step is repeated until the minimum-emission selection $a^E$ is reached or until no useful one-group change exists. The generated sequence is $\mathcal{Q} = \left(a^{(0)}, a^{(1)}, \ldots, a^{(R)}\right)$, where, for the data considered here, $a^{(0)} = a^C$ and $a^{(R)} = a^E$. The sequential non-TOPSIS procedure can be written as Algorithm~\ref{alg:pareto-set}.

\begin{algorithm}[ht]
\caption{Pareto sequence based on one-group changes}
\label{alg:pareto-set}
\footnotesize
\KwInput{Parts $\mathcal{I}$; retained type sets $\widehat{\mathcal{T}}_i$ with original labels; values $c_{it}$ and $e_{it}$}
\KwOutput{Generated Pareto sequence $\mathcal{Q}$}

\ForEach{part $i \in \mathcal{I}$}{
    $t_i^C \in \arg\min_{t \in \widehat{\mathcal{T}}_i} c_{it}$ \tcp*[r]{tie-break by smaller $e_{it}$}
    $t_i^E \in \arg\min_{t \in \widehat{\mathcal{T}}_i} e_{it}$ \tcp*[r]{tie-break by smaller $c_{it}$}
}

$a^C \leftarrow (t_1^C,\ldots,t_n^C)$, $a^E \leftarrow (t_1^E,\ldots,t_n^E)$\;
$a^{(0)} \leftarrow a^C$, $\mathcal{Q} \leftarrow (a^{(0)})$, $r \leftarrow 0$\;

\While{$a^{(r)} \neq a^E$}{
    $\mathcal{U}_0 \leftarrow \emptyset$, $\mathcal{U}_+ \leftarrow \emptyset$\;

    \ForEach{part $i \in \mathcal{I}$}{
        \ForEach{type $u \in \widehat{\mathcal{T}}_i \setminus \{t_i^{(r)}\}$}{
            $b \leftarrow a^{(r)}(i \leftarrow u)$\;
            $\Delta C \leftarrow C(b)-C(a^{(r)})$, $\Delta E \leftarrow E(a^{(r)})-E(b)$\;
            \If{$\Delta E>0$ and $\Delta C=0$}{
                Add $b$ to $\mathcal{U}_0$\;
            }
            \ElseIf{$\Delta E>0$ and $\Delta C>0$}{
                $\rho(b) \leftarrow \Delta E/\Delta C$\;
                Add $b$ to $\mathcal{U}_+$\;
            }
        }
    }

    \If{$\mathcal{U}_0 \neq \emptyset$}{
        $b^\star \leftarrow \arg\max_{b \in \mathcal{U}_0} \left(E(a^{(r)})-E(b)\right)$\;
    }
    \ElseIf{$\mathcal{U}_+ \neq \emptyset$}{
        $b^\star \leftarrow \arg\max_{b \in \mathcal{U}_+} \rho(b)$\;
    }
    \Else{
        \textbf{break}\;
    }

    $a^{(r+1)} \leftarrow b^\star$\;
    Append $a^{(r+1)}$ to $\mathcal{Q}$\;
    $r \leftarrow r+1$\;
}

\Return{$\mathcal{Q}$}\;
\end{algorithm}

\subsection*{Convexified Interpretation}

For the original discrete problem, the resulting ordered list is an ordered Pareto set (i.e., Pareto sequence), but it need not enumerate every nondominated discrete alternative. In particular, a nondominated discrete point can lie above the line segment joining two other feasible points. Such a point can be omitted by the one-group ratio rule, even though it remains nondominated as long as only integral material selections are allowed.

Let
$$
\mathcal{Y} = \left\{\bigl(C(a),E(a)\bigr): a\in\prod_{i\in\mathcal{I}}\widehat{\mathcal{T}}_i\right\}
$$
be the objective vectors of all discrete selections. Suppose now that fractional selections are admissible in the objective space, so that the attainable objective set is the convex hull $\operatorname{conv}(\mathcal{Y})$. A point on a segment between two discrete solutions then represents their convex combination.

Under the assumptions used here, namely independent material groups and additive cost and emission, the Pareto sequence $\mathcal{Q}$ gives the complete Pareto frontier of $\operatorname{conv}(\mathcal{Y})$. To see this, consider the lower convex chain of each group in the $(C, E)$ plane. Its successive emission-reduction-per-cost ratios are non-increasing. A one-group jump spanning several such edges has a ratio equal to a weighted average of their ratios, and therefore cannot exceed the first untraversed edge. Consequently, the rule maximising $\Delta E/\Delta C$ merges all local lower-convex-chain edges in decreasing ratio order. The resulting polyline is precisely the lower convex frontier of $\operatorname{conv}(\mathcal{Y})$.

If several consecutive edges have the same slope, the procedure may list additional collinear intermediate points. They lie on the same Pareto segment and do not change the frontier. Thus, for fractional selections the segments joining consecutive points of $\mathcal{Q}$ form the entire Pareto set. For the purely discrete problem, the same path remains useful but can omit nondominated discrete alternatives.

\subsection{Numerical example -- building type A1}\label{sec:numerical}

\subsubsection{Introduction}

We start with solving the numerical instance for building type A1. The data are not repeated here. Each component group is indexed by $i$, and each product variant inside a group is indexed by $t$, as described in the previous model. For the detailed $A1$ instance, the full combinatorial space contains
$$
|\mathcal{A}| = 3^{22}\cdot 4 = \text{125,524,238,436}
$$
feasible combinations. Component-wise dominated variants were removed before the search. Since both criteria are minimised, a selection that uses a dominated variant is itself dominated and can never be an efficient solution found by either method (TOPSIS or Pareto-sequence). Note, however, that the TOPSIS normalisation and the anti-ideal solution $A^-$ are computed over the enumerated set, so the closeness coefficients $S_k$ and, for intermediate weights, also the selected alternative depend on this reduction. All TOPSIS results listed below refer to the reduced set. This leaves 1728 relevant combinations. Ties with identical cost and emission are preserved. In the A1 instance, groups $i=4$ (window profile), $i=5$ (window glass) and $i=17$ (interior door) do not occur in the building. All their variants have zero cost and zero emission. The index $t$ selected for these groups is therefore an arbitrary tie and does not represent a product choice. In particular, the 1728 retained combinations consist of 64 distinct selections, each of which appears 27 times with different labels of these three groups. The same occurs for group $i=18$ (exterior door) in buildings K1 and K2, where every combination occurs three times. Every preserved type keeps its original index $t$ from the data spreadsheet. Among these discrete combinations, 20 objective vectors are nondominated. The Pareto frontier of their convex hull has 7 vertices, and these vertices are exactly the points generated by the Pareto sequence procedure below.

\subsubsection{TOPSIS solution -- equal importance of criteria}

The TOPSIS computation was carried out with equal criterion weights,
$$
w_1=w_2=0.5,
$$
where $w_1$ is the weight of total cost and $w_2$ is the weight of carbon emission. The selected alternative is
$$
a^{\star}_{\mathrm{TOPSIS}} =
(3,4,1,1,1,1,2,1,3,3,1,1,1,1,1,3,1,1,1,3,3,2,1).
$$
For this selection, $C(a^{\star}_{\mathrm{TOPSIS}})=120381.180$, $E(a^{\star}_{\mathrm{TOPSIS}})=18621.460$ and $S=0.649461$.

\paragraph{Product interpretation of the TOPSIS solution.}
{\footnotesize
\begin{itemize}
\item $i=1$, \texttt{foundation/concrete}, $t_{1}=3$ (ID 3): Ready-mix concrete, generic, C25/30 (3625/4350 PSI), XC1, 50\% fly ash, CEM IV/B pozzolanic cement.
\item $i=2$, \texttt{foundation/steel rebar}, $t_{2}=4$ (ID 7): Reinforcing steel (rebar), generic, 97\% recycled content (typical).
\item $i=3$, \texttt{foundation/membrane}, $t_{3}=1$ (ID 8): EPDM waterproofing membrane, 1.5 mm, 1.95 kg/m2.
\item $i=4$, \texttt{wall/window profile}, $t_{4}=1$ (ID 11): Aluminium profiles, 2700 kg/m3, EN AW 1050, 6060, 6063, 6005A (Aliplast Extrusion Sp. z o.o.).
\item $i=5$, \texttt{wall/window glass}, $t_{5}=1$ (ID 14): Fire-resistant safety glass with three laminated glass panes and two intumescent interlayers, 17.3 mm, 41.6 kg/m2.
\item $i=6$, \texttt{wall/stiffening}, $t_{6}=1$ (ID 17): Fibre cement boards, 1300 kg/m3 (81.16 lbs/ft3).
\item $i=7$, \texttt{wall, roof/thermal insulation p}, $t_{7}=2$ (ID 21): EPS insulation panels, L= 0.038 m2K/W, R=1 m2K/W, 38 mm (10-200 mm range), 600 x 1200mm, 0.5928 kg/m2, 15.6 kg/m3, compressive strength class 80 kN/m2, Lambda=0.038 W/(m.K), EPS 80 Graphite (EPS Sverige (BeWi Insulation, Jackon, Sundolitt).
\item $i=8$, \texttt{wall/facing e}, $t_{8}=1$ (ID 23): Fibre cement boards, 1300 kg/m3 (81.16 lbs/ft3).
\item $i=9$, \texttt{wall/permable membrane}, $t_{9}=3$ (ID 28): 3-layer vapour permeable underlay, 1.1-1.5 x 50 m, Divoroll (Braas Monier).
\item $i=10$, \texttt{wall, roof/thermal insulation w}, $t_{10}=3$ (ID 31): Stone wool insulation compact panels, unfaced, L = 0.035 W/mK, R = 2.86 m2K/W, 100 mm, 10 kg/m2, 100 kg/m3, Fire resistance class = A1, Lambda=0.035 W/(m.K), FKD-S, FKD-S C1, FKD-S C2, SmartWall S C1, SmartWall S C2, FKD-S Thermal, FKD MAX C1, FKD MAX C2, SmartWall FireGuard, FKD-T FB C2 (Knauf (2021)).
\item $i=11$, \texttt{wall, floor, roof/vapour barrier}, $t_{11}=1$ (ID 32): Polyethene vapour barrier membrane, 0.15 mm, 0.14 kg/m2.
\item $i=12$, \texttt{wall, roof/facing i}, $t_{12}=1$ (ID 35): Gypsum plaster board, regular, 10\% recycled gypsum, 6.5-25 mm (0.25-0.98 in), 10.725 kg/m2 (2.20 lbs/ft2) (for 12.5 mm/0.49 in), 858 kg/m3 (53.6 lbs/ft3.
\item $i=13$, \texttt{wall/underlay}, $t_{13}=1$ (ID 38): Gypsum plasterboard, water resistant, 12.5 mm, 8 kg/m2, 640 kg/m3, NIDA WODA 12,5 typ H2 (Siniat).
\item $i=14$, \texttt{wall, floor/tile glue}, $t_{14}=1$ (ID 41): Dry mortar, adhesive for facades and tiles, Consumption 3-4kg/m2, Weber Therm Elastik (Weber).
\item $i=15$, \texttt{wall, floor/tiles}, $t_{15}=1$ (ID 44): Ceramic glazed tile, 20 kg/m2.
\item $i=16$, \texttt{structure/steel profiles}, $t_{16}=3$ (ID 49): Structural steel profiles, generic, 100\% recycled content, I, H, U, L, and T sections, S235, S275 and S355.
\item $i=17$, \texttt{wall/door i}, $t_{17}=1$ (ID 50): Wooden and engineered wood interior doors, 1,23 m x 2,18 m, 72.8 kg/m2.
\item $i=18$, \texttt{wall/door e}, $t_{18}=1$ (ID 53): Steel door with mineral wool core, fire-rated, 900x2000 mm, 74 kg/unit.
\item $i=19$, \texttt{floor/thermal insulation p}, $t_{19}=1$ (ID 56): EPS insulation panels, graphite, L= 0.033 W/mK, R= 3.03 m2K/W, 100 mm, 3 kg/m2, 30 kg/m3, compressive strength 220 kPa, 100\% virgin polystyrene, Lambda=0.033 W/(m.K).
\item $i=20$, \texttt{floor/concrete}, $t_{20}=3$ (ID 61): Ready-mix concrete, generic, C25/30 (3625/4350 PSI), XC1, 50\% fly ash, CEM IV/B pozzolanic cement,.
\item $i=21$, \texttt{floor, roof/thermal insulation f}, $t_{21}=3$ (ID 64): PIR insulation board, aluminium facing, L=0.022 W/mK, R=3.6 m2K/W, 80 mm, 2.8 kg/m2, Lambda=0.022 W/(m.K) (UNILIN).
\item $i=22$, \texttt{floor, roof/moisture insulation}, $t_{22}=2$ (ID 66): Bitumen sheets, G 200 S4; 5,0 kg/m2, EN15804+A1, ref. year 2018.
\item $i=23$, \texttt{floor, roof/steel sheet}, $t_{23}=1$ (ID 68): Hot-dip galvanised steel sheets, Steel thickness range: 0.4-3.0 mm (0.015-0.12 in), zinc coating: 20 $\mu$m (787.4 $\mu$in) (0.28kg/m2 / 0.057 lbs/ft2 sheet steel), 0\% recycled content.
\end{itemize}
}

\paragraph{TOPSIS solutions for a weight sweep with $\varepsilon=0.1$.}
The TOPSIS procedure was also repeated for the sequence of weights $(w_1,w_2)=(1,0),(0.9,0.1),\ldots,(0.1,0.9),(0,1)$, where $w_1$ is the cost weight and $w_2$ is the emission weight. The resulting solutions are:

{\footnotesize
\begin{description}
    \item[$(1,0)$] $C=113646.530$, $E=20729.960$, $S=1.000000$,
    $a^\star=$\texttt{(3, 4, 1, 1, 1, 1, 3, 1, 3, 3, 1, 1, 1, 1, 1, 3, 1, 2, 1, 3, 1, 2, 1)}.

    \item[$(0.9,0.1)$] $C=113746.460$, $E=20128.740$, $S=0.887755$,
    $a^\star=$\texttt{(3, 4, 1, 1, 1, 1, 3, 1, 3, 3, 1, 1, 1, 1, 1, 3, 1, 2, 1, 3, 3, 2, 1)}.

    \item[$(0.8,0.2)$] $C=113746.460$, $E=20128.740$, $S=0.779439$,
    $a^\star=$\texttt{(3, 4, 1, 1, 1, 1, 3, 1, 3, 3, 1, 1, 1, 1, 1, 3, 1, 2, 1, 3, 3, 2, 1)}.

    \item[$(0.7,0.3)$] $C=114243.920$, $E=20027.360$, $S=0.675699$,
    $a^\star=$\texttt{(3, 4, 1, 1, 1, 1, 3, 1, 3, 2, 1, 1, 1, 1, 1, 3, 1, 2, 1, 3, 3, 2, 1)}.

    \item[$(0.6,0.4)$] $C=116481.820$, $E=19476.170$, $S=0.607585$,
    $a^\star=$\texttt{(3, 4, 1, 1, 1, 1, 3, 1, 3, 2, 1, 1, 1, 1, 1, 3, 1, 1, 1, 3, 3, 2, 1)}.

    \item[$(0.5,0.5)$] $C=120381.180$, $E=18621.460$, $S=0.649461$,
    $a^\star=$\texttt{(3, 4, 1, 1, 1, 1, 2, 1, 3, 3, 1, 1, 1, 1, 1, 3, 1, 1, 1, 3, 3, 2, 1)}.

    \item[$(0.4,0.6)$] $C=120878.640$, $E=18520.080$, $S=0.722461$,
    $a^\star=$\texttt{(3, 4, 1, 1, 1, 1, 2, 1, 3, 2, 1, 1, 1, 1, 1, 3, 1, 1, 1, 3, 3, 2, 1)}.

    \item[$(0.3,0.7)$] $C=122258.030$, $E=18368.570$, $S=0.788675$,
    $a^\star=$\texttt{(3, 4, 1, 1, 1, 1, 2, 1, 3, 2, 1, 1, 1, 1, 1, 3, 1, 1, 1, 3, 3, 2, 2)}.

    \item[$(0.2,0.8)$] $C=122882.910$, $E=18326.910$, $S=0.855207$,
    $a^\star=$\texttt{(3, 4, 1, 1, 1, 1, 2, 1, 3, 2, 1, 2, 1, 1, 1, 3, 1, 1, 1, 3, 3, 2, 1)}.

    \item[$(0.1,0.9)$] $C=124262.300$, $E=18175.400$, $S=0.929775$,
    $a^\star=$\texttt{(3, 4, 1, 1, 1, 1, 2, 1, 3, 2, 1, 2, 1, 1, 1, 3, 1, 1, 1, 3, 3, 2, 2)}.

    \item[$(0,1)$] $C=124262.300$, $E=18175.400$, $S=1.000000$,
    $a^\star=$\texttt{(3, 4, 1, 1, 1, 1, 2, 1, 3, 2, 1, 2, 1, 1, 1, 3, 1, 1, 1, 3, 3, 2, 2)}.

\end{description}
}

The selections obtained for $(0.7,0.3)$, $(0.6,0.4)$ and $(0.2,0.8)$ are nondominated discrete alternatives that do not belong to the Pareto sequence generated below. This illustrates the remark made above that the sequence may omit nondominated discrete points, whereas TOPSIS may select them.

\subsubsection{Pareto-sequence solution}

The Pareto sequence starts from the minimum-cost solution and then applies one-group changes that provide the largest emission reduction per additional cost unit. The starting point is $a^{(0)}=(3,4,1,1,1,1,3,1,3,3,1,1,1,1,1,3,1,2,1,3,1,2,1)$, where $C(a^{(0)})=113646.530$ and $E(a^{(0)})=20729.960$. The successive changes are summarised below.

\begin{center}
\small
\begin{tabular}{c l r r r}
\hline
Step & Change & $\Delta C$ & $\Delta E$ & $\rho=\Delta E/\Delta C$ \\
\hline
1 & $i=21:\ t_{21}=1\to 3$ & 99.930 & 601.220 & 6.016 \\
2 & $i=18:\ t_{18}=2\to 1$ & 2237.900 & 551.190 & 0.246 \\
3 & $i=7:\ t_{7}=3\to 2$ & 4396.820 & 956.090 & 0.217 \\
4 & $i=10:\ t_{10}=3\to 2$ & 497.460 & 101.380 & 0.204 \\
5 & $i=23:\ t_{23}=1\to 2$ & 1379.390 & 151.510 & 0.110 \\
6 & $i=12:\ t_{12}=1\to 2$ & 2004.270 & 193.170 & 0.096 \\
\hline
\end{tabular}
\end{center}

\paragraph{Complete list of Pareto-sequence solutions.}
The table above describes the changes. The full solution sequence generated by these changes is:

{\footnotesize
\begin{description}
\item[Step 0] $C=113646.530$, $E=20729.960$,
$a^{(0)}=$\texttt{(3, 4, 1, 1, 1, 1, 3, 1, 3, 3, 1, 1, 1, 1, 1, 3, 1, 2, 1, 3, 1, 2, 1)}.
\item[Step 1] $C=113746.460$, $E=20128.740$,
$a^{(1)}=$\texttt{(3, 4, 1, 1, 1, 1, 3, 1, 3, 3, 1, 1, 1, 1, 1, 3, 1, 2, 1, 3, 3, 2, 1)}.
\item[Step 2] $C=115984.360$, $E=19577.550$,
$a^{(2)}=$\texttt{(3, 4, 1, 1, 1, 1, 3, 1, 3, 3, 1, 1, 1, 1, 1, 3, 1, 1, 1, 3, 3, 2, 1)}.
\item[Step 3] $C=120381.180$, $E=18621.460$,
$a^{(3)}=$\texttt{(3, 4, 1, 1, 1, 1, 2, 1, 3, 3, 1, 1, 1, 1, 1, 3, 1, 1, 1, 3, 3, 2, 1)}.
\item[Step 4] $C=120878.640$, $E=18520.080$,
$a^{(4)}=$\texttt{(3, 4, 1, 1, 1, 1, 2, 1, 3, 2, 1, 1, 1, 1, 1, 3, 1, 1, 1, 3, 3, 2, 1)}.
\item[Step 5] $C=122258.030$, $E=18368.570$,
$a^{(5)}=$\texttt{(3, 4, 1, 1, 1, 1, 2, 1, 3, 2, 1, 1, 1, 1, 1, 3, 1, 1, 1, 3, 3, 2, 2)}.
\item[Step 6] $C=124262.300$, $E=18175.400$,
$a^{(6)}=$\texttt{(3, 4, 1, 1, 1, 1, 2, 1, 3, 2, 1, 2, 1, 1, 1, 3, 1, 1, 1, 3, 3, 2, 2)}.
\end{description}
}

In this instance, the equal-weight TOPSIS solution is Step 3 on the Pareto sequence. The sequence continues toward lower emissions by accepting additional cost.
\paragraph{Product interpretation of the final Pareto-sequence solution.}
{\footnotesize
\begin{itemize}
    \item $i=1$, \texttt{foundation/concrete}, $t_{1}=3$ (ID 3): Ready-mix concrete, generic, C25/30 (3625/4350 PSI), XC1, 50\% fly ash, CEM IV/B pozzolanic cement.
    \item $i=2$, \texttt{foundation/steel rebar}, $t_{2}=4$ (ID 7): Reinforcing steel (rebar), generic, 97\% recycled content (typical).
    \item $i=3$, \texttt{foundation/membrane}, $t_{3}=1$ (ID 8): EPDM waterproofing membrane, 1.5 mm, 1.95 kg/m2.
    \item $i=4$, \texttt{wall/window profile}, $t_{4}=1$ (ID 11): Aluminium profiles, 2700 kg/m3, EN AW 1050, 6060, 6063, 6005A (Aliplast Extrusion Sp. z o.o.).
    \item $i=5$, \texttt{wall/window glass}, $t_{5}=1$ (ID 14): Fire-resistant safety glass with three laminated glass panes and two intumescent interlayers, 17.3 mm, 41.6 kg/m2.
    \item $i=6$, \texttt{wall/stiffening}, $t_{6}=1$ (ID 17): Fibre cement boards, 1300 kg/m3 (81.16 lbs/ft3).
    \item $i=7$, \texttt{wall, roof/thermal insulation p}, $t_{7}=2$ (ID 21): EPS insulation panels, L= 0.038 m2K/W, R=1 m2K/W, 38 mm (10-200 mm range), 600 x 1200mm, 0.5928 kg/m2, 15.6 kg/m3, compressive strength class 80 kN/m2, Lambda=0.038 W/(m.K), EPS 80 Graphite (EPS Sverige (BeWi Insulation, Jackon, Sundolitt).
    \item $i=8$, \texttt{wall/facing e}, $t_{8}=1$ (ID 23): Fibre cement boards, 1300 kg/m3 (81.16 lbs/ft3).
    \item $i=9$, \texttt{wall/permable membrane}, $t_{9}=3$ (ID 28): 3-layer vapour permeable underlay, 1.1-1.5 x 50 m, Divoroll (Braas Monier).
    \item $i=10$, \texttt{wall, roof/thermal insulation w}, $t_{10}=2$ (ID 30): 
Stone wool insulation compact panels, uncoated, L = 0.034 W/mK, R = 2.90 m2K/W, 100 mm, 9.0 kg/m2, 90 kg/m3, Fire resistance class = A1, Lambda=0.035 W/(m.K), FKD N Thermal (Knauf) .
    \item $i=11$, \texttt{wall, floor, roof/vapour barrier}, $t_{11}=1$ (ID 32): Polyethene vapour barrier membrane, 0.15 mm, 0.14 kg/m2.
    \item $i=12$, \texttt{wall, roof/facing i}, $t_{12}=2$ (ID 36): Impregnated gypsum plasterboard, 12.5 mm, 7.27 kg/m2, Siniat Smart H (Etex Poland Sp. z o.o.).
    \item $i=13$, \texttt{wall/underlay}, $t_{13}=1$ (ID 38): Gypsum plasterboard, water resistant, 12.5 mm, 8 kg/m2, 640 kg/m3, NIDA WODA 12,5 typ H2 (Siniat).
    \item $i=14$, \texttt{wall, floor/tile glue}, $t_{14}=1$ (ID 41): Dry mortar, adhesive for facades and tiles, Consumption 3-4kg/m2, Weber Therm Elastik (Weber).
    \item $i=15$, \texttt{wall, floor/tiles}, $t_{15}=1$ (ID 44): Ceramic glazed tile, 20 kg/m2.
    \item $i=16$, \texttt{structure/steel profiles}, $t_{16}=3$ (ID 49): Structural steel profiles, generic, 100\% recycled content, I, H, U, L, and T sections, S235, S275 and S355.
    \item $i=17$, \texttt{wall/door i}, $t_{17}=1$ (ID 50): Wooden and engineered wood interior doors, 1,23 m x 2,18 m, 72.8 kg/m2.
    \item $i=18$, \texttt{wall/door e}, $t_{18}=1$ (ID 53): Steel door with mineral wool core, fire-rated, 900x2000 mm, 74 kg/unit.
    \item $i=19$, \texttt{floor/thermal insulation p}, $t_{19}=1$ (ID 56): EPS insulation panels, graphite, L= 0.033 W/mK, R= 3.03 m2K/W, 100 mm, 3 kg/m2, 30 kg/m3, compressive strength 220 kPa, 100\% virgin polystyrene, Lambda=0.033 W/(m.K).
    \item $i=20$, \texttt{floor/concrete}, $t_{20}=3$ (ID 61): Ready-mix concrete, generic, C25/30 (3625/4350 PSI), XC1, 50\% fly ash, CEM IV/B pozzolanic cement,.
    \item $i=21$, \texttt{floor, roof/thermal insulation f}, $t_{21}=3$ (ID 64): PIR insulation board, aluminium facing, L=0.022 W/mK, R=3.6 m2K/W, 80 mm, 2.8 kg/m2, Lambda=0.022 W/(m.K) (UNILIN).
    \item $i=22$, \texttt{floor, roof/moisture insulation}, $t_{22}=2$ (ID 66): Bitumen sheets, G 200 S4; 5,0 kg/m2, EN15804+A1, ref. year 2018.
    \item $i=23$, \texttt{floor, roof/steel sheet}, $t_{23}=2$ (ID 69): Hot-dip galvanised steel sheet, 7874.0 kg/m3, EN15804+A1, ref. year 2018.
\end{itemize}
}

The solutions for the remaining building types (B1, C1, D1--D4, K1 and K2), obtained in exactly the same way, are reported in Appendix~\ref{app:MT} as index sequences.

\section{Optimal selection of the building materials}

Optimisation using the TOPSIS and Pareto sequence methods was carried out for two analysis scenarios. In the first scenario, denoted as M\&T, both the materials and their transport were included, with GWP values and the total cost of materials and transport analysed. In the second scenario, denoted as M, the analysis was limited to the materials only, accounting for their GWP values and costs while excluding the transport stage. Both optimisation procedures were performed independently for each analysed building variant, i.e. A1, B1, C1, D1, D2, D3, D4, K1 and K2. All numerical results reported in Section~\ref{sec:numerical} and in Appendix~\ref{app:MT} refer to the M\&T scenario, i.e. the values $c_{it}$ and $e_{it}$ of every variant include the transport of the material to the prefabrication plant. The corresponding results for the M scenario, obtained with the same procedures from material-only data, are reported in Appendix~\ref{app:M}.

The Pareto sequence, rather than the repeated TOPSIS, was adopted as the basis for the final selection for the following reasons. A weight sweep with a fixed step $\varepsilon$ cannot give a complete picture of the Pareto set. As the results in Section~\ref{sec:numerical} show, the selection made for several consecutive weights is often the same (e.g. $(0.9,0.1)$ and $(0.8,0.2)$ for A1), so that the solution moves along one line segment. On the other hand, the weight at which the preferred selection moves to the next segment can fall between the two examined points, and, as a consequence, a nondominated solution can be missed. The number of distinct TOPSIS selections is at most $H+1$ and depends on the assumed $\varepsilon$ and, through the normalisation, on the enumerated set of alternatives. On the other hand, the Pareto sequence method requires no weights, returns every vertex of the bottom-left (i.e. Pareto-efficient) boundary of the convex hull exactly once, and gives the price of reducing the emission $\rho=\Delta E/\Delta C$ for each step explicitly. The returned points can be inspected directly, which is the basis of the geometric selection rule described below.

\begin{figure}[htbp]
\centering
\includegraphics[width = 1\textwidth]{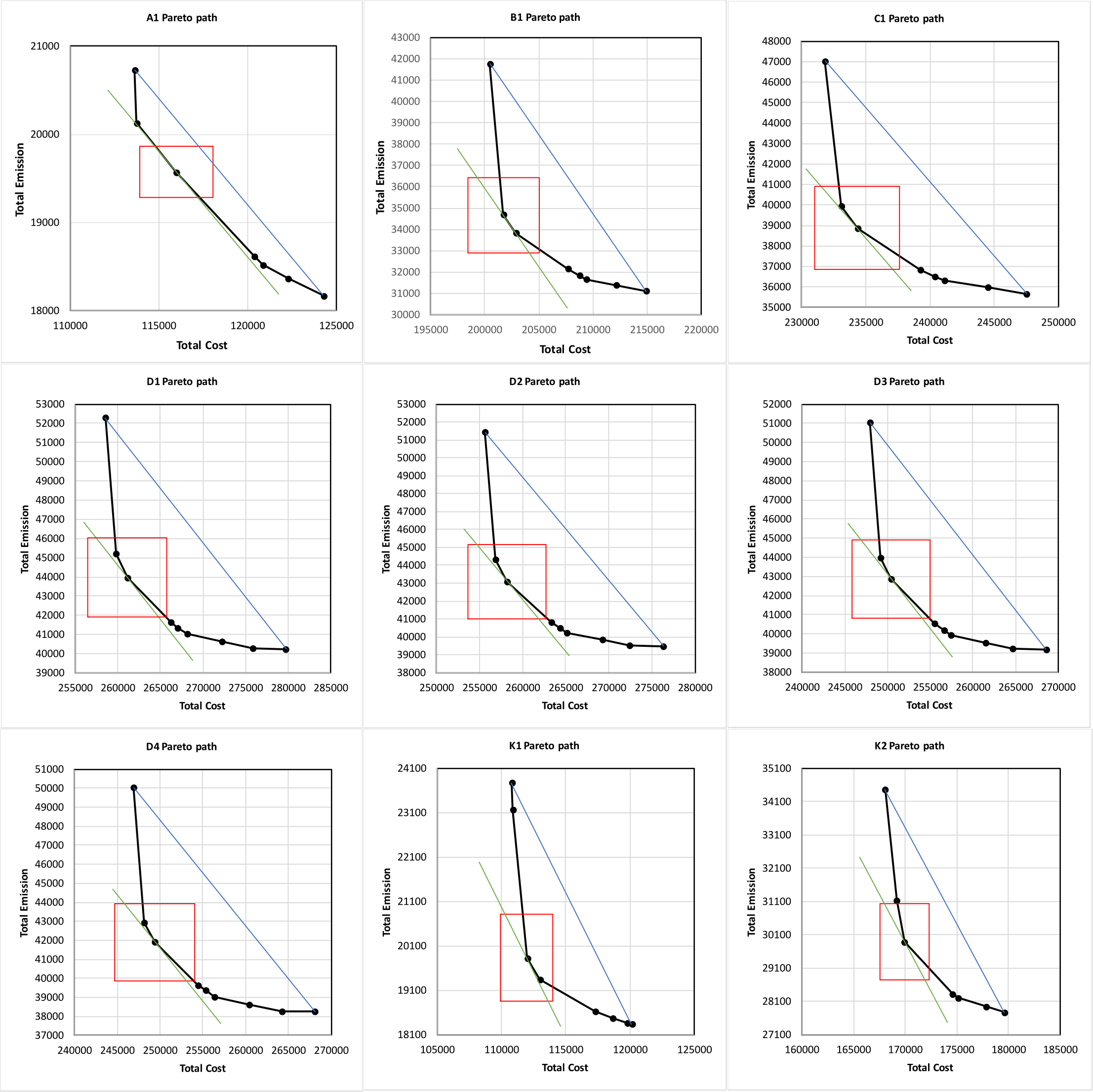}
\caption{Pareto fronts for the analysed building variants considering material- and transport-related GWP [kgCO\textsubscript{2}e] and costs [PLN] (M\&T scenario).}
Source: own.\\
\label{fig:Fig 3}
\end{figure}

\begin{figure}[htbp]
\centering
\includegraphics[width = 1\textwidth]{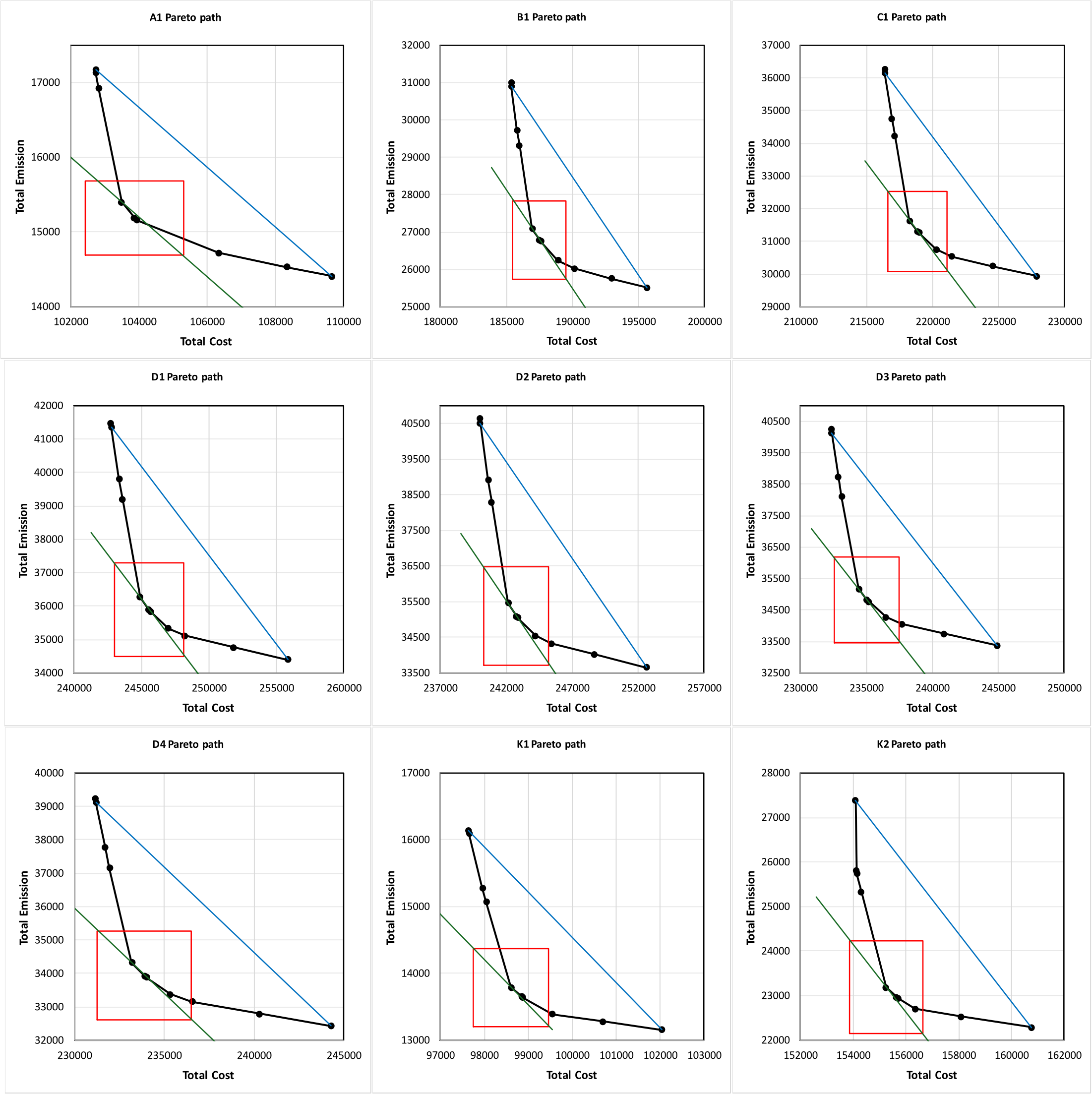}
\caption{Pareto fronts for the analysed building variants considering material-related GWP [kgCO\textsubscript{2}e] and costs [PLN] only (M scenario).}
Source: own.\\
\label{fig:Fig 4}
\end{figure}

The Pareto sequences obtained for the analysed building variants are shown in Figures \ref{fig:Fig 3} and \ref{fig:Fig 4}. To identify the recommended solution, a geometric method based on the Pareto frontier was applied. First, a straight line connecting the two extreme Pareto solutions, corresponding to the minimum-cost and minimum-emission alternatives (blue line), was drawn. Next, a line parallel to this segment (green line) was shifted towards the Pareto frontier until it first became tangent to the set of analysed solutions. The point of tangency represents the solution with the maximum perpendicular distance from the line connecting the extreme solutions and therefore corresponds to the most balanced trade-off between total cost and carbon emissions. This variant was adopted as the primary optimum solution.

\begin{figure}[htbp]
\centering
\includegraphics[width = 1\textwidth]{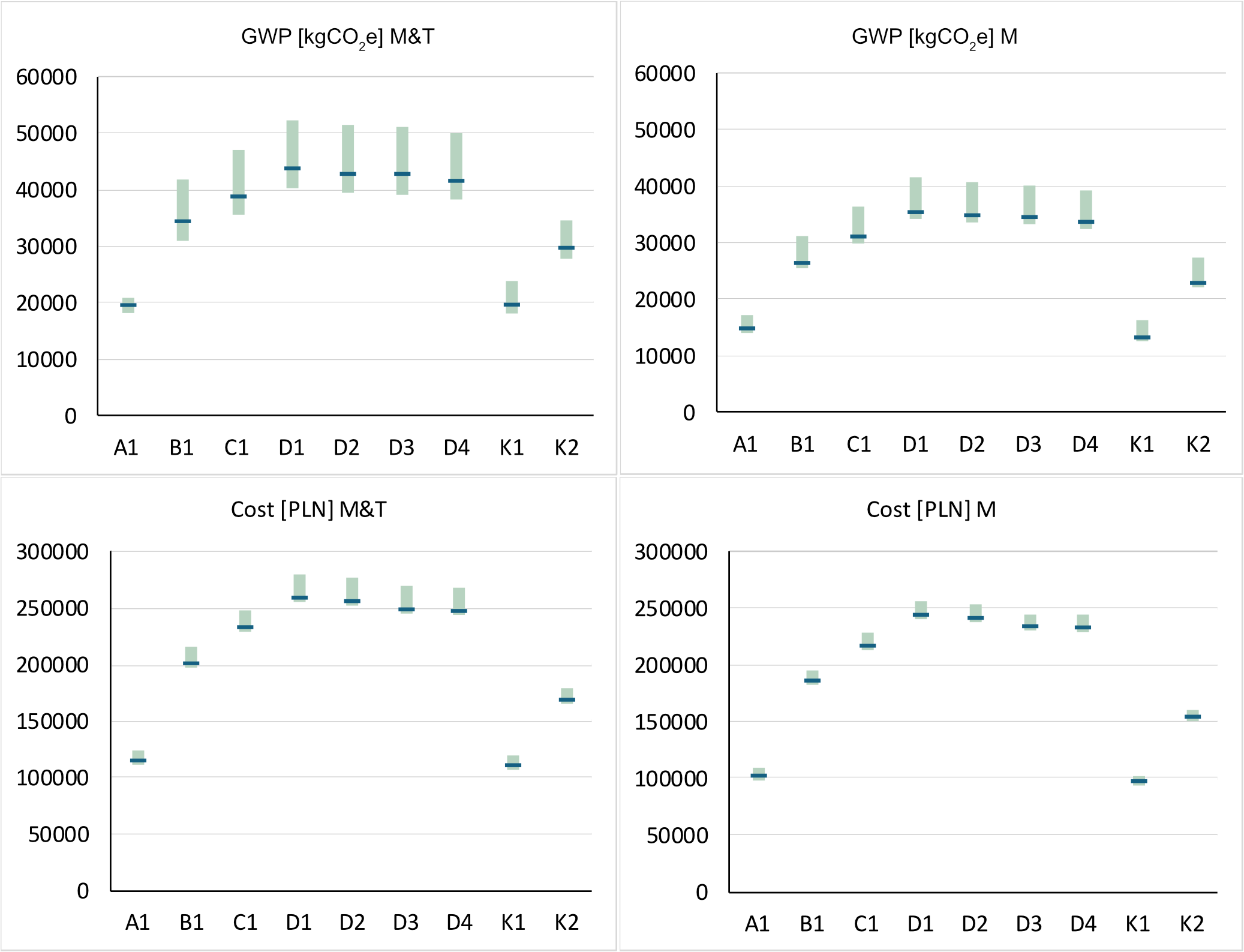}
\caption{Minimum, selected optimal, and maximum values of GWP and cost for each building variant under the M\&T and M scenarios.}
Source: own.\\
\label{fig:Fig 5}
\end{figure}

Figure \ref{fig:Fig 5} presents a comparison of the minimum-maximum ranges and the optimum solution for the individual building variants. The lower values for variant M than for M\&T result from excluding transport costs and transport-related emissions. The results for K1 are similar to those for A1. Both variants comprise two modules. Variants K2 and B1 each comprise four modules. However, B1 has more internal walls, resulting in slightly higher values than K2. All D variants (D1--D4) have similar values. The range in GWP values is greater than that for cost and averages 26.0\% for GWP in the M\&T variant, 20.2\% for GWP in the M variant, 7.9\% for cost in the M\&T variant and 5.3\% for cost in the M variant. The percentage range was calculated relative to the optimum value.

In addition to the reference solution corresponding to the tangency point, nearby Pareto-sequence solutions were considered acceptable alternatives. Their selection was based on a rectangular neighbourhood centred at the reference point $(C^\star, E^\star)$, with sides parallel to the cost and emission axes. Let $C_{\max}$ and $E_{\max}$ denote the maximum cost and maximum GWP in the Pareto sequence, attained at the minimum-emission and minimum-cost endpoints, respectively. The half-width and half-height of the rectangle were defined independently as

\begin{align*}
\Delta C &=\frac{C_{\max}-C^\star}{4},\\
\Delta E &=\frac{E_{\max}-E^\star}{4}.
\end{align*}

\noindent{}A Pareto-sequence solution $(C, E)$ was classified as acceptable if it satisfied both conditions:

\begin{align*}
|C-C^\star| &<\Delta C,\\
|E-E^\star| &<\Delta E.
\end{align*}

The red rectangles in Figures \ref{fig:Fig 3} and \ref{fig:Fig 4} illustrate these neighbourhoods. The factor of one quarter was adopted as a geometric selection parameter. The mean absolute percentage deviations of the accepted alternative solutions from their respective reference solutions were approximately 0.6\% in cost and 2.6\% in GWP for the M\&T scenario, and 0.3\% in cost and 1.0\% in GWP for the M scenario. They indicate that alternative material combinations within the selected neighbourhoods can provide similar cost and GWP outcomes. The building-material variant numbers selected across all optimal variants were compared to determine which occurred most frequently (Table \ref{tab:materials} -- 1st) and which occurred second most frequently (Table \ref{tab:materials} -- 2nd). Colours were used to indicate differences between the M\&T and M variants. The most frequently occurring materials are recommended for all building variants.

\begin{table}[ht]
\centering
\caption{Most frequently selected materials in the M\&T and M analyses. Yellow indicates differences in the first-choice selection, orange indicates differences in the second-choice selection based on a single occurrence, and blue indicates differences in the second-choice selection.}
\begin{scriptsize}
\begin{tabular}{*{24}{c}}
\hline
- & 1 & 2 & 3 & 4 & 5 & 6 & 7 & 8 & 9 & 10 & 11 & 12 & 13 & 14 & 15 & 16 & 17 & 18 & 19 & 20 & 21 & 22 & 23 \\
\hline
1st M\&T & 3 & 4 & 1 & 3 & 3 & 1 & 3 & 1 & 3 & 3 & 1 & 1 & 1 & 1 & 1 & 3 & 2 & 2 & 3 & 3 & 3 & 2 & 1 \\
2nd M\&T & - & - & - & 2 & - & - & - & - & - & - & - & - & - & - & - & - & - & \cellcolor{orange!30}1 & 1 & - & - & - & - \\
1st M & 3 & 4 & \cellcolor{yellow!30}3 & 3 & 3 & 1 & \cellcolor{yellow!30}2 & 1 & 3 & \cellcolor{yellow!30}2 & 1 & 1 & 1 & 1 & 1 & 3 & \cellcolor{yellow!30}3 & 2 & \cellcolor{yellow!30}2 & 3 & 3 & 2 & 1 \\
2nd M & - & - & - & \cellcolor{cyan!20} - & - & - & - & - & - & \cellcolor{cyan!20}3 & - & - & \cellcolor{cyan!20}2 & - & - & - & \cellcolor{cyan!20}2 & - & - & - & - & - & - \\
\hline
\end{tabular}
\end{scriptsize}
Source: own.\\
\label{tab:materials}
\end{table}

\begin{figure}[htbp]
\centering
\includegraphics[width = 1\textwidth]{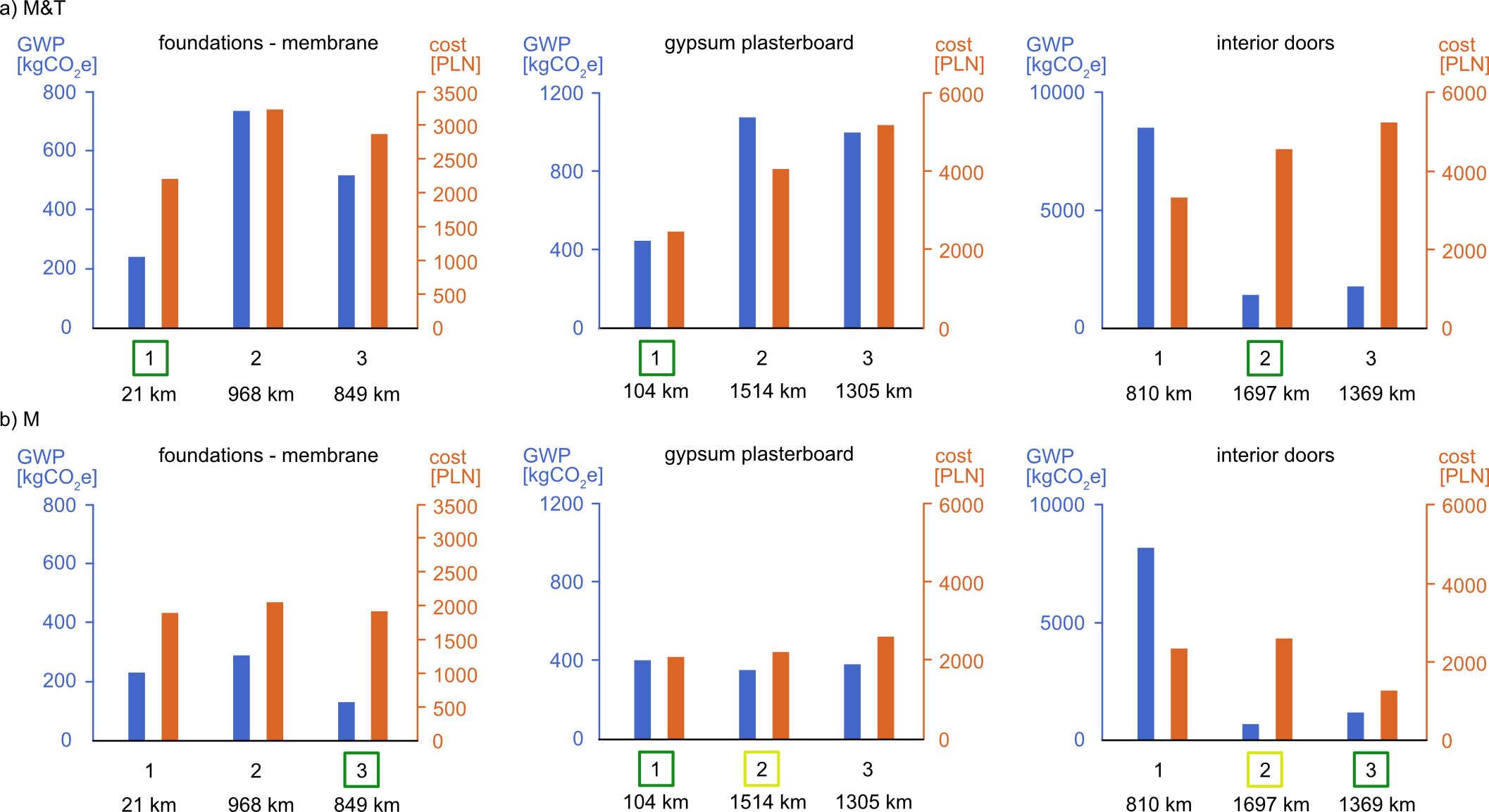}
\caption{GWP and cost of examples of the most frequently and second most frequently selected construction materials. The green outline indicates the material variant most frequently identified as optimal, the yellow outline indicates the second most frequently selected. Transport distances are shown below the material variant numbers: (a) M\&T analysis; (b) M analysis. The values are the totals for building type D4 (foundation membrane: group 3, gypsum plasterboard: wall underlay, group 13, interior doors: group 17).}
Source: own.\\
\label{fig:Fig 6}
\end{figure}

Figure \ref{fig:Fig 6} compares the GWP values and costs of the most and second most frequently chosen building materials in the M\&T and M analyses, using the example of materials for which a difference was observed.

\begin{figure}[htbp]
\centering
\includegraphics[width = 1\textwidth]{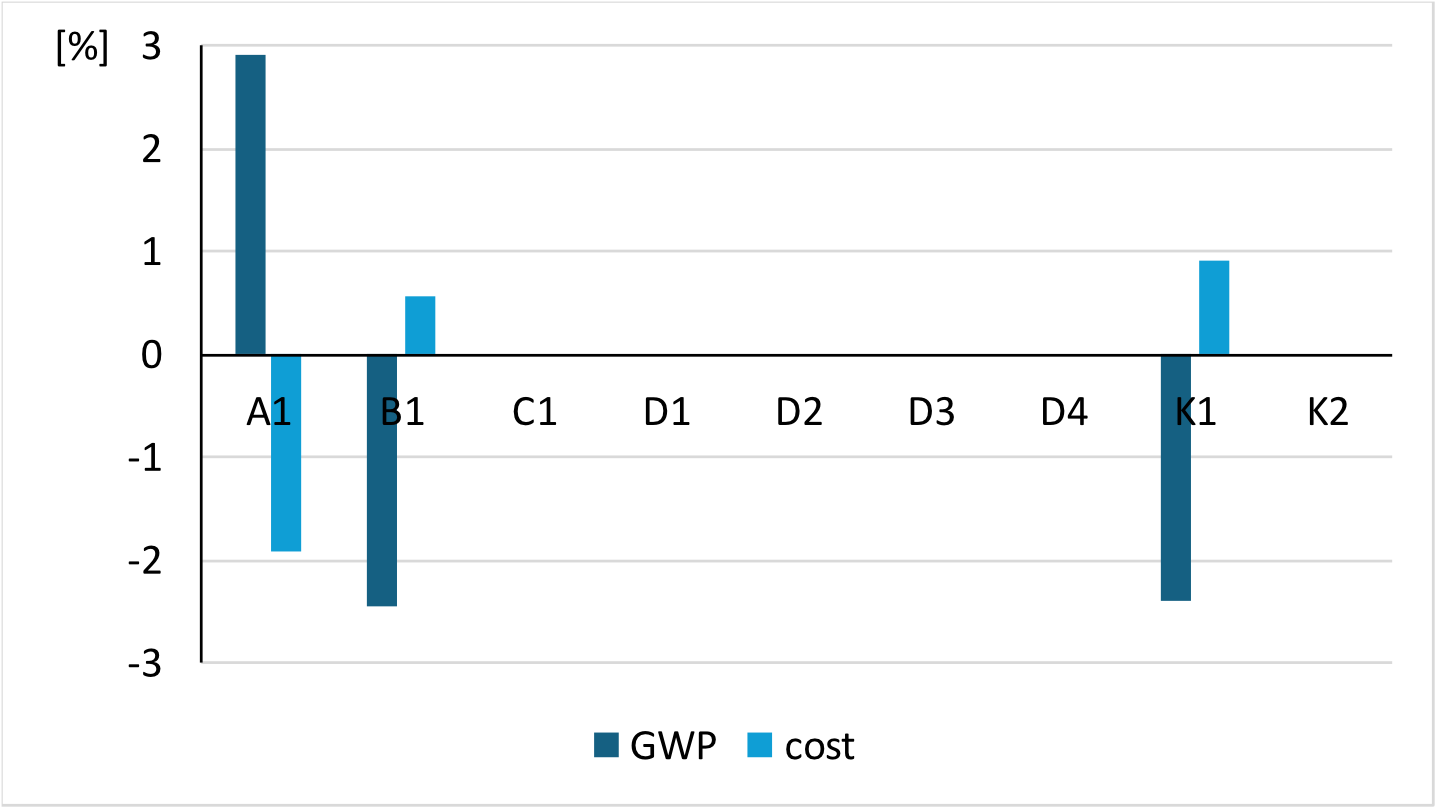}
\caption{Percentage differences in GWP [kgCO\textsubscript{2}e] and cost [PLN] between the individually optimised building variants and the corresponding variants using the recommended material set derived from the M\&T analysis. Negative values indicate lower values for the recommended-material variant.}
Source: own.\\
\label{fig:Fig 7}
\end{figure}

The primary optimum variant was compared with the variant using the recommended materials (1st M\&T). The results show small differences in GWP values (Figure \ref{fig:Fig 7}). The mean absolute difference in GWP values was 0.9\% for all building types and 2.6\% for buildings A1, B1 and K1 (these were the only cases in which the optimum and recommended variants differed). The cost differences were even smaller: the mean absolute cost difference was 0.4\% for all building types and 1.1\% for buildings A1, B1 and K1. A negative value denotes a lower value for the recommended variant than that for the optimum variant.

\begin{figure}[htbp]
\centering
\includegraphics[width = 1\textwidth]{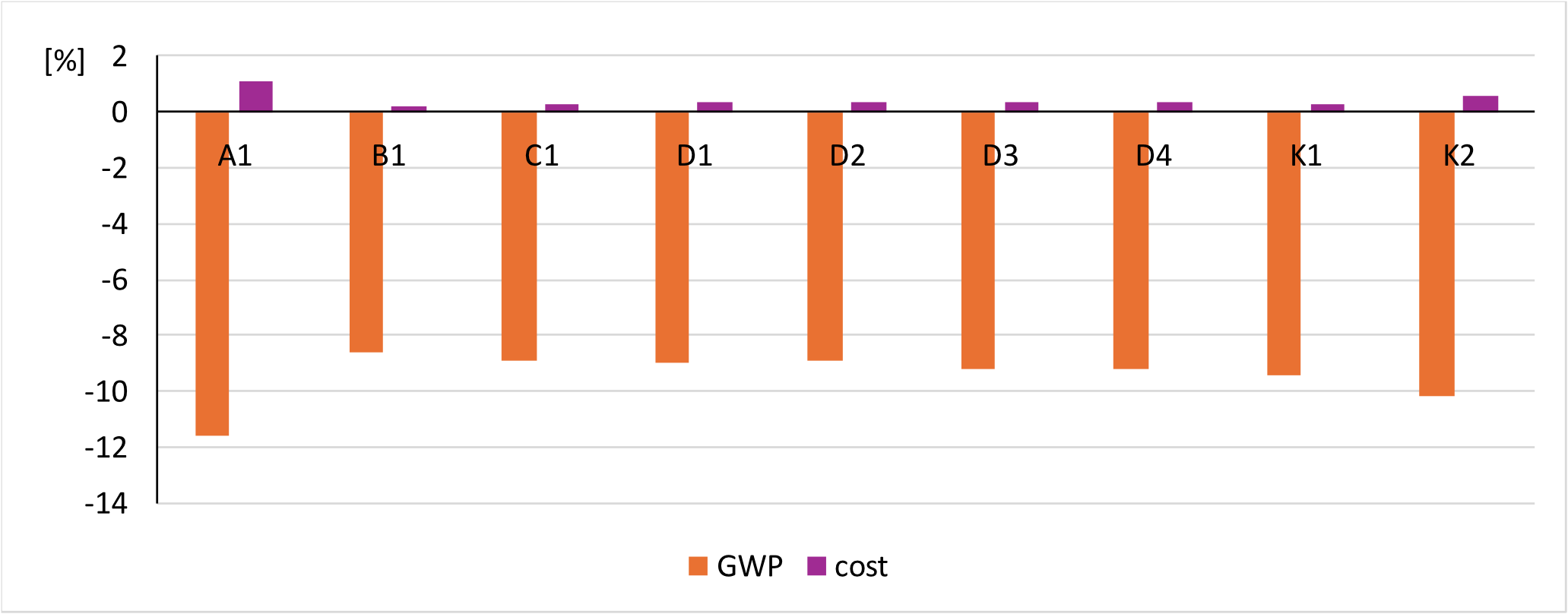}
\caption{Percentage differences in material-related GWP and cost between the recommended material sets derived from the M analysis (1st M) and the M\&T analysis (1st M\&T), evaluated for each building using material-only data. The 1st M\&T set provides the reference values. Negative values indicate lower values for the 1st M set.}
Source: own.\\
\label{fig:Fig 8}
\end{figure}

Figure \ref{fig:Fig 8} compares the recommended material sets derived from the material-only analysis (1st M) and the material-and-transport analysis (1st M\&T) for each building variant. Both sets were evaluated using material-only GWP and cost data, excluding transport. Percentage differences were calculated relative to the material-only values of the 1st M\&T set. Across all analysed buildings, the 1st M set had lower material-related GWP (-9.4\%) and a slightly higher material cost (+0.4\%) than the 1st M\&T set. Negative values indicate lower GWP or cost for the 1st M set.

\begin{figure}[htbp]
\centering
\includegraphics[width = 1\textwidth]{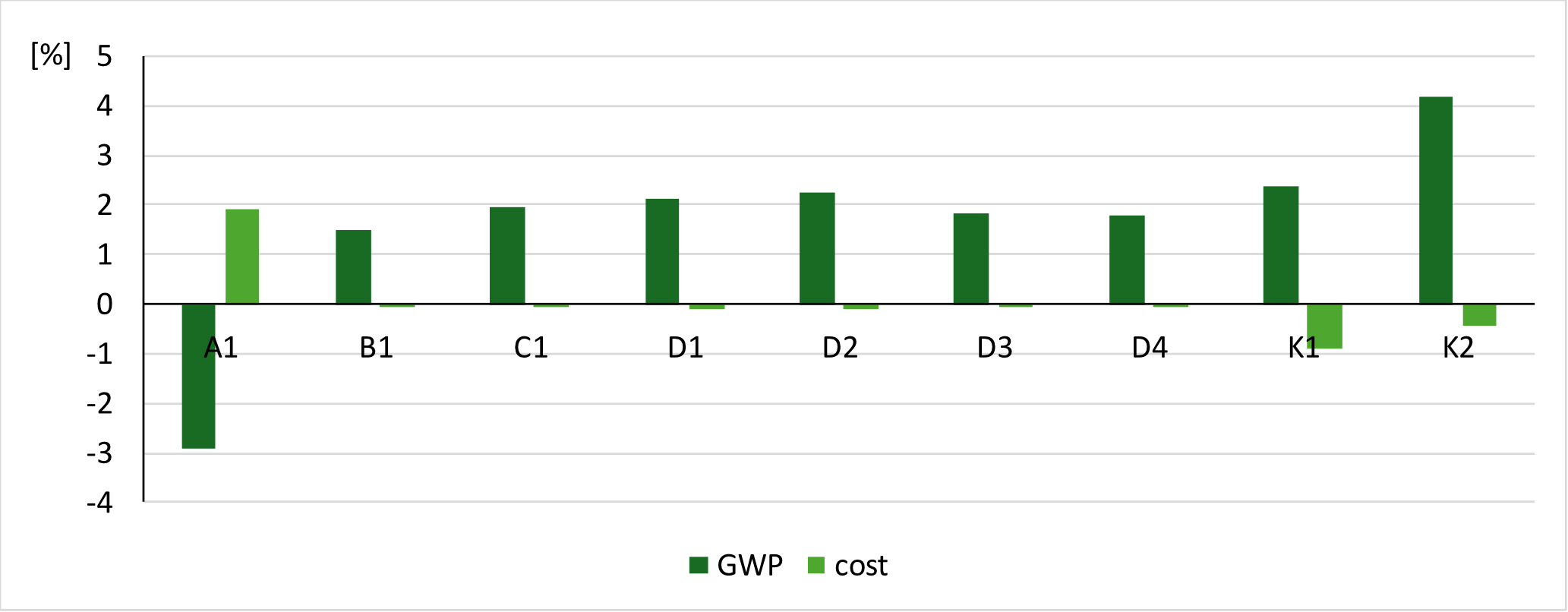}
\caption{Percentage differences in GWP [kgCO\textsubscript{2}e] and cost [PLN] between building variants using the first- and second-choice recommended material sets identified in the M\&T analysis. Negative values indicate lower values for the second-choice material set.}
Source: own.\\
\label{fig:Fig 9}
\end{figure}

The results for the first- and second-choice recommended material sets (1st and 2nd M\&T) were compared. Across all nine building variants, the mean absolute percentage difference was 2.3\% for GWP and 0.4\% for cost (Figure \ref{fig:Fig 9}). Negative values indicate lower values for the second-choice set relative to the first-choice set.

In modular construction, there is no single universal project delivery method. The Design-Bid-Build (DBB) and Design-Build (DB) approaches are among the most commonly used methods in the construction industry. The DB approach facilitates the early involvement of modular manufacturers in key design decisions concerning module dimensions, structural systems, and material selection. When interpreting cost differences between the analysed alternatives, it is important to recognise that no universal percentage threshold defines their significance. Analysis of building projects implemented using the DBB model revealed an average coefficient of variation of contractors' bid prices of 12.6\% and a median of 10.0\% \cite{nguyen_best-value_2018}. This bid dispersion was considered relatively small. In the present study, the cost variation  (Figure \ref{fig:Fig 5}) does not exceed 8\%, while the mean absolute percentage differences between the analysed alternatives range from 0.4\% to 1.1\%. These findings suggest relatively minor economic differences between the analysed solutions compared with the variability in bid prices observed in construction practice, although the two measures are not directly comparable. Building environmental impact mitigation plans target a 40\% reduction in CO2 emissions by 2030, but no thresholds were defined to consider these changes significant. The LEED multi-criteria certification system, which evaluates buildings for their sustainability, sets GWP reduction targets of at least 5\%, 10\%, and 20\% compared to a reference building (USGBC, 2020), so a reduction of less than 5\% can be considered insignificant. The range between the extreme GWP values in the analyses presented, amounting to 20.2--26.0\%, is significant. Differences between the analysed variants (Figures \ref{fig:Fig 7} and \ref{fig:Fig 9}) of 2.3--2.6\% should be regarded as relatively small. Excluding transport from the material selection procedure changes the recommended set of materials. When both recommended sets are evaluated using material-only data, the set derived from the M analysis has lower GWP for all analysed buildings, accompanied by a small increase in material cost (Figure \ref{fig:Fig 8}).

\section{Final Remarks}

The selection of construction materials for modular buildings was formulated as a bi-objective combinatorial problem with additive cost and carbon-footprint criteria. Two solution algorithms were analysed on nine buildings developed within the Rail4Earth project. The main findings are as follows. First, because both criteria are additive over independent component groups, the problem admits a simple exact description of the convexified Pareto frontier. The greedy Pareto sequence method returns all vertices of the bottom-left (i.e. Pareto-efficient) boundary of the convex hull. Second, repeatedly used TOPSIS returns a subset of nondominated solutions that depends on the weight step and, through the normalisation, on the enumerated set of alternatives. It may select nondominated discrete points that the Pareto sequence omits, but it cannot guarantee a complete picture of the cost-emission trade-off. Third, the geometric rule applied to the Pareto sequence provides a stable set of recommended materials. The balanced solutions of the individual buildings differ from the common recommended set by at most 3\% in carbon footprint and by at most 2\% in cost, and the cost spread between the extreme solutions (below 8\%) is small compared with the dispersion observed in construction bids. On the other hand, the corresponding spread in carbon footprint (20-26\%) is practically significant. Finally, including the transport stage changes the recommended material set, confirming that transport-related emissions and costs should be considered in selecting materials for buildings assembled from products sourced across Europe.

The study has several limitations. The environmental data combine product-specific EPDs with generic datasets, and part of the cost data relies on calibrated estimates generated by LLMs. Both sources introduce uncertainty that was not propagated through the optimisation. Only two criteria were considered, and the discrete Pareto set was approximated by its convex hull. The environmental assessment combines material production impacts for stages A1--A3 with direct fuel-combustion emissions from delivery journeys to the prefabrication plant. The transport component excludes upstream fuel-supply emissions and empty return journeys and therefore does not represent the full life-cycle greenhouse gas impact of transport. Future work could extend the model to further life-cycle stages and impact categories, incorporate uncertainty in the input data, and examine compatibility constraints between component groups, which were not assumed in the current research.

\section*{Acknowledgment}
This research was co-funded by the European Union under Horizon Europe FP4-Rail4EARTH project `Sustainable and green rail systems' (GA 101101917) and co-funded by the Polish Ministry of Science and Higher Education (former Ministry of Education and Science) under the `International Co-Financed Projects' programme (Agreement No. 5400/HE/2023/2). The research was also partially supported by the SUMOC project co-financed by the National Science Centre (NCN) through grant 2024/55/I/HS4/01297, the Slovenian Research and Innovation Agency (ARIS) through grant N2-0434 and the Austrian Science Fund (FWF) through grant 10.55776/PIN4864224.

\section*{Disclosure}
The authors used various AI tools to support the literature review and imputation of missing data, to verify the correctness of computations and the consistency of notions and notations used in the text, and to correct language errors in the final version of the paper. The authors take full responsibility for the contents of the paper.

\bibliographystyle{apalike}
\bibliography{staz}

\begin{appendices}

\section{Numerical results for the remaining building types (M\&T scenario)}\label{app:MT}

This appendix reports, for the building types B1, C1, D1--D4, K1 and K2, the TOPSIS weight-sweep solutions and the Pareto sequences obtained with the procedures presented in Section~\ref{sec:numerical} from the data of the M\&T scenario (materials and transport). The notation, the dominance filtering and the labelling of the types are the same as for building type A1 in Section~\ref{sec:numerical}. The Pareto sequences are those shown in Figure~\ref{fig:Fig 3}.

\subsection*{Building type B1}

\paragraph{TOPSIS solutions for a weight sweep with $\varepsilon=0.1$.}
The results are ordered from the cost-only to the emission-only weighting.

{\footnotesize
\begin{description}
    \item[$(1,0)$] $C=200485.840$, $E=41755.060$, $S=1.000000$,
    $a^\star=$\texttt{(3, 4, 1, 2, 3, 1, 3, 1, 3, 3, 1, 1, 1, 1, 1, 3, 1, 2, 3, 3, 3, 2, 1)}.

    \item[$(0.9,0.1)$] $C=201719.590$, $E=34672.470$, $S=0.844234$,
    $a^\star=$\texttt{(3, 4, 1, 2, 3, 1, 3, 1, 3, 3, 1, 1, 1, 1, 1, 3, 2, 2, 3, 3, 3, 2, 1)}.

    \item[$(0.8,0.2)$] $C=202878.880$, $E=33827.010$, $S=0.783376$,
    $a^\star=$\texttt{(3, 4, 1, 3, 3, 1, 3, 1, 3, 3, 1, 1, 1, 1, 1, 3, 2, 2, 3, 3, 3, 2, 1)}.

    \item[$(0.7,0.3)$] $C=204636.180$, $E=33306.810$, $S=0.771701$,
    $a^\star=$\texttt{(3, 4, 1, 3, 3, 1, 3, 1, 3, 2, 1, 1, 1, 1, 1, 3, 2, 1, 3, 3, 3, 2, 1)}.

    \item[$(0.6,0.4)$] $C=207632.930$, $E=32160.910$, $S=0.818486$,
    $a^\star=$\texttt{(3, 4, 1, 3, 3, 1, 2, 1, 3, 3, 1, 1, 1, 1, 1, 3, 2, 2, 3, 3, 3, 2, 1)}.

    \item[$(0.5,0.5)$] $C=208701.770$, $E=31830.850$, $S=0.860028$,
    $a^\star=$\texttt{(3, 4, 1, 3, 3, 1, 2, 1, 3, 3, 1, 1, 1, 1, 1, 3, 2, 1, 3, 3, 3, 2, 1)}.

    \item[$(0.4,0.6)$] $C=209390.230$, $E=31640.710$, $S=0.895652$,
    $a^\star=$\texttt{(3, 4, 1, 3, 3, 1, 2, 1, 3, 2, 1, 1, 1, 1, 1, 3, 2, 1, 3, 3, 3, 2, 1)}.

    \item[$(0.3,0.7)$] $C=209390.230$, $E=31640.710$, $S=0.921087$,
    $a^\star=$\texttt{(3, 4, 1, 3, 3, 1, 2, 1, 3, 2, 1, 1, 1, 1, 1, 3, 2, 1, 3, 3, 3, 2, 1)}.

    \item[$(0.2,0.8)$] $C=212100.250$, $E=31365.330$, $S=0.947165$,
    $a^\star=$\texttt{(3, 4, 1, 3, 3, 1, 2, 1, 3, 2, 1, 1, 1, 1, 1, 3, 2, 1, 3, 3, 3, 2, 2)}.

    \item[$(0.1,0.9)$] $C=214924.740$, $E=31093.110$, $S=0.974167$,
    $a^\star=$\texttt{(3, 4, 1, 3, 3, 1, 2, 1, 3, 2, 1, 2, 1, 1, 1, 3, 2, 1, 3, 3, 3, 2, 2)}.

    \item[$(0,1)$] $C=214924.740$, $E=31093.110$, $S=1.000000$,
    $a^\star=$\texttt{(3, 4, 1, 3, 3, 1, 2, 1, 3, 2, 1, 2, 1, 1, 1, 3, 2, 1, 3, 3, 3, 2, 2)}.

\end{description}
}

\paragraph{Pareto-sequence solution.}
The successive one-group changes are:

\begin{center}
\small
\begin{tabular}{c l r r r}
\hline
Step & Change & $\Delta C$ & $\Delta E$ & $\rho=\Delta E/\Delta C$ \\
\hline
1 & $i=17:\ t_{17}=1\to 2$ & 1233.750 & 7082.590 & 5.741 \\
2 & $i=4:\ t_{4}=2\to 3$ & 1159.290 & 845.460 & 0.729 \\
3 & $i=7:\ t_{7}=3\to 2$ & 4754.050 & 1666.100 & 0.350 \\
4 & $i=18:\ t_{18}=2\to 1$ & 1068.840 & 330.060 & 0.309 \\
5 & $i=10:\ t_{10}=3\to 2$ & 688.460 & 190.140 & 0.276 \\
6 & $i=23:\ t_{23}=1\to 2$ & 2710.020 & 275.380 & 0.102 \\
7 & $i=12:\ t_{12}=1\to 2$ & 2824.490 & 272.220 & 0.096 \\
\hline
\end{tabular}
\end{center}

\paragraph{Complete list of Pareto-sequence solutions.}
The corresponding solution vectors are:

{\footnotesize
\begin{description}
    \item[Step 0] $C=200485.840$, $E=41755.060$,
    $a^{(0)}=$\texttt{(3, 4, 1, 2, 3, 1, 3, 1, 3, 3, 1, 1, 1, 1, 1, 3, 1, 2, 3, 3, 3, 2, 1)}.

    \item[Step 1] $C=201719.590$, $E=34672.470$,
    $a^{(1)}=$\texttt{(3, 4, 1, 2, 3, 1, 3, 1, 3, 3, 1, 1, 1, 1, 1, 3, 2, 2, 3, 3, 3, 2, 1)}.

    \item[Step 2] $C=202878.880$, $E=33827.010$,
    $a^{(2)}=$\texttt{(3, 4, 1, 3, 3, 1, 3, 1, 3, 3, 1, 1, 1, 1, 1, 3, 2, 2, 3, 3, 3, 2, 1)}.

    \item[Step 3] $C=207632.930$, $E=32160.910$,
    $a^{(3)}=$\texttt{(3, 4, 1, 3, 3, 1, 2, 1, 3, 3, 1, 1, 1, 1, 1, 3, 2, 2, 3, 3, 3, 2, 1)}.

    \item[Step 4] $C=208701.770$, $E=31830.850$,
    $a^{(4)}=$\texttt{(3, 4, 1, 3, 3, 1, 2, 1, 3, 3, 1, 1, 1, 1, 1, 3, 2, 1, 3, 3, 3, 2, 1)}.

    \item[Step 5] $C=209390.230$, $E=31640.710$,
    $a^{(5)}=$\texttt{(3, 4, 1, 3, 3, 1, 2, 1, 3, 2, 1, 1, 1, 1, 1, 3, 2, 1, 3, 3, 3, 2, 1)}.

    \item[Step 6] $C=212100.250$, $E=31365.330$,
    $a^{(6)}=$\texttt{(3, 4, 1, 3, 3, 1, 2, 1, 3, 2, 1, 1, 1, 1, 1, 3, 2, 1, 3, 3, 3, 2, 2)}.

    \item[Step 7] $C=214924.740$, $E=31093.110$,
    $a^{(7)}=$\texttt{(3, 4, 1, 3, 3, 1, 2, 1, 3, 2, 1, 2, 1, 1, 1, 3, 2, 1, 3, 3, 3, 2, 2)}.

\end{description}
}

\subsection*{Building type C1}

\paragraph{TOPSIS solutions for a weight sweep with $\varepsilon=0.1$.}
The results are ordered from the cost-only to the emission-only weighting.

{\footnotesize
\begin{description}
    \item[$(1,0)$] $C=231841.220$, $E=47037.600$, $S=1.000000$,
    $a^\star=$\texttt{(3, 4, 1, 2, 3, 1, 3, 1, 3, 3, 1, 1, 1, 1, 1, 3, 1, 2, 3, 3, 3, 2, 1)}.

    \item[$(0.9,0.1)$] $C=233074.970$, $E=39955.010$, $S=0.834105$,
    $a^\star=$\texttt{(3, 4, 1, 2, 3, 1, 3, 1, 3, 3, 1, 1, 1, 1, 1, 3, 2, 2, 3, 3, 3, 2, 1)}.

    \item[$(0.8,0.2)$] $C=234333.490$, $E=38876.570$, $S=0.770803$,
    $a^\star=$\texttt{(3, 4, 1, 3, 3, 1, 3, 1, 3, 3, 1, 1, 1, 1, 1, 3, 2, 2, 3, 3, 3, 2, 1)}.

    \item[$(0.7,0.3)$] $C=239267.820$, $E=36852.160$, $S=0.768604$,
    $a^\star=$\texttt{(3, 4, 1, 3, 3, 1, 2, 1, 3, 3, 1, 1, 1, 1, 1, 3, 2, 2, 3, 3, 3, 2, 1)}.

    \item[$(0.6,0.4)$] $C=239267.820$, $E=36852.160$, $S=0.821296$,
    $a^\star=$\texttt{(3, 4, 1, 3, 3, 1, 2, 1, 3, 3, 1, 1, 1, 1, 1, 3, 2, 2, 3, 3, 3, 2, 1)}.

    \item[$(0.5,0.5)$] $C=241096.500$, $E=36298.790$, $S=0.862342$,
    $a^\star=$\texttt{(3, 4, 1, 3, 3, 1, 2, 1, 3, 2, 1, 1, 1, 1, 1, 3, 2, 1, 3, 3, 3, 2, 1)}.

    \item[$(0.4,0.6)$] $C=241096.500$, $E=36298.790$, $S=0.896550$,
    $a^\star=$\texttt{(3, 4, 1, 3, 3, 1, 2, 1, 3, 2, 1, 1, 1, 1, 1, 3, 2, 1, 3, 3, 3, 2, 1)}.

    \item[$(0.3,0.7)$] $C=241096.500$, $E=36298.790$, $S=0.920045$,
    $a^\star=$\texttt{(3, 4, 1, 3, 3, 1, 2, 1, 3, 2, 1, 1, 1, 1, 1, 3, 2, 1, 3, 3, 3, 2, 1)}.

    \item[$(0.2,0.8)$] $C=244471.850$, $E=35961.480$, $S=0.946855$,
    $a^\star=$\texttt{(3, 4, 1, 3, 3, 1, 2, 1, 3, 2, 1, 1, 1, 1, 1, 3, 2, 1, 3, 3, 3, 2, 2)}.

    \item[$(0.1,0.9)$] $C=247539.710$, $E=35665.820$, $S=0.974120$,
    $a^\star=$\texttt{(3, 4, 1, 3, 3, 1, 2, 1, 3, 2, 1, 2, 1, 1, 1, 3, 2, 1, 3, 3, 3, 2, 2)}.

    \item[$(0,1)$] $C=247539.710$, $E=35665.820$, $S=1.000000$,
    $a^\star=$\texttt{(3, 4, 1, 3, 3, 1, 2, 1, 3, 2, 1, 2, 1, 1, 1, 3, 2, 1, 3, 3, 3, 2, 2)}.

\end{description}
}

\paragraph{Pareto-sequence solution.}
The successive one-group changes are:

\begin{center}
\small
\begin{tabular}{c l r r r}
\hline
Step & Change & $\Delta C$ & $\Delta E$ & $\rho=\Delta E/\Delta C$ \\
\hline
1 & $i=17:\ t_{17}=1\to 2$ & 1233.750 & 7082.590 & 5.741 \\
2 & $i=4:\ t_{4}=2\to 3$ & 1258.520 & 1078.440 & 0.857 \\
3 & $i=7:\ t_{7}=3\to 2$ & 4934.330 & 2024.410 & 0.410 \\
4 & $i=18:\ t_{18}=2\to 1$ & 1068.840 & 330.060 & 0.309 \\
5 & $i=10:\ t_{10}=3\to 2$ & 759.840 & 223.310 & 0.294 \\
6 & $i=23:\ t_{23}=1\to 2$ & 3375.350 & 337.310 & 0.100 \\
7 & $i=12:\ t_{12}=1\to 2$ & 3067.860 & 295.660 & 0.096 \\
\hline
\end{tabular}
\end{center}

\paragraph{Complete list of Pareto-sequence solutions.}
The corresponding solution vectors are:

{\footnotesize
\begin{description}
    \item[Step 0] $C=231841.220$, $E=47037.600$,
    $a^{(0)}=$\texttt{(3, 4, 1, 2, 3, 1, 3, 1, 3, 3, 1, 1, 1, 1, 1, 3, 1, 2, 3, 3, 3, 2, 1)}.

    \item[Step 1] $C=233074.970$, $E=39955.010$,
    $a^{(1)}=$\texttt{(3, 4, 1, 2, 3, 1, 3, 1, 3, 3, 1, 1, 1, 1, 1, 3, 2, 2, 3, 3, 3, 2, 1)}.

    \item[Step 2] $C=234333.490$, $E=38876.570$,
    $a^{(2)}=$\texttt{(3, 4, 1, 3, 3, 1, 3, 1, 3, 3, 1, 1, 1, 1, 1, 3, 2, 2, 3, 3, 3, 2, 1)}.

    \item[Step 3] $C=239267.820$, $E=36852.160$,
    $a^{(3)}=$\texttt{(3, 4, 1, 3, 3, 1, 2, 1, 3, 3, 1, 1, 1, 1, 1, 3, 2, 2, 3, 3, 3, 2, 1)}.

    \item[Step 4] $C=240336.660$, $E=36522.100$,
    $a^{(4)}=$\texttt{(3, 4, 1, 3, 3, 1, 2, 1, 3, 3, 1, 1, 1, 1, 1, 3, 2, 1, 3, 3, 3, 2, 1)}.

    \item[Step 5] $C=241096.500$, $E=36298.790$,
    $a^{(5)}=$\texttt{(3, 4, 1, 3, 3, 1, 2, 1, 3, 2, 1, 1, 1, 1, 1, 3, 2, 1, 3, 3, 3, 2, 1)}.

    \item[Step 6] $C=244471.850$, $E=35961.480$,
    $a^{(6)}=$\texttt{(3, 4, 1, 3, 3, 1, 2, 1, 3, 2, 1, 1, 1, 1, 1, 3, 2, 1, 3, 3, 3, 2, 2)}.

    \item[Step 7] $C=247539.710$, $E=35665.820$,
    $a^{(7)}=$\texttt{(3, 4, 1, 3, 3, 1, 2, 1, 3, 2, 1, 2, 1, 1, 1, 3, 2, 1, 3, 3, 3, 2, 2)}.

\end{description}
}

\subsection*{Building type D1}

\paragraph{TOPSIS solutions for a weight sweep with $\varepsilon=0.1$.}
The results are ordered from the cost-only to the emission-only weighting.

{\footnotesize
\begin{description}
    \item[$(1,0)$] $C=258562.670$, $E=52286.850$, $S=1.000000$,
    $a^\star=$\texttt{(3, 4, 1, 2, 3, 1, 3, 1, 3, 3, 1, 1, 1, 1, 1, 3, 1, 2, 3, 3, 3, 2, 1)}.

    \item[$(0.9,0.1)$] $C=259796.420$, $E=45204.260$, $S=0.856289$,
    $a^\star=$\texttt{(3, 4, 1, 2, 3, 1, 3, 1, 3, 3, 1, 1, 1, 1, 1, 3, 2, 2, 3, 3, 3, 2, 1)}.

    \item[$(0.8,0.2)$] $C=261122.150$, $E=43968.050$, $S=0.787303$,
    $a^\star=$\texttt{(3, 4, 1, 3, 3, 1, 3, 1, 3, 3, 1, 1, 1, 1, 1, 3, 2, 2, 3, 3, 3, 2, 1)}.

    \item[$(0.7,0.3)$] $C=266216.230$, $E=41626.130$, $S=0.779748$,
    $a^\star=$\texttt{(3, 4, 1, 3, 3, 1, 2, 1, 3, 3, 1, 1, 1, 1, 1, 3, 2, 2, 3, 3, 3, 2, 1)}.

    \item[$(0.6,0.4)$] $C=266216.230$, $E=41626.130$, $S=0.822406$,
    $a^\star=$\texttt{(3, 4, 1, 3, 3, 1, 2, 1, 3, 3, 1, 1, 1, 1, 1, 3, 2, 2, 3, 3, 3, 2, 1)}.

    \item[$(0.5,0.5)$] $C=268159.600$, $E=41019.460$, $S=0.861871$,
    $a^\star=$\texttt{(3, 4, 1, 3, 3, 1, 2, 1, 3, 2, 1, 1, 1, 1, 1, 3, 2, 1, 3, 3, 3, 2, 1)}.

    \item[$(0.4,0.6)$] $C=268159.600$, $E=41019.460$, $S=0.893149$,
    $a^\star=$\texttt{(3, 4, 1, 3, 3, 1, 2, 1, 3, 2, 1, 1, 1, 1, 1, 3, 2, 1, 3, 3, 3, 2, 1)}.

    \item[$(0.3,0.7)$] $C=271793.910$, $E=40669.200$, $S=0.915591$,
    $a^\star=$\texttt{(3, 4, 1, 3, 3, 1, 2, 1, 3, 2, 1, 2, 1, 1, 1, 3, 2, 1, 3, 3, 3, 2, 1)}.

    \item[$(0.2,0.8)$] $C=272200.260$, $E=40620.210$, $S=0.942931$,
    $a^\star=$\texttt{(3, 4, 1, 3, 3, 1, 2, 1, 3, 2, 1, 1, 1, 1, 1, 3, 2, 1, 3, 3, 3, 2, 2)}.

    \item[$(0.1,0.9)$] $C=275834.570$, $E=40269.950$, $S=0.973000$,
    $a^\star=$\texttt{(3, 4, 1, 3, 3, 1, 2, 1, 3, 2, 1, 2, 1, 1, 1, 3, 2, 1, 3, 3, 3, 2, 2)}.

    \item[$(0,1)$] $C=279752.360$, $E=40228.530$, $S=1.000000$,
    $a^\star=$\texttt{(3, 4, 1, 3, 3, 1, 2, 1, 3, 2, 1, 2, 1, 1, 1, 3, 2, 1, 2, 3, 3, 2, 2)}.

\end{description}
}

\paragraph{Pareto-sequence solution.}
The successive one-group changes are:

\begin{center}
\small
\begin{tabular}{c l r r r}
\hline
Step & Change & $\Delta C$ & $\Delta E$ & $\rho=\Delta E/\Delta C$ \\
\hline
1 & $i=17:\ t_{17}=1\to 2$ & 1233.750 & 7082.590 & 5.741 \\
2 & $i=4:\ t_{4}=2\to 3$ & 1325.730 & 1236.210 & 0.932 \\
3 & $i=7:\ t_{7}=3\to 2$ & 5094.080 & 2341.920 & 0.460 \\
4 & $i=10:\ t_{10}=3\to 2$ & 874.530 & 276.610 & 0.316 \\
5 & $i=18:\ t_{18}=2\to 1$ & 1068.840 & 330.060 & 0.309 \\
6 & $i=23:\ t_{23}=1\to 2$ & 4040.660 & 399.250 & 0.099 \\
7 & $i=12:\ t_{12}=1\to 2$ & 3634.310 & 350.260 & 0.096 \\
8 & $i=19:\ t_{19}=3\to 2$ & 3917.790 & 41.420 & 0.011 \\
\hline
\end{tabular}
\end{center}

\paragraph{Complete list of Pareto-sequence solutions.}
The corresponding solution vectors are:

{\footnotesize
\begin{description}
    \item[Step 0] $C=258562.670$, $E=52286.850$,
    $a^{(0)}=$\texttt{(3, 4, 1, 2, 3, 1, 3, 1, 3, 3, 1, 1, 1, 1, 1, 3, 1, 2, 3, 3, 3, 2, 1)}.

    \item[Step 1] $C=259796.420$, $E=45204.260$,
    $a^{(1)}=$\texttt{(3, 4, 1, 2, 3, 1, 3, 1, 3, 3, 1, 1, 1, 1, 1, 3, 2, 2, 3, 3, 3, 2, 1)}.

    \item[Step 2] $C=261122.150$, $E=43968.050$,
    $a^{(2)}=$\texttt{(3, 4, 1, 3, 3, 1, 3, 1, 3, 3, 1, 1, 1, 1, 1, 3, 2, 2, 3, 3, 3, 2, 1)}.

    \item[Step 3] $C=266216.230$, $E=41626.130$,
    $a^{(3)}=$\texttt{(3, 4, 1, 3, 3, 1, 2, 1, 3, 3, 1, 1, 1, 1, 1, 3, 2, 2, 3, 3, 3, 2, 1)}.

    \item[Step 4] $C=267090.760$, $E=41349.520$,
    $a^{(4)}=$\texttt{(3, 4, 1, 3, 3, 1, 2, 1, 3, 2, 1, 1, 1, 1, 1, 3, 2, 2, 3, 3, 3, 2, 1)}.

    \item[Step 5] $C=268159.600$, $E=41019.460$,
    $a^{(5)}=$\texttt{(3, 4, 1, 3, 3, 1, 2, 1, 3, 2, 1, 1, 1, 1, 1, 3, 2, 1, 3, 3, 3, 2, 1)}.

    \item[Step 6] $C=272200.260$, $E=40620.210$,
    $a^{(6)}=$\texttt{(3, 4, 1, 3, 3, 1, 2, 1, 3, 2, 1, 1, 1, 1, 1, 3, 2, 1, 3, 3, 3, 2, 2)}.

    \item[Step 7] $C=275834.570$, $E=40269.950$,
    $a^{(7)}=$\texttt{(3, 4, 1, 3, 3, 1, 2, 1, 3, 2, 1, 2, 1, 1, 1, 3, 2, 1, 3, 3, 3, 2, 2)}.

    \item[Step 8] $C=279752.360$, $E=40228.530$,
    $a^{(8)}=$\texttt{(3, 4, 1, 3, 3, 1, 2, 1, 3, 2, 1, 2, 1, 1, 1, 3, 2, 1, 2, 3, 3, 2, 2)}.

\end{description}
}

\subsection*{Building type D2}

\paragraph{TOPSIS solutions for a weight sweep with $\varepsilon=0.1$.}
The results are ordered from the cost-only to the emission-only weighting.

{\footnotesize
\begin{description}
    \item[$(1,0)$] $C=255653.920$, $E=51425.990$, $S=1.000000$,
    $a^\star=$\texttt{(3, 4, 1, 2, 3, 1, 3, 1, 3, 3, 1, 1, 1, 1, 1, 3, 1, 2, 3, 3, 3, 2, 1)}.

    \item[$(0.9,0.1)$] $C=256887.670$, $E=44343.400$, $S=0.855129$,
    $a^\star=$\texttt{(3, 4, 1, 2, 3, 1, 3, 1, 3, 3, 1, 1, 1, 1, 1, 3, 2, 2, 3, 3, 3, 2, 1)}.

    \item[$(0.8,0.2)$] $C=258229.410$, $E=43069.580$, $S=0.788985$,
    $a^\star=$\texttt{(3, 4, 1, 3, 3, 1, 3, 1, 3, 3, 1, 1, 1, 1, 1, 3, 2, 2, 3, 3, 3, 2, 1)}.

    \item[$(0.7,0.3)$] $C=263285.870$, $E=40802.460$, $S=0.780528$,
    $a^\star=$\texttt{(3, 4, 1, 3, 3, 1, 2, 1, 3, 3, 1, 1, 1, 1, 1, 3, 2, 2, 3, 3, 3, 2, 1)}.

    \item[$(0.6,0.4)$] $C=263285.870$, $E=40802.460$, $S=0.824625$,
    $a^\star=$\texttt{(3, 4, 1, 3, 3, 1, 2, 1, 3, 3, 1, 1, 1, 1, 1, 3, 2, 2, 3, 3, 3, 2, 1)}.

    \item[$(0.5,0.5)$] $C=265162.900$, $E=40226.630$, $S=0.863554$,
    $a^\star=$\texttt{(3, 4, 1, 3, 3, 1, 2, 1, 3, 2, 1, 1, 1, 1, 1, 3, 2, 1, 3, 3, 3, 2, 1)}.

    \item[$(0.4,0.6)$] $C=265162.900$, $E=40226.630$, $S=0.895052$,
    $a^\star=$\texttt{(3, 4, 1, 3, 3, 1, 2, 1, 3, 2, 1, 1, 1, 1, 1, 3, 2, 1, 3, 3, 3, 2, 1)}.

    \item[$(0.3,0.7)$] $C=268395.610$, $E=39915.080$, $S=0.917624$,
    $a^\star=$\texttt{(3, 4, 1, 3, 3, 1, 2, 1, 3, 2, 1, 2, 1, 1, 1, 3, 2, 1, 3, 3, 3, 2, 1)}.

    \item[$(0.2,0.8)$] $C=269203.560$, $E=39827.380$, $S=0.944670$,
    $a^\star=$\texttt{(3, 4, 1, 3, 3, 1, 2, 1, 3, 2, 1, 1, 1, 1, 1, 3, 2, 1, 3, 3, 3, 2, 2)}.

    \item[$(0.1,0.9)$] $C=272436.270$, $E=39515.830$, $S=0.973649$,
    $a^\star=$\texttt{(3, 4, 1, 3, 3, 1, 2, 1, 3, 2, 1, 2, 1, 1, 1, 3, 2, 1, 3, 3, 3, 2, 2)}.

    \item[$(0,1)$] $C=276354.060$, $E=39474.410$, $S=1.000000$,
    $a^\star=$\texttt{(3, 4, 1, 3, 3, 1, 2, 1, 3, 2, 1, 2, 1, 1, 1, 3, 2, 1, 2, 3, 3, 2, 2)}.

\end{description}
}

\paragraph{Pareto-sequence solution.}
The successive one-group changes are:

\begin{center}
\small
\begin{tabular}{c l r r r}
\hline
Step & Change & $\Delta C$ & $\Delta E$ & $\rho=\Delta E/\Delta C$ \\
\hline
1 & $i=17:\ t_{17}=1\to 2$ & 1233.750 & 7082.590 & 5.741 \\
2 & $i=4:\ t_{4}=2\to 3$ & 1341.740 & 1273.820 & 0.949 \\
3 & $i=7:\ t_{7}=3\to 2$ & 5056.460 & 2267.120 & 0.448 \\
4 & $i=18:\ t_{18}=2\to 1$ & 1068.840 & 330.060 & 0.309 \\
5 & $i=10:\ t_{10}=3\to 2$ & 808.190 & 245.770 & 0.304 \\
6 & $i=23:\ t_{23}=1\to 2$ & 4040.660 & 399.250 & 0.099 \\
7 & $i=12:\ t_{12}=1\to 2$ & 3232.710 & 311.550 & 0.096 \\
8 & $i=19:\ t_{19}=3\to 2$ & 3917.790 & 41.420 & 0.011 \\
\hline
\end{tabular}
\end{center}

\paragraph{Complete list of Pareto-sequence solutions.}
The corresponding solution vectors are:

{\footnotesize
\begin{description}
    \item[Step 0] $C=255653.920$, $E=51425.990$,
    $a^{(0)}=$\texttt{(3, 4, 1, 2, 3, 1, 3, 1, 3, 3, 1, 1, 1, 1, 1, 3, 1, 2, 3, 3, 3, 2, 1)}.

    \item[Step 1] $C=256887.670$, $E=44343.400$,
    $a^{(1)}=$\texttt{(3, 4, 1, 2, 3, 1, 3, 1, 3, 3, 1, 1, 1, 1, 1, 3, 2, 2, 3, 3, 3, 2, 1)}.

    \item[Step 2] $C=258229.410$, $E=43069.580$,
    $a^{(2)}=$\texttt{(3, 4, 1, 3, 3, 1, 3, 1, 3, 3, 1, 1, 1, 1, 1, 3, 2, 2, 3, 3, 3, 2, 1)}.

    \item[Step 3] $C=263285.870$, $E=40802.460$,
    $a^{(3)}=$\texttt{(3, 4, 1, 3, 3, 1, 2, 1, 3, 3, 1, 1, 1, 1, 1, 3, 2, 2, 3, 3, 3, 2, 1)}.

    \item[Step 4] $C=264354.710$, $E=40472.400$,
    $a^{(4)}=$\texttt{(3, 4, 1, 3, 3, 1, 2, 1, 3, 3, 1, 1, 1, 1, 1, 3, 2, 1, 3, 3, 3, 2, 1)}.

    \item[Step 5] $C=265162.900$, $E=40226.630$,
    $a^{(5)}=$\texttt{(3, 4, 1, 3, 3, 1, 2, 1, 3, 2, 1, 1, 1, 1, 1, 3, 2, 1, 3, 3, 3, 2, 1)}.

    \item[Step 6] $C=269203.560$, $E=39827.380$,
    $a^{(6)}=$\texttt{(3, 4, 1, 3, 3, 1, 2, 1, 3, 2, 1, 1, 1, 1, 1, 3, 2, 1, 3, 3, 3, 2, 2)}.

    \item[Step 7] $C=272436.270$, $E=39515.830$,
    $a^{(7)}=$\texttt{(3, 4, 1, 3, 3, 1, 2, 1, 3, 2, 1, 2, 1, 1, 1, 3, 2, 1, 3, 3, 3, 2, 2)}.

    \item[Step 8] $C=276354.060$, $E=39474.410$,
    $a^{(8)}=$\texttt{(3, 4, 1, 3, 3, 1, 2, 1, 3, 2, 1, 2, 1, 1, 1, 3, 2, 1, 2, 3, 3, 2, 2)}.

\end{description}
}

\subsection*{Building type D3}

\paragraph{TOPSIS solutions for a weight sweep with $\varepsilon=0.1$.}
The results are ordered from the cost-only to the emission-only weighting.

{\footnotesize
\begin{description}
    \item[$(1,0)$] $C=247941.260$, $E=51040.340$, $S=1.000000$,
    $a^\star=$\texttt{(3, 4, 1, 2, 3, 1, 3, 1, 3, 3, 1, 1, 1, 1, 1, 3, 1, 2, 3, 3, 3, 2, 1)}.

    \item[$(0.9,0.1)$] $C=249175.010$, $E=43957.750$, $S=0.859253$,
    $a^\star=$\texttt{(3, 4, 1, 2, 3, 1, 3, 1, 3, 3, 1, 1, 1, 1, 1, 3, 2, 2, 3, 3, 3, 2, 1)}.

    \item[$(0.8,0.2)$] $C=250433.530$, $E=42879.310$, $S=0.787491$,
    $a^\star=$\texttt{(3, 4, 1, 3, 3, 1, 3, 1, 3, 3, 1, 1, 1, 1, 1, 3, 2, 2, 3, 3, 3, 2, 1)}.

    \item[$(0.7,0.3)$] $C=255535.780$, $E=40521.170$, $S=0.777366$,
    $a^\star=$\texttt{(3, 4, 1, 3, 3, 1, 2, 1, 3, 3, 1, 1, 1, 1, 1, 3, 2, 2, 3, 3, 3, 2, 1)}.

    \item[$(0.6,0.4)$] $C=255535.780$, $E=40521.170$, $S=0.821898$,
    $a^\star=$\texttt{(3, 4, 1, 3, 3, 1, 2, 1, 3, 3, 1, 1, 1, 1, 1, 3, 2, 2, 3, 3, 3, 2, 1)}.

    \item[$(0.5,0.5)$] $C=257417.040$, $E=39943.360$, $S=0.860805$,
    $a^\star=$\texttt{(3, 4, 1, 3, 3, 1, 2, 1, 3, 2, 1, 1, 1, 1, 1, 3, 2, 1, 3, 3, 3, 2, 1)}.

    \item[$(0.4,0.6)$] $C=257417.040$, $E=39943.360$, $S=0.893137$,
    $a^\star=$\texttt{(3, 4, 1, 3, 3, 1, 2, 1, 3, 2, 1, 1, 1, 1, 1, 3, 2, 1, 3, 3, 3, 2, 1)}.

    \item[$(0.3,0.7)$] $C=260616.840$, $E=39634.980$, $S=0.915859$,
    $a^\star=$\texttt{(3, 4, 1, 3, 3, 1, 2, 1, 3, 2, 1, 2, 1, 1, 1, 3, 2, 1, 3, 3, 3, 2, 1)}.

    \item[$(0.2,0.8)$] $C=261457.700$, $E=39544.110$, $S=0.943561$,
    $a^\star=$\texttt{(3, 4, 1, 3, 3, 1, 2, 1, 3, 2, 1, 1, 1, 1, 1, 3, 2, 1, 3, 3, 3, 2, 2)}.

    \item[$(0.1,0.9)$] $C=264657.500$, $E=39235.730$, $S=0.972950$,
    $a^\star=$\texttt{(3, 4, 1, 3, 3, 1, 2, 1, 3, 2, 1, 2, 1, 1, 1, 3, 2, 1, 3, 3, 3, 2, 2)}.

    \item[$(0,1)$] $C=268575.290$, $E=39194.310$, $S=1.000000$,
    $a^\star=$\texttt{(3, 4, 1, 3, 3, 1, 2, 1, 3, 2, 1, 2, 1, 1, 1, 3, 2, 1, 2, 3, 3, 2, 2)}.

\end{description}
}

\paragraph{Pareto-sequence solution.}
The successive one-group changes are:

\begin{center}
\small
\begin{tabular}{c l r r r}
\hline
Step & Change & $\Delta C$ & $\Delta E$ & $\rho=\Delta E/\Delta C$ \\
\hline
1 & $i=17:\ t_{17}=1\to 2$ & 1233.750 & 7082.590 & 5.741 \\
2 & $i=4:\ t_{4}=2\to 3$ & 1258.520 & 1078.440 & 0.857 \\
3 & $i=7:\ t_{7}=3\to 2$ & 5102.250 & 2358.140 & 0.462 \\
4 & $i=18:\ t_{18}=2\to 1$ & 1068.840 & 330.060 & 0.309 \\
5 & $i=10:\ t_{10}=3\to 2$ & 812.420 & 247.750 & 0.305 \\
6 & $i=23:\ t_{23}=1\to 2$ & 4040.660 & 399.250 & 0.099 \\
7 & $i=12:\ t_{12}=1\to 2$ & 3199.800 & 308.380 & 0.096 \\
8 & $i=19:\ t_{19}=3\to 2$ & 3917.790 & 41.420 & 0.011 \\
\hline
\end{tabular}
\end{center}

\paragraph{Complete list of Pareto-sequence solutions.}
The corresponding solution vectors are:

{\footnotesize
\begin{description}
    \item[Step 0] $C=247941.260$, $E=51040.340$,
    $a^{(0)}=$\texttt{(3, 4, 1, 2, 3, 1, 3, 1, 3, 3, 1, 1, 1, 1, 1, 3, 1, 2, 3, 3, 3, 2, 1)}.

    \item[Step 1] $C=249175.010$, $E=43957.750$,
    $a^{(1)}=$\texttt{(3, 4, 1, 2, 3, 1, 3, 1, 3, 3, 1, 1, 1, 1, 1, 3, 2, 2, 3, 3, 3, 2, 1)}.

    \item[Step 2] $C=250433.530$, $E=42879.310$,
    $a^{(2)}=$\texttt{(3, 4, 1, 3, 3, 1, 3, 1, 3, 3, 1, 1, 1, 1, 1, 3, 2, 2, 3, 3, 3, 2, 1)}.

    \item[Step 3] $C=255535.780$, $E=40521.170$,
    $a^{(3)}=$\texttt{(3, 4, 1, 3, 3, 1, 2, 1, 3, 3, 1, 1, 1, 1, 1, 3, 2, 2, 3, 3, 3, 2, 1)}.

    \item[Step 4] $C=256604.620$, $E=40191.110$,
    $a^{(4)}=$\texttt{(3, 4, 1, 3, 3, 1, 2, 1, 3, 3, 1, 1, 1, 1, 1, 3, 2, 1, 3, 3, 3, 2, 1)}.

    \item[Step 5] $C=257417.040$, $E=39943.360$,
    $a^{(5)}=$\texttt{(3, 4, 1, 3, 3, 1, 2, 1, 3, 2, 1, 1, 1, 1, 1, 3, 2, 1, 3, 3, 3, 2, 1)}.

    \item[Step 6] $C=261457.700$, $E=39544.110$,
    $a^{(6)}=$\texttt{(3, 4, 1, 3, 3, 1, 2, 1, 3, 2, 1, 1, 1, 1, 1, 3, 2, 1, 3, 3, 3, 2, 2)}.

    \item[Step 7] $C=264657.500$, $E=39235.730$,
    $a^{(7)}=$\texttt{(3, 4, 1, 3, 3, 1, 2, 1, 3, 2, 1, 2, 1, 1, 1, 3, 2, 1, 3, 3, 3, 2, 2)}.

    \item[Step 8] $C=268575.290$, $E=39194.310$,
    $a^{(8)}=$\texttt{(3, 4, 1, 3, 3, 1, 2, 1, 3, 2, 1, 2, 1, 1, 1, 3, 2, 1, 2, 3, 3, 2, 2)}.

\end{description}
}

\subsection*{Building type D4}

\paragraph{TOPSIS solutions for a weight sweep with $\varepsilon=0.1$.}
The results are ordered from the cost-only to the emission-only weighting.

{\footnotesize
\begin{description}
    \item[$(1,0)$] $C=246920.860$, $E=50035.080$, $S=1.000000$,
    $a^\star=$\texttt{(3, 4, 1, 2, 3, 1, 3, 1, 3, 3, 1, 1, 1, 1, 1, 3, 1, 2, 3, 3, 3, 2, 1)}.

    \item[$(0.9,0.1)$] $C=248154.610$, $E=42952.490$, $S=0.861810$,
    $a^\star=$\texttt{(3, 4, 1, 2, 3, 1, 3, 1, 3, 3, 1, 1, 1, 1, 1, 3, 2, 2, 3, 3, 3, 2, 1)}.

    \item[$(0.8,0.2)$] $C=249397.110$, $E=41911.660$, $S=0.790028$,
    $a^\star=$\texttt{(3, 4, 1, 3, 3, 1, 3, 1, 3, 3, 1, 1, 1, 1, 1, 3, 2, 2, 3, 3, 3, 2, 1)}.

    \item[$(0.7,0.3)$] $C=254448.570$, $E=39654.480$, $S=0.778757$,
    $a^\star=$\texttt{(3, 4, 1, 3, 3, 1, 2, 1, 3, 3, 1, 1, 1, 1, 1, 3, 2, 2, 3, 3, 3, 2, 1)}.

    \item[$(0.6,0.4)$] $C=254448.570$, $E=39654.480$, $S=0.820233$,
    $a^\star=$\texttt{(3, 4, 1, 3, 3, 1, 2, 1, 3, 3, 1, 1, 1, 1, 1, 3, 2, 2, 3, 3, 3, 2, 1)}.

    \item[$(0.5,0.5)$] $C=256402.560$, $E=39042.860$, $S=0.860088$,
    $a^\star=$\texttt{(3, 4, 1, 3, 3, 1, 2, 1, 3, 2, 1, 1, 1, 1, 1, 3, 2, 1, 3, 3, 3, 2, 1)}.

    \item[$(0.4,0.6)$] $C=256402.560$, $E=39042.860$, $S=0.891208$,
    $a^\star=$\texttt{(3, 4, 1, 3, 3, 1, 2, 1, 3, 2, 1, 1, 1, 1, 1, 3, 2, 1, 3, 3, 3, 2, 1)}.

    \item[$(0.3,0.7)$] $C=260166.900$, $E=38680.070$, $S=0.913900$,
    $a^\star=$\texttt{(3, 4, 1, 3, 3, 1, 2, 1, 3, 2, 1, 2, 1, 1, 1, 3, 2, 1, 3, 3, 3, 2, 1)}.

    \item[$(0.2,0.8)$] $C=260443.220$, $E=38643.610$, $S=0.941498$,
    $a^\star=$\texttt{(3, 4, 1, 3, 3, 1, 2, 1, 3, 2, 1, 1, 1, 1, 1, 3, 2, 1, 3, 3, 3, 2, 2)}.

    \item[$(0.1,0.9)$] $C=264207.560$, $E=38280.820$, $S=0.972456$,
    $a^\star=$\texttt{(3, 4, 1, 3, 3, 1, 2, 1, 3, 2, 1, 2, 1, 1, 1, 3, 2, 1, 3, 3, 3, 2, 2)}.

    \item[$(0,1)$] $C=268125.350$, $E=38239.400$, $S=1.000000$,
    $a^\star=$\texttt{(3, 4, 1, 3, 3, 1, 2, 1, 3, 2, 1, 2, 1, 1, 1, 3, 2, 1, 2, 3, 3, 2, 2)}.

\end{description}
}

\paragraph{Pareto-sequence solution.}
The successive one-group changes are:

\begin{center}
\small
\begin{tabular}{c l r r r}
\hline
Step & Change & $\Delta C$ & $\Delta E$ & $\rho=\Delta E/\Delta C$ \\
\hline
1 & $i=17:\ t_{17}=1\to 2$ & 1233.750 & 7082.590 & 5.741 \\
2 & $i=4:\ t_{4}=2\to 3$ & 1242.500 & 1040.830 & 0.838 \\
3 & $i=7:\ t_{7}=3\to 2$ & 5051.460 & 2257.180 & 0.447 \\
4 & $i=10:\ t_{10}=3\to 2$ & 885.150 & 281.560 & 0.318 \\
5 & $i=18:\ t_{18}=2\to 1$ & 1068.840 & 330.060 & 0.309 \\
6 & $i=23:\ t_{23}=1\to 2$ & 4040.660 & 399.250 & 0.099 \\
7 & $i=12:\ t_{12}=1\to 2$ & 3764.340 & 362.790 & 0.096 \\
8 & $i=19:\ t_{19}=3\to 2$ & 3917.790 & 41.420 & 0.011 \\
\hline
\end{tabular}
\end{center}

\paragraph{Complete list of Pareto-sequence solutions.}
The corresponding solution vectors are:

{\footnotesize
\begin{description}
    \item[Step 0] $C=246920.860$, $E=50035.080$,
    $a^{(0)}=$\texttt{(3, 4, 1, 2, 3, 1, 3, 1, 3, 3, 1, 1, 1, 1, 1, 3, 1, 2, 3, 3, 3, 2, 1)}.

    \item[Step 1] $C=248154.610$, $E=42952.490$,
    $a^{(1)}=$\texttt{(3, 4, 1, 2, 3, 1, 3, 1, 3, 3, 1, 1, 1, 1, 1, 3, 2, 2, 3, 3, 3, 2, 1)}.

    \item[Step 2] $C=249397.110$, $E=41911.660$,
    $a^{(2)}=$\texttt{(3, 4, 1, 3, 3, 1, 3, 1, 3, 3, 1, 1, 1, 1, 1, 3, 2, 2, 3, 3, 3, 2, 1)}.

    \item[Step 3] $C=254448.570$, $E=39654.480$,
    $a^{(3)}=$\texttt{(3, 4, 1, 3, 3, 1, 2, 1, 3, 3, 1, 1, 1, 1, 1, 3, 2, 2, 3, 3, 3, 2, 1)}.

    \item[Step 4] $C=255333.720$, $E=39372.920$,
    $a^{(4)}=$\texttt{(3, 4, 1, 3, 3, 1, 2, 1, 3, 2, 1, 1, 1, 1, 1, 3, 2, 2, 3, 3, 3, 2, 1)}.

    \item[Step 5] $C=256402.560$, $E=39042.860$,
    $a^{(5)}=$\texttt{(3, 4, 1, 3, 3, 1, 2, 1, 3, 2, 1, 1, 1, 1, 1, 3, 2, 1, 3, 3, 3, 2, 1)}.

    \item[Step 6] $C=260443.220$, $E=38643.610$,
    $a^{(6)}=$\texttt{(3, 4, 1, 3, 3, 1, 2, 1, 3, 2, 1, 1, 1, 1, 1, 3, 2, 1, 3, 3, 3, 2, 2)}.

    \item[Step 7] $C=264207.560$, $E=38280.820$,
    $a^{(7)}=$\texttt{(3, 4, 1, 3, 3, 1, 2, 1, 3, 2, 1, 2, 1, 1, 1, 3, 2, 1, 3, 3, 3, 2, 2)}.

    \item[Step 8] $C=268125.350$, $E=38239.400$,
    $a^{(8)}=$\texttt{(3, 4, 1, 3, 3, 1, 2, 1, 3, 2, 1, 2, 1, 1, 1, 3, 2, 1, 2, 3, 3, 2, 2)}.

\end{description}
}

\subsection*{Building type K1}

\paragraph{TOPSIS solutions for a weight sweep with $\varepsilon=0.1$.}
The results are ordered from the cost-only to the emission-only weighting.

{\footnotesize
\begin{description}
    \item[$(1,0)$] $C=110751.520$, $E=23772.660$, $S=1.000000$,
    $a^\star=$\texttt{(3, 4, 1, 2, 3, 1, 3, 1, 3, 3, 1, 1, 1, 1, 1, 3, 1, 1, 1, 3, 1, 2, 1)}.

    \item[$(0.9,0.1)$] $C=111960.980$, $E=19831.460$, $S=0.848929$,
    $a^\star=$\texttt{(3, 4, 1, 2, 3, 1, 3, 1, 3, 3, 1, 1, 1, 1, 1, 3, 2, 1, 1, 3, 3, 2, 1)}.

    \item[$(0.8,0.2)$] $C=111960.980$, $E=19831.460$, $S=0.804583$,
    $a^\star=$\texttt{(3, 4, 1, 2, 3, 1, 3, 1, 3, 3, 1, 1, 1, 1, 1, 3, 2, 1, 1, 3, 3, 2, 1)}.

    \item[$(0.7,0.3)$] $C=112969.870$, $E=19339.150$, $S=0.796031$,
    $a^\star=$\texttt{(3, 4, 1, 3, 3, 1, 3, 1, 3, 3, 1, 1, 1, 1, 1, 3, 2, 1, 1, 3, 3, 2, 1)}.

    \item[$(0.6,0.4)$] $C=112969.870$, $E=19339.150$, $S=0.805229$,
    $a^\star=$\texttt{(3, 4, 1, 3, 3, 1, 3, 1, 3, 3, 1, 1, 1, 1, 1, 3, 2, 1, 1, 3, 3, 2, 1)}.

    \item[$(0.5,0.5)$] $C=114349.260$, $E=19187.640$, $S=0.813922$,
    $a^\star=$\texttt{(3, 4, 1, 3, 3, 1, 3, 1, 3, 3, 1, 1, 1, 1, 1, 3, 2, 1, 1, 3, 3, 2, 2)}.

    \item[$(0.4,0.6)$] $C=117241.340$, $E=18632.220$, $S=0.859112$,
    $a^\star=$\texttt{(3, 4, 1, 3, 3, 1, 2, 1, 3, 3, 1, 1, 1, 1, 1, 3, 2, 1, 1, 3, 3, 2, 1)}.

    \item[$(0.3,0.7)$] $C=117241.340$, $E=18632.220$, $S=0.897387$,
    $a^\star=$\texttt{(3, 4, 1, 3, 3, 1, 2, 1, 3, 3, 1, 1, 1, 1, 1, 3, 2, 1, 1, 3, 3, 2, 1)}.

    \item[$(0.2,0.8)$] $C=118620.730$, $E=18480.710$, $S=0.931760$,
    $a^\star=$\texttt{(3, 4, 1, 3, 3, 1, 2, 1, 3, 3, 1, 1, 1, 1, 1, 3, 2, 1, 1, 3, 3, 2, 2)}.

    \item[$(0.1,0.9)$] $C=119765.180$, $E=18369.170$, $S=0.966713$,
    $a^\star=$\texttt{(3, 4, 1, 3, 3, 1, 2, 1, 3, 3, 1, 2, 1, 1, 1, 3, 2, 1, 1, 3, 3, 2, 2)}.

    \item[$(0,1)$] $C=120113.620$, $E=18337.050$, $S=1.000000$,
    $a^\star=$\texttt{(3, 4, 1, 3, 3, 1, 2, 1, 3, 2, 1, 2, 1, 1, 1, 3, 2, 1, 1, 3, 3, 2, 2)}.

\end{description}
}

\paragraph{Pareto-sequence solution.}
The successive one-group changes are:

\begin{center}
\small
\begin{tabular}{c l r r r}
\hline
Step & Change & $\Delta C$ & $\Delta E$ & $\rho=\Delta E/\Delta C$ \\
\hline
1 & $i=21:\ t_{21}=1\to 3$ & 99.930 & 601.220 & 6.016 \\
2 & $i=17:\ t_{17}=1\to 2$ & 1109.530 & 3339.980 & 3.010 \\
3 & $i=4:\ t_{4}=2\to 3$ & 1008.890 & 492.310 & 0.488 \\
4 & $i=7:\ t_{7}=3\to 2$ & 4271.470 & 706.930 & 0.166 \\
5 & $i=23:\ t_{23}=1\to 2$ & 1379.390 & 151.510 & 0.110 \\
6 & $i=12:\ t_{12}=1\to 2$ & 1144.450 & 111.540 & 0.097 \\
7 & $i=10:\ t_{10}=3\to 2$ & 348.440 & 32.120 & 0.092 \\
\hline
\end{tabular}
\end{center}

\paragraph{Complete list of Pareto-sequence solutions.}
The corresponding solution vectors are:

{\footnotesize
\begin{description}
    \item[Step 0] $C=110751.520$, $E=23772.660$,
    $a^{(0)}=$\texttt{(3, 4, 1, 2, 3, 1, 3, 1, 3, 3, 1, 1, 1, 1, 1, 3, 1, 1, 1, 3, 1, 2, 1)}.

    \item[Step 1] $C=110851.450$, $E=23171.440$,
    $a^{(1)}=$\texttt{(3, 4, 1, 2, 3, 1, 3, 1, 3, 3, 1, 1, 1, 1, 1, 3, 1, 1, 1, 3, 3, 2, 1)}.

    \item[Step 2] $C=111960.980$, $E=19831.460$,
    $a^{(2)}=$\texttt{(3, 4, 1, 2, 3, 1, 3, 1, 3, 3, 1, 1, 1, 1, 1, 3, 2, 1, 1, 3, 3, 2, 1)}.

    \item[Step 3] $C=112969.870$, $E=19339.150$,
    $a^{(3)}=$\texttt{(3, 4, 1, 3, 3, 1, 3, 1, 3, 3, 1, 1, 1, 1, 1, 3, 2, 1, 1, 3, 3, 2, 1)}.

    \item[Step 4] $C=117241.340$, $E=18632.220$,
    $a^{(4)}=$\texttt{(3, 4, 1, 3, 3, 1, 2, 1, 3, 3, 1, 1, 1, 1, 1, 3, 2, 1, 1, 3, 3, 2, 1)}.

    \item[Step 5] $C=118620.730$, $E=18480.710$,
    $a^{(5)}=$\texttt{(3, 4, 1, 3, 3, 1, 2, 1, 3, 3, 1, 1, 1, 1, 1, 3, 2, 1, 1, 3, 3, 2, 2)}.

    \item[Step 6] $C=119765.180$, $E=18369.170$,
    $a^{(6)}=$\texttt{(3, 4, 1, 3, 3, 1, 2, 1, 3, 3, 1, 2, 1, 1, 1, 3, 2, 1, 1, 3, 3, 2, 2)}.

    \item[Step 7] $C=120113.620$, $E=18337.050$,
    $a^{(7)}=$\texttt{(3, 4, 1, 3, 3, 1, 2, 1, 3, 2, 1, 2, 1, 1, 1, 3, 2, 1, 1, 3, 3, 2, 2)}.

\end{description}
}

\subsection*{Building type K2}

\paragraph{TOPSIS solutions for a weight sweep with $\varepsilon=0.1$.}
The results are ordered from the cost-only to the emission-only weighting.

{\footnotesize
\begin{description}
    \item[$(1,0)$] $C=168060.300$, $E=34466.200$, $S=1.000000$,
    $a^\star=$\texttt{(3, 4, 1, 2, 3, 1, 3, 1, 3, 3, 1, 1, 1, 1, 1, 3, 1, 1, 3, 3, 3, 2, 1)}.

    \item[$(0.9,0.1)$] $C=169169.800$, $E=31126.300$, $S=0.818749$,
    $a^\star=$\texttt{(3, 4, 1, 2, 3, 1, 3, 1, 3, 3, 1, 1, 1, 1, 1, 3, 2, 1, 3, 3, 3, 2, 1)}.

    \item[$(0.8,0.2)$] $C=169945.800$, $E=29890.100$, $S=0.767560$,
    $a^\star=$\texttt{(3, 4, 1, 3, 3, 1, 3, 1, 3, 3, 1, 1, 1, 1, 1, 3, 2, 1, 3, 3, 3, 2, 1)}.

    \item[$(0.7,0.3)$] $C=169945.800$, $E=29890.100$, $S=0.729521$,
    $a^\star=$\texttt{(3, 4, 1, 3, 3, 1, 3, 1, 3, 3, 1, 1, 1, 1, 1, 3, 2, 1, 3, 3, 3, 2, 1)}.

    \item[$(0.6,0.4)$] $C=174650.700$, $E=28321.700$, $S=0.771960$,
    $a^\star=$\texttt{(3, 4, 1, 3, 3, 1, 2, 1, 3, 3, 1, 1, 1, 1, 1, 3, 2, 1, 3, 3, 3, 2, 1)}.

    \item[$(0.5,0.5)$] $C=174650.700$, $E=28321.700$, $S=0.826239$,
    $a^\star=$\texttt{(3, 4, 1, 3, 3, 1, 2, 1, 3, 3, 1, 1, 1, 1, 1, 3, 2, 1, 3, 3, 3, 2, 1)}.

    \item[$(0.4,0.6)$] $C=175170.800$, $E=28209.800$, $S=0.867538$,
    $a^\star=$\texttt{(3, 4, 1, 3, 3, 1, 2, 1, 3, 2, 1, 1, 1, 1, 1, 3, 2, 1, 3, 3, 3, 2, 1)}.

    \item[$(0.3,0.7)$] $C=175170.800$, $E=28209.800$, $S=0.899029$,
    $a^\star=$\texttt{(3, 4, 1, 3, 3, 1, 2, 1, 3, 2, 1, 1, 1, 1, 1, 3, 2, 1, 3, 3, 3, 2, 1)}.

    \item[$(0.2,0.8)$] $C=177880.800$, $E=27934.400$, $S=0.932717$,
    $a^\star=$\texttt{(3, 4, 1, 3, 3, 1, 2, 1, 3, 2, 1, 1, 1, 1, 1, 3, 2, 1, 3, 3, 3, 2, 2)}.

    \item[$(0.1,0.9)$] $C=179623.400$, $E=27766.400$, $S=0.966756$,
    $a^\star=$\texttt{(3, 4, 1, 3, 3, 1, 2, 1, 3, 2, 1, 2, 1, 1, 1, 3, 2, 1, 3, 3, 3, 2, 2)}.

    \item[$(0,1)$] $C=179623.400$, $E=27766.400$, $S=1.000000$,
    $a^\star=$\texttt{(3, 4, 1, 3, 3, 1, 2, 1, 3, 2, 1, 2, 1, 1, 1, 3, 2, 1, 3, 3, 3, 2, 2)}.

\end{description}
}

\paragraph{Pareto-sequence solution.}
The successive one-group changes are:

\begin{center}
\small
\begin{tabular}{c l r r r}
\hline
Step & Change & $\Delta C$ & $\Delta E$ & $\rho=\Delta E/\Delta C$ \\
\hline
1 & $i=17:\ t_{17}=1\to 2$ & 1109.500 & 3339.900 & 3.010 \\
2 & $i=4:\ t_{4}=2\to 3$ & 776.000 & 1236.200 & 1.593 \\
3 & $i=7:\ t_{7}=3\to 2$ & 4704.900 & 1568.400 & 0.333 \\
4 & $i=10:\ t_{10}=3\to 2$ & 520.100 & 111.900 & 0.215 \\
5 & $i=23:\ t_{23}=1\to 2$ & 2710.000 & 275.400 & 0.102 \\
6 & $i=12:\ t_{12}=1\to 2$ & 1742.600 & 168.000 & 0.096 \\
\hline
\end{tabular}
\end{center}

\paragraph{Complete list of Pareto-sequence solutions.}
The corresponding solution vectors are:

{\footnotesize
\begin{description}
    \item[Step 0] $C=168060.300$, $E=34466.200$,
    $a^{(0)}=$\texttt{(3, 4, 1, 2, 3, 1, 3, 1, 3, 3, 1, 1, 1, 1, 1, 3, 1, 1, 3, 3, 3, 2, 1)}.

    \item[Step 1] $C=169169.800$, $E=31126.300$,
    $a^{(1)}=$\texttt{(3, 4, 1, 2, 3, 1, 3, 1, 3, 3, 1, 1, 1, 1, 1, 3, 2, 1, 3, 3, 3, 2, 1)}.

    \item[Step 2] $C=169945.800$, $E=29890.100$,
    $a^{(2)}=$\texttt{(3, 4, 1, 3, 3, 1, 3, 1, 3, 3, 1, 1, 1, 1, 1, 3, 2, 1, 3, 3, 3, 2, 1)}.

    \item[Step 3] $C=174650.700$, $E=28321.700$,
    $a^{(3)}=$\texttt{(3, 4, 1, 3, 3, 1, 2, 1, 3, 3, 1, 1, 1, 1, 1, 3, 2, 1, 3, 3, 3, 2, 1)}.

    \item[Step 4] $C=175170.800$, $E=28209.800$,
    $a^{(4)}=$\texttt{(3, 4, 1, 3, 3, 1, 2, 1, 3, 2, 1, 1, 1, 1, 1, 3, 2, 1, 3, 3, 3, 2, 1)}.

    \item[Step 5] $C=177880.800$, $E=27934.400$,
    $a^{(5)}=$\texttt{(3, 4, 1, 3, 3, 1, 2, 1, 3, 2, 1, 1, 1, 1, 1, 3, 2, 1, 3, 3, 3, 2, 2)}.

    \item[Step 6] $C=179623.400$, $E=27766.400$,
    $a^{(6)}=$\texttt{(3, 4, 1, 3, 3, 1, 2, 1, 3, 2, 1, 2, 1, 1, 1, 3, 2, 1, 3, 3, 3, 2, 2)}.

\end{description}
}

\section{Numerical results for the M scenario (materials only)}\label{app:M}

This appendix reports, for the building types A1, B1, C1, D1--D4, K1 and K2, the TOPSIS weight-sweep solutions and the Pareto sequences obtained with the procedures presented in Section~\ref{sec:numerical} from the data of the M scenario (materials only). The notation, the dominance filtering and the labelling of the types are the same as in Section~\ref{sec:numerical}. The Pareto sequences are those shown in Figure~\ref{fig:Fig 4}.

\subsection*{Building type A1}

\paragraph{TOPSIS solutions for a weight sweep with $\varepsilon=0.1$.}
The results are ordered from the cost-only to the emission-only weighting.

{\footnotesize
\begin{description}
    \item[$(1,0)$] $C=102727.758$, $E=17174.929$, $S=1.000000$,
    $a^\star=$\texttt{(3, 4, 1, 1, 1, 1, 3, 1, 3, 3, 1, 1, 1, 1, 1, 3, 1, 2, 3, 3, 3, 2, 1)}.

    \item[$(0.9,0.1)$] $C=103496.601$, $E=15400.244$, $S=0.855152$,
    $a^\star=$\texttt{(3, 4, 3, 1, 1, 1, 2, 1, 3, 3, 1, 1, 1, 1, 1, 3, 1, 2, 2, 3, 3, 2, 1)}.

    \item[$(0.8,0.2)$] $C=103863.801$, $E=15189.799$, $S=0.793577$,
    $a^\star=$\texttt{(3, 4, 3, 1, 1, 1, 2, 1, 3, 2, 1, 1, 1, 1, 1, 3, 1, 2, 2, 3, 3, 2, 1)}.

    \item[$(0.7,0.3)$] $C=103943.421$, $E=15159.159$, $S=0.765892$,
    $a^\star=$\texttt{(3, 4, 3, 1, 1, 1, 2, 1, 3, 2, 1, 1, 2, 1, 1, 3, 1, 2, 2, 3, 3, 2, 1)}.

    \item[$(0.6,0.4)$] $C=103943.421$, $E=15159.159$, $S=0.748272$,
    $a^\star=$\texttt{(3, 4, 3, 1, 1, 1, 2, 1, 3, 2, 1, 1, 2, 1, 1, 3, 1, 2, 2, 3, 3, 2, 1)}.

    \item[$(0.5,0.5)$] $C=106334.910$, $E=14719.159$, $S=0.800871$,
    $a^\star=$\texttt{(3, 4, 3, 1, 1, 1, 2, 1, 3, 2, 1, 1, 2, 1, 1, 3, 1, 1, 2, 3, 3, 2, 1)}.

    \item[$(0.4,0.6)$] $C=106334.910$, $E=14719.159$, $S=0.838456$,
    $a^\star=$\texttt{(3, 4, 3, 1, 1, 1, 2, 1, 3, 2, 1, 1, 2, 1, 1, 3, 1, 1, 2, 3, 3, 2, 1)}.

    \item[$(0.3,0.7)$] $C=107669.211$, $E=14595.441$, $S=0.874923$,
    $a^\star=$\texttt{(3, 4, 3, 1, 1, 1, 2, 1, 3, 2, 1, 1, 2, 1, 1, 3, 1, 1, 2, 3, 3, 2, 2)}.

    \item[$(0.2,0.8)$] $C=108339.176$, $E=14526.945$, $S=0.916049$,
    $a^\star=$\texttt{(3, 4, 3, 1, 1, 1, 2, 1, 3, 2, 1, 2, 2, 1, 1, 3, 1, 1, 2, 3, 3, 2, 1)}.

    \item[$(0.1,0.9)$] $C=109673.478$, $E=14403.226$, $S=0.960204$,
    $a^\star=$\texttt{(3, 4, 3, 1, 1, 1, 2, 1, 3, 2, 1, 2, 2, 1, 1, 3, 1, 1, 2, 3, 3, 2, 2)}.

    \item[$(0,1)$] $C=109673.478$, $E=14403.226$, $S=1.000000$,
    $a^\star=$\texttt{(3, 4, 3, 1, 1, 1, 2, 1, 3, 2, 1, 2, 2, 1, 1, 3, 1, 1, 2, 3, 3, 2, 2)}.

\end{description}
}

\paragraph{Pareto-sequence solution.}
The successive one-group changes are:

\begin{center}
\small
\begin{tabular}{c l r r r}
\hline
Step & Change & $\Delta C$ & $\Delta E$ & $\rho=\Delta E/\Delta C$ \\
\hline
1 & $i=3:\ t_{3}=1\to 3$ & 10.085 & 43.295 & 4.293 \\
2 & $i=19:\ t_{19}=3\to 2$ & 83.689 & 203.596 & 2.433 \\
3 & $i=7:\ t_{7}=3\to 2$ & 675.068 & 1527.793 & 2.263 \\
4 & $i=10:\ t_{10}=3\to 2$ & 367.201 & 210.445 & 0.573 \\
5 & $i=13:\ t_{13}=1\to 2$ & 79.620 & 30.640 & 0.385 \\
6 & $i=18:\ t_{18}=2\to 1$ & 2391.489 & 440.000 & 0.184 \\
7 & $i=12:\ t_{12}=1\to 2$ & 2004.266 & 192.214 & 0.096 \\
8 & $i=23:\ t_{23}=1\to 2$ & 1334.301 & 123.718 & 0.093 \\
\hline
\end{tabular}
\end{center}

\paragraph{Complete list of Pareto-sequence solutions.}
The corresponding solution vectors are:

{\footnotesize
\begin{description}
    \item[Step 0] $C=102727.758$, $E=17174.929$,
    $a^{(0)}=$\texttt{(3, 4, 1, 1, 1, 1, 3, 1, 3, 3, 1, 1, 1, 1, 1, 3, 1, 2, 3, 3, 3, 2, 1)}.

    \item[Step 1] $C=102737.844$, $E=17131.633$,
    $a^{(1)}=$\texttt{(3, 4, 3, 1, 1, 1, 3, 1, 3, 3, 1, 1, 1, 1, 1, 3, 1, 2, 3, 3, 3, 2, 1)}.

    \item[Step 2] $C=102821.532$, $E=16928.037$,
    $a^{(2)}=$\texttt{(3, 4, 3, 1, 1, 1, 3, 1, 3, 3, 1, 1, 1, 1, 1, 3, 1, 2, 2, 3, 3, 2, 1)}.

    \item[Step 3] $C=103496.601$, $E=15400.244$,
    $a^{(3)}=$\texttt{(3, 4, 3, 1, 1, 1, 2, 1, 3, 3, 1, 1, 1, 1, 1, 3, 1, 2, 2, 3, 3, 2, 1)}.

    \item[Step 4] $C=103863.801$, $E=15189.799$,
    $a^{(4)}=$\texttt{(3, 4, 3, 1, 1, 1, 2, 1, 3, 2, 1, 1, 1, 1, 1, 3, 1, 2, 2, 3, 3, 2, 1)}.

    \item[Step 5] $C=103943.421$, $E=15159.159$,
    $a^{(5)}=$\texttt{(3, 4, 3, 1, 1, 1, 2, 1, 3, 2, 1, 1, 2, 1, 1, 3, 1, 2, 2, 3, 3, 2, 1)}.

    \item[Step 6] $C=106334.910$, $E=14719.159$,
    $a^{(6)}=$\texttt{(3, 4, 3, 1, 1, 1, 2, 1, 3, 2, 1, 1, 2, 1, 1, 3, 1, 1, 2, 3, 3, 2, 1)}.

    \item[Step 7] $C=108339.176$, $E=14526.945$,
    $a^{(7)}=$\texttt{(3, 4, 3, 1, 1, 1, 2, 1, 3, 2, 1, 2, 2, 1, 1, 3, 1, 1, 2, 3, 3, 2, 1)}.

    \item[Step 8] $C=109673.478$, $E=14403.226$,
    $a^{(8)}=$\texttt{(3, 4, 3, 1, 1, 1, 2, 1, 3, 2, 1, 2, 2, 1, 1, 3, 1, 1, 2, 3, 3, 2, 2)}.

\end{description}
}

\subsection*{Building type B1}

\paragraph{TOPSIS solutions for a weight sweep with $\varepsilon=0.1$.}
The results are ordered from the cost-only to the emission-only weighting.

{\footnotesize
\begin{description}
    \item[$(1,0)$] $C=185356.274$, $E=30988.746$, $S=1.000000$,
    $a^\star=$\texttt{(3, 4, 1, 2, 3, 1, 3, 1, 3, 3, 1, 1, 1, 1, 1, 3, 3, 2, 3, 3, 3, 2, 1)}.

    \item[$(0.9,0.1)$] $C=186968.567$, $E=27085.583$, $S=0.820232$,
    $a^\star=$\texttt{(3, 4, 3, 3, 3, 1, 2, 1, 3, 3, 1, 1, 1, 1, 1, 3, 3, 2, 2, 3, 3, 2, 1)}.

    \item[$(0.8,0.2)$] $C=187493.127$, $E=26784.954$, $S=0.780293$,
    $a^\star=$\texttt{(3, 4, 3, 3, 3, 1, 2, 1, 3, 2, 1, 1, 1, 1, 1, 3, 3, 2, 2, 3, 3, 2, 1)}.

    \item[$(0.7,0.3)$] $C=188817.699$, $E=26284.598$, $S=0.786377$,
    $a^\star=$\texttt{(3, 4, 3, 3, 3, 1, 2, 1, 3, 2, 1, 1, 1, 1, 1, 3, 2, 2, 2, 3, 3, 2, 1)}.

    \item[$(0.6,0.4)$] $C=188917.899$, $E=26246.038$, $S=0.821036$,
    $a^\star=$\texttt{(3, 4, 3, 3, 3, 1, 2, 1, 3, 2, 1, 1, 2, 1, 1, 3, 2, 2, 2, 3, 3, 2, 1)}.

    \item[$(0.5,0.5)$] $C=190113.643$, $E=26026.038$, $S=0.851653$,
    $a^\star=$\texttt{(3, 4, 3, 3, 3, 1, 2, 1, 3, 2, 1, 1, 2, 1, 1, 3, 2, 1, 2, 3, 3, 2, 1)}.

    \item[$(0.4,0.6)$] $C=190113.643$, $E=26026.038$, $S=0.876992$,
    $a^\star=$\texttt{(3, 4, 3, 3, 3, 1, 2, 1, 3, 2, 1, 1, 2, 1, 1, 3, 2, 1, 2, 3, 3, 2, 1)}.

    \item[$(0.3,0.7)$] $C=192938.138$, $E=25755.162$, $S=0.906179$,
    $a^\star=$\texttt{(3, 4, 3, 3, 3, 1, 2, 1, 3, 2, 1, 2, 2, 1, 1, 3, 2, 1, 2, 3, 3, 2, 1)}.

    \item[$(0.2,0.8)$] $C=195606.741$, $E=25507.726$, $S=0.935096$,
    $a^\star=$\texttt{(3, 4, 3, 3, 3, 1, 2, 1, 3, 2, 1, 2, 2, 1, 1, 3, 2, 1, 2, 3, 3, 2, 2)}.

    \item[$(0.1,0.9)$] $C=195606.741$, $E=25507.726$, $S=0.970075$,
    $a^\star=$\texttt{(3, 4, 3, 3, 3, 1, 2, 1, 3, 2, 1, 2, 2, 1, 1, 3, 2, 1, 2, 3, 3, 2, 2)}.

    \item[$(0,1)$] $C=195606.741$, $E=25507.726$, $S=1.000000$,
    $a^\star=$\texttt{(3, 4, 3, 3, 3, 1, 2, 1, 3, 2, 1, 2, 2, 1, 1, 3, 2, 1, 2, 3, 3, 2, 2)}.

\end{description}
}

\paragraph{Pareto-sequence solution.}
The successive one-group changes are:

\begin{center}
\small
\begin{tabular}{c l r r r}
\hline
Step & Change & $\Delta C$ & $\Delta E$ & $\rho=\Delta E/\Delta C$ \\
\hline
1 & $i=3:\ t_{3}=1\to 3$ & 18.973 & 81.448 & 4.293 \\
2 & $i=4:\ t_{4}=2\to 3$ & 436.358 & 1174.925 & 2.693 \\
3 & $i=19:\ t_{19}=3\to 2$ & 167.377 & 407.192 & 2.433 \\
4 & $i=7:\ t_{7}=3\to 2$ & 989.585 & 2239.598 & 2.263 \\
5 & $i=10:\ t_{10}=3\to 2$ & 524.560 & 300.628 & 0.573 \\
6 & $i=13:\ t_{13}=1\to 2$ & 100.200 & 38.560 & 0.385 \\
7 & $i=17:\ t_{17}=3\to 2$ & 1324.571 & 500.356 & 0.378 \\
8 & $i=18:\ t_{18}=2\to 1$ & 1195.745 & 220.000 & 0.184 \\
9 & $i=12:\ t_{12}=1\to 2$ & 2824.495 & 270.876 & 0.096 \\
10 & $i=23:\ t_{23}=1\to 2$ & 2668.603 & 247.436 & 0.093 \\
\hline
\end{tabular}
\end{center}

\paragraph{Complete list of Pareto-sequence solutions.}
The corresponding solution vectors are:

{\footnotesize
\begin{description}
    \item[Step 0] $C=185356.274$, $E=30988.746$,
    $a^{(0)}=$\texttt{(3, 4, 1, 2, 3, 1, 3, 1, 3, 3, 1, 1, 1, 1, 1, 3, 3, 2, 3, 3, 3, 2, 1)}.

    \item[Step 1] $C=185375.247$, $E=30907.298$,
    $a^{(1)}=$\texttt{(3, 4, 3, 2, 3, 1, 3, 1, 3, 3, 1, 1, 1, 1, 1, 3, 3, 2, 3, 3, 3, 2, 1)}.

    \item[Step 2] $C=185811.605$, $E=29732.373$,
    $a^{(2)}=$\texttt{(3, 4, 3, 3, 3, 1, 3, 1, 3, 3, 1, 1, 1, 1, 1, 3, 3, 2, 3, 3, 3, 2, 1)}.

    \item[Step 3] $C=185978.982$, $E=29325.181$,
    $a^{(3)}=$\texttt{(3, 4, 3, 3, 3, 1, 3, 1, 3, 3, 1, 1, 1, 1, 1, 3, 3, 2, 2, 3, 3, 2, 1)}.

    \item[Step 4] $C=186968.567$, $E=27085.583$,
    $a^{(4)}=$\texttt{(3, 4, 3, 3, 3, 1, 2, 1, 3, 3, 1, 1, 1, 1, 1, 3, 3, 2, 2, 3, 3, 2, 1)}.

    \item[Step 5] $C=187493.127$, $E=26784.954$,
    $a^{(5)}=$\texttt{(3, 4, 3, 3, 3, 1, 2, 1, 3, 2, 1, 1, 1, 1, 1, 3, 3, 2, 2, 3, 3, 2, 1)}.

    \item[Step 6] $C=187593.327$, $E=26746.394$,
    $a^{(6)}=$\texttt{(3, 4, 3, 3, 3, 1, 2, 1, 3, 2, 1, 1, 2, 1, 1, 3, 3, 2, 2, 3, 3, 2, 1)}.

    \item[Step 7] $C=188917.899$, $E=26246.038$,
    $a^{(7)}=$\texttt{(3, 4, 3, 3, 3, 1, 2, 1, 3, 2, 1, 1, 2, 1, 1, 3, 2, 2, 2, 3, 3, 2, 1)}.

    \item[Step 8] $C=190113.643$, $E=26026.038$,
    $a^{(8)}=$\texttt{(3, 4, 3, 3, 3, 1, 2, 1, 3, 2, 1, 1, 2, 1, 1, 3, 2, 1, 2, 3, 3, 2, 1)}.

    \item[Step 9] $C=192938.138$, $E=25755.162$,
    $a^{(9)}=$\texttt{(3, 4, 3, 3, 3, 1, 2, 1, 3, 2, 1, 2, 2, 1, 1, 3, 2, 1, 2, 3, 3, 2, 1)}.

    \item[Step 10] $C=195606.741$, $E=25507.726$,
    $a^{(10)}=$\texttt{(3, 4, 3, 3, 3, 1, 2, 1, 3, 2, 1, 2, 2, 1, 1, 3, 2, 1, 2, 3, 3, 2, 2)}.

\end{description}
}

\subsection*{Building type C1}

\paragraph{TOPSIS solutions for a weight sweep with $\varepsilon=0.1$.}
The results are ordered from the cost-only to the emission-only weighting.

{\footnotesize
\begin{description}
    \item[$(1,0)$] $C=216363.233$, $E=36252.943$, $S=1.000000$,
    $a^\star=$\texttt{(3, 4, 1, 2, 3, 1, 3, 1, 3, 3, 1, 1, 1, 1, 1, 3, 3, 2, 3, 3, 3, 2, 1)}.

    \item[$(0.9,0.1)$] $C=218267.392$, $E=31635.628$, $S=0.817049$,
    $a^\star=$\texttt{(3, 4, 3, 3, 3, 1, 2, 1, 3, 3, 1, 1, 1, 1, 1, 3, 3, 2, 2, 3, 3, 2, 1)}.

    \item[$(0.8,0.2)$] $C=218850.758$, $E=31301.297$, $S=0.783934$,
    $a^\star=$\texttt{(3, 4, 3, 3, 3, 1, 2, 1, 3, 2, 1, 1, 1, 1, 1, 3, 3, 2, 2, 3, 3, 2, 1)}.

    \item[$(0.7,0.3)$] $C=220175.330$, $E=30800.941$, $S=0.793196$,
    $a^\star=$\texttt{(3, 4, 3, 3, 3, 1, 2, 1, 3, 2, 1, 1, 1, 1, 1, 3, 2, 2, 2, 3, 3, 2, 1)}.

    \item[$(0.6,0.4)$] $C=220285.820$, $E=30758.421$, $S=0.826968$,
    $a^\star=$\texttt{(3, 4, 3, 3, 3, 1, 2, 1, 3, 2, 1, 1, 2, 1, 1, 3, 2, 2, 2, 3, 3, 2, 1)}.

    \item[$(0.5,0.5)$] $C=221481.565$, $E=30538.421$, $S=0.856196$,
    $a^\star=$\texttt{(3, 4, 3, 3, 3, 1, 2, 1, 3, 2, 1, 1, 2, 1, 1, 3, 2, 1, 2, 3, 3, 2, 1)}.

    \item[$(0.4,0.6)$] $C=221481.565$, $E=30538.421$, $S=0.879278$,
    $a^\star=$\texttt{(3, 4, 3, 3, 3, 1, 2, 1, 3, 2, 1, 1, 2, 1, 1, 3, 2, 1, 2, 3, 3, 2, 1)}.

    \item[$(0.3,0.7)$] $C=224549.428$, $E=30244.205$, $S=0.908189$,
    $a^\star=$\texttt{(3, 4, 3, 3, 3, 1, 2, 1, 3, 2, 1, 2, 2, 1, 1, 3, 2, 1, 2, 3, 3, 2, 1)}.

    \item[$(0.2,0.8)$] $C=227885.181$, $E=29934.910$, $S=0.936333$,
    $a^\star=$\texttt{(3, 4, 3, 3, 3, 1, 2, 1, 3, 2, 1, 2, 2, 1, 1, 3, 2, 1, 2, 3, 3, 2, 2)}.

    \item[$(0.1,0.9)$] $C=227885.181$, $E=29934.910$, $S=0.970666$,
    $a^\star=$\texttt{(3, 4, 3, 3, 3, 1, 2, 1, 3, 2, 1, 2, 2, 1, 1, 3, 2, 1, 2, 3, 3, 2, 2)}.

    \item[$(0,1)$] $C=227885.181$, $E=29934.910$, $S=1.000000$,
    $a^\star=$\texttt{(3, 4, 3, 3, 3, 1, 2, 1, 3, 2, 1, 2, 2, 1, 1, 3, 2, 1, 2, 3, 3, 2, 2)}.

\end{description}
}

\paragraph{Pareto-sequence solution.}
The successive one-group changes are:

\begin{center}
\small
\begin{tabular}{c l r r r}
\hline
Step & Change & $\Delta C$ & $\Delta E$ & $\rho=\Delta E/\Delta C$ \\
\hline
1 & $i=3:\ t_{3}=1\to 3$ & 23.542 & 101.062 & 4.293 \\
2 & $i=4:\ t_{4}=2\to 3$ & 523.084 & 1408.442 & 2.693 \\
3 & $i=19:\ t_{19}=3\to 2$ & 209.221 & 508.990 & 2.433 \\
4 & $i=7:\ t_{7}=3\to 2$ & 1148.311 & 2598.822 & 2.263 \\
5 & $i=10:\ t_{10}=3\to 2$ & 583.366 & 334.330 & 0.573 \\
6 & $i=13:\ t_{13}=1\to 2$ & 110.490 & 42.520 & 0.385 \\
7 & $i=17:\ t_{17}=3\to 2$ & 1324.571 & 500.356 & 0.378 \\
8 & $i=18:\ t_{18}=2\to 1$ & 1195.745 & 220.000 & 0.184 \\
9 & $i=12:\ t_{12}=1\to 2$ & 3067.863 & 294.216 & 0.096 \\
10 & $i=23:\ t_{23}=1\to 2$ & 3335.753 & 309.296 & 0.093 \\
\hline
\end{tabular}
\end{center}

\paragraph{Complete list of Pareto-sequence solutions.}
The corresponding solution vectors are:

{\footnotesize
\begin{description}
    \item[Step 0] $C=216363.233$, $E=36252.943$,
    $a^{(0)}=$\texttt{(3, 4, 1, 2, 3, 1, 3, 1, 3, 3, 1, 1, 1, 1, 1, 3, 3, 2, 3, 3, 3, 2, 1)}.

    \item[Step 1] $C=216386.775$, $E=36151.881$,
    $a^{(1)}=$\texttt{(3, 4, 3, 2, 3, 1, 3, 1, 3, 3, 1, 1, 1, 1, 1, 3, 3, 2, 3, 3, 3, 2, 1)}.

    \item[Step 2] $C=216909.860$, $E=34743.440$,
    $a^{(2)}=$\texttt{(3, 4, 3, 3, 3, 1, 3, 1, 3, 3, 1, 1, 1, 1, 1, 3, 3, 2, 3, 3, 3, 2, 1)}.

    \item[Step 3] $C=217119.081$, $E=34234.450$,
    $a^{(3)}=$\texttt{(3, 4, 3, 3, 3, 1, 3, 1, 3, 3, 1, 1, 1, 1, 1, 3, 3, 2, 2, 3, 3, 2, 1)}.

    \item[Step 4] $C=218267.392$, $E=31635.628$,
    $a^{(4)}=$\texttt{(3, 4, 3, 3, 3, 1, 2, 1, 3, 3, 1, 1, 1, 1, 1, 3, 3, 2, 2, 3, 3, 2, 1)}.

    \item[Step 5] $C=218850.758$, $E=31301.297$,
    $a^{(5)}=$\texttt{(3, 4, 3, 3, 3, 1, 2, 1, 3, 2, 1, 1, 1, 1, 1, 3, 3, 2, 2, 3, 3, 2, 1)}.

    \item[Step 6] $C=218961.249$, $E=31258.777$,
    $a^{(6)}=$\texttt{(3, 4, 3, 3, 3, 1, 2, 1, 3, 2, 1, 1, 2, 1, 1, 3, 3, 2, 2, 3, 3, 2, 1)}.

    \item[Step 7] $C=220285.820$, $E=30758.421$,
    $a^{(7)}=$\texttt{(3, 4, 3, 3, 3, 1, 2, 1, 3, 2, 1, 1, 2, 1, 1, 3, 2, 2, 2, 3, 3, 2, 1)}.

    \item[Step 8] $C=221481.565$, $E=30538.421$,
    $a^{(8)}=$\texttt{(3, 4, 3, 3, 3, 1, 2, 1, 3, 2, 1, 1, 2, 1, 1, 3, 2, 1, 2, 3, 3, 2, 1)}.

    \item[Step 9] $C=224549.428$, $E=30244.205$,
    $a^{(9)}=$\texttt{(3, 4, 3, 3, 3, 1, 2, 1, 3, 2, 1, 2, 2, 1, 1, 3, 2, 1, 2, 3, 3, 2, 1)}.

    \item[Step 10] $C=227885.181$, $E=29934.910$,
    $a^{(10)}=$\texttt{(3, 4, 3, 3, 3, 1, 2, 1, 3, 2, 1, 2, 2, 1, 1, 3, 2, 1, 2, 3, 3, 2, 2)}.

\end{description}
}

\subsection*{Building type D1}

\paragraph{TOPSIS solutions for a weight sweep with $\varepsilon=0.1$.}
The results are ordered from the cost-only to the emission-only weighting.

{\footnotesize
\begin{description}
    \item[$(1,0)$] $C=242729.133$, $E=41483.101$, $S=1.000000$,
    $a^\star=$\texttt{(3, 4, 1, 2, 3, 1, 3, 1, 3, 3, 1, 1, 1, 1, 1, 3, 3, 2, 3, 3, 3, 2, 1)}.

    \item[$(0.9,0.1)$] $C=244878.917$, $E=36268.646$, $S=0.819878$,
    $a^\star=$\texttt{(3, 4, 3, 3, 3, 1, 2, 1, 3, 3, 1, 1, 1, 1, 1, 3, 3, 2, 2, 3, 3, 2, 1)}.

    \item[$(0.8,0.2)$] $C=244878.917$, $E=36268.646$, $S=0.787397$,
    $a^\star=$\texttt{(3, 4, 3, 3, 3, 1, 2, 1, 3, 3, 1, 1, 1, 1, 1, 3, 3, 2, 2, 3, 3, 2, 1)}.

    \item[$(0.7,0.3)$] $C=246881.344$, $E=35379.807$, $S=0.794986$,
    $a^\star=$\texttt{(3, 4, 3, 3, 3, 1, 2, 1, 3, 2, 1, 1, 1, 1, 1, 3, 2, 2, 2, 3, 3, 2, 1)}.

    \item[$(0.6,0.4)$] $C=247002.124$, $E=35333.327$, $S=0.826878$,
    $a^\star=$\texttt{(3, 4, 3, 3, 3, 1, 2, 1, 3, 2, 1, 1, 2, 1, 1, 3, 2, 2, 2, 3, 3, 2, 1)}.

    \item[$(0.5,0.5)$] $C=248197.869$, $E=35113.327$, $S=0.854685$,
    $a^\star=$\texttt{(3, 4, 3, 3, 3, 1, 2, 1, 3, 2, 1, 1, 2, 1, 1, 3, 2, 1, 2, 3, 3, 2, 1)}.

    \item[$(0.4,0.6)$] $C=248197.869$, $E=35113.327$, $S=0.875909$,
    $a^\star=$\texttt{(3, 4, 3, 3, 3, 1, 2, 1, 3, 2, 1, 1, 2, 1, 1, 3, 2, 1, 2, 3, 3, 2, 1)}.

    \item[$(0.3,0.7)$] $C=251832.171$, $E=34764.788$, $S=0.905619$,
    $a^\star=$\texttt{(3, 4, 3, 3, 3, 1, 2, 1, 3, 2, 1, 2, 2, 1, 1, 3, 2, 1, 2, 3, 3, 2, 1)}.

    \item[$(0.2,0.8)$] $C=255835.075$, $E=34393.633$, $S=0.934225$,
    $a^\star=$\texttt{(3, 4, 3, 3, 3, 1, 2, 1, 3, 2, 1, 2, 2, 1, 1, 3, 2, 1, 2, 3, 3, 2, 2)}.

    \item[$(0.1,0.9)$] $C=255835.075$, $E=34393.633$, $S=0.969658$,
    $a^\star=$\texttt{(3, 4, 3, 3, 3, 1, 2, 1, 3, 2, 1, 2, 2, 1, 1, 3, 2, 1, 2, 3, 3, 2, 2)}.

    \item[$(0,1)$] $C=255835.075$, $E=34393.633$, $S=1.000000$,
    $a^\star=$\texttt{(3, 4, 3, 3, 3, 1, 2, 1, 3, 2, 1, 2, 2, 1, 1, 3, 2, 1, 2, 3, 3, 2, 2)}.

\end{description}
}

\paragraph{Pareto-sequence solution.}
The successive one-group changes are:

\begin{center}
\small
\begin{tabular}{c l r r r}
\hline
Step & Change & $\Delta C$ & $\Delta E$ & $\rho=\Delta E/\Delta C$ \\
\hline
1 & $i=3:\ t_{3}=1\to 3$ & 27.944 & 119.959 & 4.293 \\
2 & $i=4:\ t_{4}=2\to 3$ & 581.811 & 1566.566 & 2.693 \\
3 & $i=19:\ t_{19}=3\to 2$ & 251.066 & 610.788 & 2.433 \\
4 & $i=7:\ t_{7}=3\to 2$ & 1288.964 & 2917.141 & 2.263 \\
5 & $i=10:\ t_{10}=3\to 2$ & 677.856 & 388.483 & 0.573 \\
6 & $i=13:\ t_{13}=1\to 2$ & 120.781 & 46.480 & 0.385 \\
7 & $i=17:\ t_{17}=3\to 2$ & 1324.571 & 500.356 & 0.378 \\
8 & $i=18:\ t_{18}=2\to 1$ & 1195.745 & 220.000 & 0.184 \\
9 & $i=12:\ t_{12}=1\to 2$ & 3634.302 & 348.539 & 0.096 \\
10 & $i=23:\ t_{23}=1\to 2$ & 4002.904 & 371.155 & 0.093 \\
\hline
\end{tabular}
\end{center}

\paragraph{Complete list of Pareto-sequence solutions.}
The corresponding solution vectors are:

{\footnotesize
\begin{description}
    \item[Step 0] $C=242729.133$, $E=41483.101$,
    $a^{(0)}=$\texttt{(3, 4, 1, 2, 3, 1, 3, 1, 3, 3, 1, 1, 1, 1, 1, 3, 3, 2, 3, 3, 3, 2, 1)}.

    \item[Step 1] $C=242757.077$, $E=41363.142$,
    $a^{(1)}=$\texttt{(3, 4, 3, 2, 3, 1, 3, 1, 3, 3, 1, 1, 1, 1, 1, 3, 3, 2, 3, 3, 3, 2, 1)}.

    \item[Step 2] $C=243338.887$, $E=39796.576$,
    $a^{(2)}=$\texttt{(3, 4, 3, 3, 3, 1, 3, 1, 3, 3, 1, 1, 1, 1, 1, 3, 3, 2, 3, 3, 3, 2, 1)}.

    \item[Step 3] $C=243589.953$, $E=39185.787$,
    $a^{(3)}=$\texttt{(3, 4, 3, 3, 3, 1, 3, 1, 3, 3, 1, 1, 1, 1, 1, 3, 3, 2, 2, 3, 3, 2, 1)}.

    \item[Step 4] $C=244878.917$, $E=36268.646$,
    $a^{(4)}=$\texttt{(3, 4, 3, 3, 3, 1, 2, 1, 3, 3, 1, 1, 1, 1, 1, 3, 3, 2, 2, 3, 3, 2, 1)}.

    \item[Step 5] $C=245556.772$, $E=35880.163$,
    $a^{(5)}=$\texttt{(3, 4, 3, 3, 3, 1, 2, 1, 3, 2, 1, 1, 1, 1, 1, 3, 3, 2, 2, 3, 3, 2, 1)}.

    \item[Step 6] $C=245677.553$, $E=35833.683$,
    $a^{(6)}=$\texttt{(3, 4, 3, 3, 3, 1, 2, 1, 3, 2, 1, 1, 2, 1, 1, 3, 3, 2, 2, 3, 3, 2, 1)}.

    \item[Step 7] $C=247002.124$, $E=35333.327$,
    $a^{(7)}=$\texttt{(3, 4, 3, 3, 3, 1, 2, 1, 3, 2, 1, 1, 2, 1, 1, 3, 2, 2, 2, 3, 3, 2, 1)}.

    \item[Step 8] $C=248197.869$, $E=35113.327$,
    $a^{(8)}=$\texttt{(3, 4, 3, 3, 3, 1, 2, 1, 3, 2, 1, 1, 2, 1, 1, 3, 2, 1, 2, 3, 3, 2, 1)}.

    \item[Step 9] $C=251832.171$, $E=34764.788$,
    $a^{(9)}=$\texttt{(3, 4, 3, 3, 3, 1, 2, 1, 3, 2, 1, 2, 2, 1, 1, 3, 2, 1, 2, 3, 3, 2, 1)}.

    \item[Step 10] $C=255835.075$, $E=34393.633$,
    $a^{(10)}=$\texttt{(3, 4, 3, 3, 3, 1, 2, 1, 3, 2, 1, 2, 2, 1, 1, 3, 2, 1, 2, 3, 3, 2, 2)}.

\end{description}
}

\subsection*{Building type D2}

\paragraph{TOPSIS solutions for a weight sweep with $\varepsilon=0.1$.}
The results are ordered from the cost-only to the emission-only weighting.

{\footnotesize
\begin{description}
    \item[$(1,0)$] $C=240000.441$, $E=40628.959$, $S=1.000000$,
    $a^\star=$\texttt{(3, 4, 1, 2, 3, 1, 3, 1, 3, 3, 1, 1, 1, 1, 1, 3, 3, 2, 3, 3, 3, 2, 1)}.

    \item[$(0.9,0.1)$] $C=242130.922$, $E=35452.517$, $S=0.816462$,
    $a^\star=$\texttt{(3, 4, 3, 3, 3, 1, 2, 1, 3, 3, 1, 1, 1, 1, 1, 3, 3, 2, 2, 3, 3, 2, 1)}.

    \item[$(0.8,0.2)$] $C=242130.922$, $E=35452.517$, $S=0.786882$,
    $a^\star=$\texttt{(3, 4, 3, 3, 3, 1, 2, 1, 3, 3, 1, 1, 1, 1, 1, 3, 3, 2, 2, 3, 3, 2, 1)}.

    \item[$(0.7,0.3)$] $C=244078.692$, $E=34595.002$, $S=0.796719$,
    $a^\star=$\texttt{(3, 4, 3, 3, 3, 1, 2, 1, 3, 2, 1, 1, 1, 1, 1, 3, 2, 2, 2, 3, 3, 2, 1)}.

    \item[$(0.6,0.4)$] $C=244199.473$, $E=34548.522$, $S=0.829507$,
    $a^\star=$\texttt{(3, 4, 3, 3, 3, 1, 2, 1, 3, 2, 1, 1, 2, 1, 1, 3, 2, 2, 2, 3, 3, 2, 1)}.

    \item[$(0.5,0.5)$] $C=245395.218$, $E=34328.522$, $S=0.857621$,
    $a^\star=$\texttt{(3, 4, 3, 3, 3, 1, 2, 1, 3, 2, 1, 1, 2, 1, 1, 3, 2, 1, 2, 3, 3, 2, 1)}.

    \item[$(0.4,0.6)$] $C=245395.218$, $E=34328.522$, $S=0.879252$,
    $a^\star=$\texttt{(3, 4, 3, 3, 3, 1, 2, 1, 3, 2, 1, 1, 2, 1, 1, 3, 2, 1, 2, 3, 3, 2, 1)}.

    \item[$(0.3,0.7)$] $C=248627.932$, $E=34018.496$, $S=0.908012$,
    $a^\star=$\texttt{(3, 4, 3, 3, 3, 1, 2, 1, 3, 2, 1, 2, 2, 1, 1, 3, 2, 1, 2, 3, 3, 2, 1)}.

    \item[$(0.2,0.8)$] $C=252630.836$, $E=33647.341$, $S=0.936102$,
    $a^\star=$\texttt{(3, 4, 3, 3, 3, 1, 2, 1, 3, 2, 1, 2, 2, 1, 1, 3, 2, 1, 2, 3, 3, 2, 2)}.

    \item[$(0.1,0.9)$] $C=252630.836$, $E=33647.341$, $S=0.970556$,
    $a^\star=$\texttt{(3, 4, 3, 3, 3, 1, 2, 1, 3, 2, 1, 2, 2, 1, 1, 3, 2, 1, 2, 3, 3, 2, 2)}.

    \item[$(0,1)$] $C=252630.836$, $E=33647.341$, $S=1.000000$,
    $a^\star=$\texttt{(3, 4, 3, 3, 3, 1, 2, 1, 3, 2, 1, 2, 2, 1, 1, 3, 2, 1, 2, 3, 3, 2, 2)}.

\end{description}
}

\paragraph{Pareto-sequence solution.}
The successive one-group changes are:

\begin{center}
\small
\begin{tabular}{c l r r r}
\hline
Step & Change & $\Delta C$ & $\Delta E$ & $\rho=\Delta E/\Delta C$ \\
\hline
1 & $i=3:\ t_{3}=1\to 3$ & 27.777 & 119.241 & 4.293 \\
2 & $i=4:\ t_{4}=2\to 3$ & 595.811 & 1604.262 & 2.693 \\
3 & $i=19:\ t_{19}=3\to 2$ & 251.066 & 610.788 & 2.433 \\
4 & $i=7:\ t_{7}=3\to 2$ & 1255.828 & 2842.150 & 2.263 \\
5 & $i=10:\ t_{10}=3\to 2$ & 623.199 & 357.159 & 0.573 \\
6 & $i=13:\ t_{13}=1\to 2$ & 120.781 & 46.480 & 0.385 \\
7 & $i=17:\ t_{17}=3\to 2$ & 1324.571 & 500.356 & 0.378 \\
8 & $i=18:\ t_{18}=2\to 1$ & 1195.745 & 220.000 & 0.184 \\
9 & $i=12:\ t_{12}=1\to 2$ & 3232.714 & 310.026 & 0.096 \\
10 & $i=23:\ t_{23}=1\to 2$ & 4002.904 & 371.155 & 0.093 \\
\hline
\end{tabular}
\end{center}

\paragraph{Complete list of Pareto-sequence solutions.}
The corresponding solution vectors are:

{\footnotesize
\begin{description}
    \item[Step 0] $C=240000.441$, $E=40628.959$,
    $a^{(0)}=$\texttt{(3, 4, 1, 2, 3, 1, 3, 1, 3, 3, 1, 1, 1, 1, 1, 3, 3, 2, 3, 3, 3, 2, 1)}.

    \item[Step 1] $C=240028.217$, $E=40509.718$,
    $a^{(1)}=$\texttt{(3, 4, 3, 2, 3, 1, 3, 1, 3, 3, 1, 1, 1, 1, 1, 3, 3, 2, 3, 3, 3, 2, 1)}.

    \item[Step 2] $C=240624.028$, $E=38905.455$,
    $a^{(2)}=$\texttt{(3, 4, 3, 3, 3, 1, 3, 1, 3, 3, 1, 1, 1, 1, 1, 3, 3, 2, 3, 3, 3, 2, 1)}.

    \item[Step 3] $C=240875.094$, $E=38294.667$,
    $a^{(3)}=$\texttt{(3, 4, 3, 3, 3, 1, 3, 1, 3, 3, 1, 1, 1, 1, 1, 3, 3, 2, 2, 3, 3, 2, 1)}.

    \item[Step 4] $C=242130.922$, $E=35452.517$,
    $a^{(4)}=$\texttt{(3, 4, 3, 3, 3, 1, 2, 1, 3, 3, 1, 1, 1, 1, 1, 3, 3, 2, 2, 3, 3, 2, 1)}.

    \item[Step 5] $C=242754.121$, $E=35095.358$,
    $a^{(5)}=$\texttt{(3, 4, 3, 3, 3, 1, 2, 1, 3, 2, 1, 1, 1, 1, 1, 3, 3, 2, 2, 3, 3, 2, 1)}.

    \item[Step 6] $C=242874.902$, $E=35048.878$,
    $a^{(6)}=$\texttt{(3, 4, 3, 3, 3, 1, 2, 1, 3, 2, 1, 1, 2, 1, 1, 3, 3, 2, 2, 3, 3, 2, 1)}.

    \item[Step 7] $C=244199.473$, $E=34548.522$,
    $a^{(7)}=$\texttt{(3, 4, 3, 3, 3, 1, 2, 1, 3, 2, 1, 1, 2, 1, 1, 3, 2, 2, 2, 3, 3, 2, 1)}.

    \item[Step 8] $C=245395.218$, $E=34328.522$,
    $a^{(8)}=$\texttt{(3, 4, 3, 3, 3, 1, 2, 1, 3, 2, 1, 1, 2, 1, 1, 3, 2, 1, 2, 3, 3, 2, 1)}.

    \item[Step 9] $C=248627.932$, $E=34018.496$,
    $a^{(9)}=$\texttt{(3, 4, 3, 3, 3, 1, 2, 1, 3, 2, 1, 2, 2, 1, 1, 3, 2, 1, 2, 3, 3, 2, 1)}.

    \item[Step 10] $C=252630.836$, $E=33647.341$,
    $a^{(10)}=$\texttt{(3, 4, 3, 3, 3, 1, 2, 1, 3, 2, 1, 2, 2, 1, 1, 3, 2, 1, 2, 3, 3, 2, 2)}.

\end{description}
}

\subsection*{Building type D3}

\paragraph{TOPSIS solutions for a weight sweep with $\varepsilon=0.1$.}
The results are ordered from the cost-only to the emission-only weighting.

{\footnotesize
\begin{description}
    \item[$(1,0)$] $C=232325.669$, $E=40244.961$, $S=1.000000$,
    $a^\star=$\texttt{(3, 4, 1, 2, 3, 1, 3, 1, 3, 3, 1, 1, 1, 1, 1, 3, 3, 2, 3, 3, 3, 2, 1)}.

    \item[$(0.9,0.1)$] $C=234423.742$, $E=35173.092$, $S=0.817993$,
    $a^\star=$\texttt{(3, 4, 3, 3, 3, 1, 2, 1, 3, 3, 1, 1, 1, 1, 1, 3, 3, 2, 2, 3, 3, 2, 1)}.

    \item[$(0.8,0.2)$] $C=234423.742$, $E=35173.092$, $S=0.787238$,
    $a^\star=$\texttt{(3, 4, 3, 3, 3, 1, 2, 1, 3, 3, 1, 1, 1, 1, 1, 3, 3, 2, 2, 3, 3, 2, 1)}.

    \item[$(0.7,0.3)$] $C=236375.004$, $E=34313.576$, $S=0.793503$,
    $a^\star=$\texttt{(3, 4, 3, 3, 3, 1, 2, 1, 3, 2, 1, 1, 1, 1, 1, 3, 2, 2, 2, 3, 3, 2, 1)}.

    \item[$(0.6,0.4)$] $C=236495.785$, $E=34267.096$, $S=0.826588$,
    $a^\star=$\texttt{(3, 4, 3, 3, 3, 1, 2, 1, 3, 2, 1, 1, 2, 1, 1, 3, 2, 2, 2, 3, 3, 2, 1)}.

    \item[$(0.5,0.5)$] $C=237691.529$, $E=34047.096$, $S=0.854684$,
    $a^\star=$\texttt{(3, 4, 3, 3, 3, 1, 2, 1, 3, 2, 1, 1, 2, 1, 1, 3, 2, 1, 2, 3, 3, 2, 1)}.

    \item[$(0.4,0.6)$] $C=237691.529$, $E=34047.096$, $S=0.877113$,
    $a^\star=$\texttt{(3, 4, 3, 3, 3, 1, 2, 1, 3, 2, 1, 1, 2, 1, 1, 3, 2, 1, 2, 3, 3, 2, 1)}.

    \item[$(0.3,0.7)$] $C=240891.321$, $E=33740.228$, $S=0.905735$,
    $a^\star=$\texttt{(3, 4, 3, 3, 3, 1, 2, 1, 3, 2, 1, 2, 2, 1, 1, 3, 2, 1, 2, 3, 3, 2, 1)}.

    \item[$(0.2,0.8)$] $C=244894.225$, $E=33369.073$, $S=0.934088$,
    $a^\star=$\texttt{(3, 4, 3, 3, 3, 1, 2, 1, 3, 2, 1, 2, 2, 1, 1, 3, 2, 1, 2, 3, 3, 2, 2)}.

    \item[$(0.1,0.9)$] $C=244894.225$, $E=33369.073$, $S=0.969592$,
    $a^\star=$\texttt{(3, 4, 3, 3, 3, 1, 2, 1, 3, 2, 1, 2, 2, 1, 1, 3, 2, 1, 2, 3, 3, 2, 2)}.

    \item[$(0,1)$] $C=244894.225$, $E=33369.073$, $S=1.000000$,
    $a^\star=$\texttt{(3, 4, 3, 3, 3, 1, 2, 1, 3, 2, 1, 2, 2, 1, 1, 3, 2, 1, 2, 3, 3, 2, 2)}.

\end{description}
}

\paragraph{Pareto-sequence solution.}
The successive one-group changes are:

\begin{center}
\small
\begin{tabular}{c l r r r}
\hline
Step & Change & $\Delta C$ & $\Delta E$ & $\rho=\Delta E/\Delta C$ \\
\hline
1 & $i=3:\ t_{3}=1\to 3$ & 27.777 & 119.241 & 4.293 \\
2 & $i=4:\ t_{4}=2\to 3$ & 523.084 & 1408.442 & 2.693 \\
3 & $i=19:\ t_{19}=3\to 2$ & 251.066 & 610.788 & 2.433 \\
4 & $i=7:\ t_{7}=3\to 2$ & 1296.147 & 2933.398 & 2.263 \\
5 & $i=10:\ t_{10}=3\to 2$ & 626.690 & 359.160 & 0.573 \\
6 & $i=13:\ t_{13}=1\to 2$ & 120.781 & 46.480 & 0.385 \\
7 & $i=17:\ t_{17}=3\to 2$ & 1324.571 & 500.356 & 0.378 \\
8 & $i=18:\ t_{18}=2\to 1$ & 1195.745 & 220.000 & 0.184 \\
9 & $i=12:\ t_{12}=1\to 2$ & 3199.791 & 306.868 & 0.096 \\
10 & $i=23:\ t_{23}=1\to 2$ & 4002.904 & 371.155 & 0.093 \\
\hline
\end{tabular}
\end{center}

\paragraph{Complete list of Pareto-sequence solutions.}
The corresponding solution vectors are:

{\footnotesize
\begin{description}
    \item[Step 0] $C=232325.669$, $E=40244.961$,
    $a^{(0)}=$\texttt{(3, 4, 1, 2, 3, 1, 3, 1, 3, 3, 1, 1, 1, 1, 1, 3, 3, 2, 3, 3, 3, 2, 1)}.

    \item[Step 1] $C=232353.445$, $E=40125.720$,
    $a^{(1)}=$\texttt{(3, 4, 3, 2, 3, 1, 3, 1, 3, 3, 1, 1, 1, 1, 1, 3, 3, 2, 3, 3, 3, 2, 1)}.

    \item[Step 2] $C=232876.530$, $E=38717.278$,
    $a^{(2)}=$\texttt{(3, 4, 3, 3, 3, 1, 3, 1, 3, 3, 1, 1, 1, 1, 1, 3, 3, 2, 3, 3, 3, 2, 1)}.

    \item[Step 3] $C=233127.595$, $E=38106.490$,
    $a^{(3)}=$\texttt{(3, 4, 3, 3, 3, 1, 3, 1, 3, 3, 1, 1, 1, 1, 1, 3, 3, 2, 2, 3, 3, 2, 1)}.

    \item[Step 4] $C=234423.742$, $E=35173.092$,
    $a^{(4)}=$\texttt{(3, 4, 3, 3, 3, 1, 2, 1, 3, 3, 1, 1, 1, 1, 1, 3, 3, 2, 2, 3, 3, 2, 1)}.

    \item[Step 5] $C=235050.433$, $E=34813.932$,
    $a^{(5)}=$\texttt{(3, 4, 3, 3, 3, 1, 2, 1, 3, 2, 1, 1, 1, 1, 1, 3, 3, 2, 2, 3, 3, 2, 1)}.

    \item[Step 6] $C=235171.213$, $E=34767.452$,
    $a^{(6)}=$\texttt{(3, 4, 3, 3, 3, 1, 2, 1, 3, 2, 1, 1, 2, 1, 1, 3, 3, 2, 2, 3, 3, 2, 1)}.

    \item[Step 7] $C=236495.785$, $E=34267.096$,
    $a^{(7)}=$\texttt{(3, 4, 3, 3, 3, 1, 2, 1, 3, 2, 1, 1, 2, 1, 1, 3, 2, 2, 2, 3, 3, 2, 1)}.

    \item[Step 8] $C=237691.529$, $E=34047.096$,
    $a^{(8)}=$\texttt{(3, 4, 3, 3, 3, 1, 2, 1, 3, 2, 1, 1, 2, 1, 1, 3, 2, 1, 2, 3, 3, 2, 1)}.

    \item[Step 9] $C=240891.321$, $E=33740.228$,
    $a^{(9)}=$\texttt{(3, 4, 3, 3, 3, 1, 2, 1, 3, 2, 1, 2, 2, 1, 1, 3, 2, 1, 2, 3, 3, 2, 1)}.

    \item[Step 10] $C=244894.225$, $E=33369.073$,
    $a^{(10)}=$\texttt{(3, 4, 3, 3, 3, 1, 2, 1, 3, 2, 1, 2, 2, 1, 1, 3, 2, 1, 2, 3, 3, 2, 2)}.

\end{description}
}

\subsection*{Building type D4}

\paragraph{TOPSIS solutions for a weight sweep with $\varepsilon=0.1$.}
The results are ordered from the cost-only to the emission-only weighting.

{\footnotesize
\begin{description}
    \item[$(1,0)$] $C=231164.936$, $E=39237.067$, $S=1.000000$,
    $a^\star=$\texttt{(3, 4, 1, 2, 3, 1, 3, 1, 3, 3, 1, 1, 1, 1, 1, 3, 3, 2, 3, 3, 3, 2, 1)}.

    \item[$(0.9,0.1)$] $C=233199.970$, $E=34322.644$, $S=0.825526$,
    $a^\star=$\texttt{(3, 4, 3, 3, 3, 1, 2, 1, 3, 3, 1, 1, 1, 1, 1, 3, 3, 2, 2, 3, 3, 2, 1)}.

    \item[$(0.8,0.2)$] $C=233886.582$, $E=33929.142$, $S=0.787154$,
    $a^\star=$\texttt{(3, 4, 3, 3, 3, 1, 2, 1, 3, 2, 1, 1, 1, 1, 1, 3, 3, 2, 2, 3, 3, 2, 1)}.

    \item[$(0.7,0.3)$] $C=235211.154$, $E=33428.786$, $S=0.792157$,
    $a^\star=$\texttt{(3, 4, 3, 3, 3, 1, 2, 1, 3, 2, 1, 1, 1, 1, 1, 3, 2, 2, 2, 3, 3, 2, 1)}.

    \item[$(0.6,0.4)$] $C=235331.934$, $E=33382.306$, $S=0.822175$,
    $a^\star=$\texttt{(3, 4, 3, 3, 3, 1, 2, 1, 3, 2, 1, 1, 2, 1, 1, 3, 2, 2, 2, 3, 3, 2, 1)}.

    \item[$(0.5,0.5)$] $C=236527.679$, $E=33162.306$, $S=0.850065$,
    $a^\star=$\texttt{(3, 4, 3, 3, 3, 1, 2, 1, 3, 2, 1, 1, 2, 1, 1, 3, 2, 1, 2, 3, 3, 2, 1)}.

    \item[$(0.4,0.6)$] $C=236527.679$, $E=33162.306$, $S=0.870715$,
    $a^\star=$\texttt{(3, 4, 3, 3, 3, 1, 2, 1, 3, 2, 1, 1, 2, 1, 1, 3, 2, 1, 2, 3, 3, 2, 1)}.

    \item[$(0.3,0.7)$] $C=240292.014$, $E=32801.297$, $S=0.902329$,
    $a^\star=$\texttt{(3, 4, 3, 3, 3, 1, 2, 1, 3, 2, 1, 2, 2, 1, 1, 3, 2, 1, 2, 3, 3, 2, 1)}.

    \item[$(0.2,0.8)$] $C=244294.918$, $E=32430.142$, $S=0.932175$,
    $a^\star=$\texttt{(3, 4, 3, 3, 3, 1, 2, 1, 3, 2, 1, 2, 2, 1, 1, 3, 2, 1, 2, 3, 3, 2, 2)}.

    \item[$(0.1,0.9)$] $C=244294.918$, $E=32430.142$, $S=0.968675$,
    $a^\star=$\texttt{(3, 4, 3, 3, 3, 1, 2, 1, 3, 2, 1, 2, 2, 1, 1, 3, 2, 1, 2, 3, 3, 2, 2)}.

    \item[$(0,1)$] $C=244294.918$, $E=32430.142$, $S=1.000000$,
    $a^\star=$\texttt{(3, 4, 3, 3, 3, 1, 2, 1, 3, 2, 1, 2, 2, 1, 1, 3, 2, 1, 2, 3, 3, 2, 2)}.

\end{description}
}

\paragraph{Pareto-sequence solution.}
The successive one-group changes are:

\begin{center}
\small
\begin{tabular}{c l r r r}
\hline
Step & Change & $\Delta C$ & $\Delta E$ & $\rho=\Delta E/\Delta C$ \\
\hline
1 & $i=3:\ t_{3}=1\to 3$ & 23.458 & 100.703 & 4.293 \\
2 & $i=4:\ t_{4}=2\to 3$ & 509.084 & 1370.746 & 2.693 \\
3 & $i=19:\ t_{19}=3\to 2$ & 251.066 & 610.788 & 2.433 \\
4 & $i=7:\ t_{7}=3\to 2$ & 1251.425 & 2832.186 & 2.263 \\
5 & $i=10:\ t_{10}=3\to 2$ & 686.612 & 393.502 & 0.573 \\
6 & $i=13:\ t_{13}=1\to 2$ & 120.781 & 46.480 & 0.385 \\
7 & $i=17:\ t_{17}=3\to 2$ & 1324.571 & 500.356 & 0.378 \\
8 & $i=18:\ t_{18}=2\to 1$ & 1195.745 & 220.000 & 0.184 \\
9 & $i=12:\ t_{12}=1\to 2$ & 3764.335 & 361.009 & 0.096 \\
10 & $i=23:\ t_{23}=1\to 2$ & 4002.904 & 371.155 & 0.093 \\
\hline
\end{tabular}
\end{center}

\paragraph{Complete list of Pareto-sequence solutions.}
The corresponding solution vectors are:

{\footnotesize
\begin{description}
    \item[Step 0] $C=231164.936$, $E=39237.067$,
    $a^{(0)}=$\texttt{(3, 4, 1, 2, 3, 1, 3, 1, 3, 3, 1, 1, 1, 1, 1, 3, 3, 2, 3, 3, 3, 2, 1)}.

    \item[Step 1] $C=231188.395$, $E=39136.364$,
    $a^{(1)}=$\texttt{(3, 4, 3, 2, 3, 1, 3, 1, 3, 3, 1, 1, 1, 1, 1, 3, 3, 2, 3, 3, 3, 2, 1)}.

    \item[Step 2] $C=231697.479$, $E=37765.618$,
    $a^{(2)}=$\texttt{(3, 4, 3, 3, 3, 1, 3, 1, 3, 3, 1, 1, 1, 1, 1, 3, 3, 2, 3, 3, 3, 2, 1)}.

    \item[Step 3] $C=231948.545$, $E=37154.830$,
    $a^{(3)}=$\texttt{(3, 4, 3, 3, 3, 1, 3, 1, 3, 3, 1, 1, 1, 1, 1, 3, 3, 2, 2, 3, 3, 2, 1)}.

    \item[Step 4] $C=233199.970$, $E=34322.644$,
    $a^{(4)}=$\texttt{(3, 4, 3, 3, 3, 1, 2, 1, 3, 3, 1, 1, 1, 1, 1, 3, 3, 2, 2, 3, 3, 2, 1)}.

    \item[Step 5] $C=233886.582$, $E=33929.142$,
    $a^{(5)}=$\texttt{(3, 4, 3, 3, 3, 1, 2, 1, 3, 2, 1, 1, 1, 1, 1, 3, 3, 2, 2, 3, 3, 2, 1)}.

    \item[Step 6] $C=234007.363$, $E=33882.662$,
    $a^{(6)}=$\texttt{(3, 4, 3, 3, 3, 1, 2, 1, 3, 2, 1, 1, 2, 1, 1, 3, 3, 2, 2, 3, 3, 2, 1)}.

    \item[Step 7] $C=235331.934$, $E=33382.306$,
    $a^{(7)}=$\texttt{(3, 4, 3, 3, 3, 1, 2, 1, 3, 2, 1, 1, 2, 1, 1, 3, 2, 2, 2, 3, 3, 2, 1)}.

    \item[Step 8] $C=236527.679$, $E=33162.306$,
    $a^{(8)}=$\texttt{(3, 4, 3, 3, 3, 1, 2, 1, 3, 2, 1, 1, 2, 1, 1, 3, 2, 1, 2, 3, 3, 2, 1)}.

    \item[Step 9] $C=240292.014$, $E=32801.297$,
    $a^{(9)}=$\texttt{(3, 4, 3, 3, 3, 1, 2, 1, 3, 2, 1, 2, 2, 1, 1, 3, 2, 1, 2, 3, 3, 2, 1)}.

    \item[Step 10] $C=244294.918$, $E=32430.142$,
    $a^{(10)}=$\texttt{(3, 4, 3, 3, 3, 1, 2, 1, 3, 2, 1, 2, 2, 1, 1, 3, 2, 1, 2, 3, 3, 2, 2)}.

\end{description}
}

\subsection*{Building type K1}

\paragraph{TOPSIS solutions for a weight sweep with $\varepsilon=0.1$.}
The results are ordered from the cost-only to the emission-only weighting.

{\footnotesize
\begin{description}
    \item[$(1,0)$] $C=97645.270$, $E=16135.556$, $S=1.000000$,
    $a^\star=$\texttt{(3, 4, 1, 2, 3, 1, 3, 1, 3, 3, 1, 1, 1, 1, 1, 3, 3, 1, 3, 3, 3, 2, 1)}.

    \item[$(0.9,0.1)$] $C=98608.643$, $E=13789.688$, $S=0.781889$,
    $a^\star=$\texttt{(3, 4, 3, 3, 3, 1, 2, 1, 3, 3, 1, 1, 1, 1, 1, 3, 3, 1, 2, 3, 3, 2, 1)}.

    \item[$(0.8,0.2)$] $C=98608.643$, $E=13789.688$, $S=0.783984$,
    $a^\star=$\texttt{(3, 4, 3, 3, 3, 1, 2, 1, 3, 3, 1, 1, 1, 1, 1, 3, 3, 1, 2, 3, 3, 2, 1)}.

    \item[$(0.7,0.3)$] $C=98884.055$, $E=13637.679$, $S=0.808542$,
    $a^\star=$\texttt{(3, 4, 3, 3, 3, 1, 2, 1, 3, 2, 1, 1, 2, 1, 1, 3, 3, 1, 2, 3, 3, 2, 1)}.

    \item[$(0.6,0.4)$] $C=99546.341$, $E=13387.501$, $S=0.853732$,
    $a^\star=$\texttt{(3, 4, 3, 3, 3, 1, 2, 1, 3, 2, 1, 1, 2, 1, 1, 3, 2, 1, 2, 3, 3, 2, 1)}.

    \item[$(0.5,0.5)$] $C=99546.341$, $E=13387.501$, $S=0.883723$,
    $a^\star=$\texttt{(3, 4, 3, 3, 3, 1, 2, 1, 3, 2, 1, 1, 2, 1, 1, 3, 2, 1, 2, 3, 3, 2, 1)}.

    \item[$(0.4,0.6)$] $C=99546.341$, $E=13387.501$, $S=0.902051$,
    $a^\star=$\texttt{(3, 4, 3, 3, 3, 1, 2, 1, 3, 2, 1, 1, 2, 1, 1, 3, 2, 1, 2, 3, 3, 2, 1)}.

    \item[$(0.3,0.7)$] $C=100703.733$, $E=13276.504$, $S=0.925919$,
    $a^\star=$\texttt{(3, 4, 3, 3, 3, 1, 2, 1, 3, 2, 1, 2, 2, 1, 1, 3, 2, 1, 2, 3, 3, 2, 1)}.

    \item[$(0.2,0.8)$] $C=102007.059$, $E=13164.706$, $S=0.948713$,
    $a^\star=$\texttt{(3, 4, 3, 3, 3, 1, 2, 1, 3, 2, 1, 2, 1, 1, 1, 3, 2, 1, 2, 3, 3, 2, 2)}.

    \item[$(0.1,0.9)$] $C=102038.034$, $E=13152.786$, $S=0.976530$,
    $a^\star=$\texttt{(3, 4, 3, 3, 3, 1, 2, 1, 3, 2, 1, 2, 2, 1, 1, 3, 2, 1, 2, 3, 3, 2, 2)}.

    \item[$(0,1)$] $C=102038.034$, $E=13152.786$, $S=1.000000$,
    $a^\star=$\texttt{(3, 4, 3, 3, 3, 1, 2, 1, 3, 2, 1, 2, 2, 1, 1, 3, 2, 1, 2, 3, 3, 2, 2)}.

\end{description}
}

\paragraph{Pareto-sequencev solution.}
The successive one-group changes are:

\begin{center}
\small
\begin{tabular}{c l r r r}
\hline
Step & Change & $\Delta C$ & $\Delta E$ & $\rho=\Delta E/\Delta C$ \\
\hline
1 & $i=3:\ t_{3}=1\to 3$ & 10.085 & 43.295 & 4.293 \\
2 & $i=4:\ t_{4}=2\to 3$ & 304.905 & 820.979 & 2.693 \\
3 & $i=19:\ t_{19}=3\to 2$ & 83.689 & 203.596 & 2.433 \\
4 & $i=7:\ t_{7}=3\to 2$ & 564.694 & 1277.997 & 2.263 \\
5 & $i=10:\ t_{10}=3\to 2$ & 244.438 & 140.089 & 0.573 \\
6 & $i=13:\ t_{13}=1\to 2$ & 30.975 & 11.920 & 0.385 \\
7 & $i=17:\ t_{17}=3\to 2$ & 662.286 & 250.178 & 0.378 \\
8 & $i=12:\ t_{12}=1\to 2$ & 1157.391 & 110.997 & 0.096 \\
9 & $i=23:\ t_{23}=1\to 2$ & 1334.301 & 123.718 & 0.093 \\
\hline
\end{tabular}
\end{center}

\paragraph{Complete list of Pareto-sequence solutions.}
The corresponding solution vectors are:

{\footnotesize
\begin{description}
    \item[Step 0] $C=97645.270$, $E=16135.556$,
    $a^{(0)}=$\texttt{(3, 4, 1, 2, 3, 1, 3, 1, 3, 3, 1, 1, 1, 1, 1, 3, 3, 1, 3, 3, 3, 2, 1)}.

    \item[Step 1] $C=97655.355$, $E=16092.261$,
    $a^{(1)}=$\texttt{(3, 4, 3, 2, 3, 1, 3, 1, 3, 3, 1, 1, 1, 1, 1, 3, 3, 1, 3, 3, 3, 2, 1)}.

    \item[Step 2] $C=97960.260$, $E=15271.281$,
    $a^{(2)}=$\texttt{(3, 4, 3, 3, 3, 1, 3, 1, 3, 3, 1, 1, 1, 1, 1, 3, 3, 1, 3, 3, 3, 2, 1)}.

    \item[Step 3] $C=98043.949$, $E=15067.685$,
    $a^{(3)}=$\texttt{(3, 4, 3, 3, 3, 1, 3, 1, 3, 3, 1, 1, 1, 1, 1, 3, 3, 1, 2, 3, 3, 2, 1)}.

    \item[Step 4] $C=98608.643$, $E=13789.688$,
    $a^{(4)}=$\texttt{(3, 4, 3, 3, 3, 1, 2, 1, 3, 3, 1, 1, 1, 1, 1, 3, 3, 1, 2, 3, 3, 2, 1)}.

    \item[Step 5] $C=98853.081$, $E=13649.599$,
    $a^{(5)}=$\texttt{(3, 4, 3, 3, 3, 1, 2, 1, 3, 2, 1, 1, 1, 1, 1, 3, 3, 1, 2, 3, 3, 2, 1)}.

    \item[Step 6] $C=98884.055$, $E=13637.679$,
    $a^{(6)}=$\texttt{(3, 4, 3, 3, 3, 1, 2, 1, 3, 2, 1, 1, 2, 1, 1, 3, 3, 1, 2, 3, 3, 2, 1)}.

    \item[Step 7] $C=99546.341$, $E=13387.501$,
    $a^{(7)}=$\texttt{(3, 4, 3, 3, 3, 1, 2, 1, 3, 2, 1, 1, 2, 1, 1, 3, 2, 1, 2, 3, 3, 2, 1)}.

    \item[Step 8] $C=100703.733$, $E=13276.504$,
    $a^{(8)}=$\texttt{(3, 4, 3, 3, 3, 1, 2, 1, 3, 2, 1, 2, 2, 1, 1, 3, 2, 1, 2, 3, 3, 2, 1)}.

    \item[Step 9] $C=102038.034$, $E=13152.786$,
    $a^{(9)}=$\texttt{(3, 4, 3, 3, 3, 1, 2, 1, 3, 2, 1, 2, 2, 1, 1, 3, 2, 1, 2, 3, 3, 2, 2)}.

\end{description}
}

\subsection*{Building type K2}

\paragraph{TOPSIS solutions for a weight sweep with $\varepsilon=0.1$.}
The results are ordered from the cost-only to the emission-only weighting.

{\footnotesize
\begin{description}
    \item[$(1,0)$] $C=154085.144$, $E=27375.199$, $S=1.000000$,
    $a^\star=$\texttt{(3, 4, 1, 2, 3, 1, 3, 1, 3, 3, 1, 1, 1, 1, 1, 3, 3, 1, 3, 3, 3, 2, 1)}.

    \item[$(0.9,0.1)$] $C=155247.399$, $E=23188.930$, $S=0.824997$,
    $a^\star=$\texttt{(3, 4, 3, 3, 3, 1, 2, 1, 3, 3, 1, 1, 1, 1, 1, 3, 3, 1, 2, 3, 3, 2, 1)}.

    \item[$(0.8,0.2)$] $C=155247.399$, $E=23188.930$, $S=0.823629$,
    $a^\star=$\texttt{(3, 4, 3, 3, 3, 1, 2, 1, 3, 3, 1, 1, 1, 1, 1, 3, 3, 1, 2, 3, 3, 2, 1)}.

    \item[$(0.7,0.3)$] $C=155684.784$, $E=22947.968$, $S=0.843969$,
    $a^\star=$\texttt{(3, 4, 3, 3, 3, 1, 2, 1, 3, 2, 1, 1, 2, 1, 1, 3, 3, 1, 2, 3, 3, 2, 1)}.

    \item[$(0.6,0.4)$] $C=156347.069$, $E=22697.790$, $S=0.876142$,
    $a^\star=$\texttt{(3, 4, 3, 3, 3, 1, 2, 1, 3, 2, 1, 1, 2, 1, 1, 3, 2, 1, 2, 3, 3, 2, 1)}.

    \item[$(0.5,0.5)$] $C=156347.069$, $E=22697.790$, $S=0.896278$,
    $a^\star=$\texttt{(3, 4, 3, 3, 3, 1, 2, 1, 3, 2, 1, 1, 2, 1, 1, 3, 2, 1, 2, 3, 3, 2, 1)}.

    \item[$(0.4,0.6)$] $C=158089.611$, $E=22530.676$, $S=0.908471$,
    $a^\star=$\texttt{(3, 4, 3, 3, 3, 1, 2, 1, 3, 2, 1, 2, 2, 1, 1, 3, 2, 1, 2, 3, 3, 2, 1)}.

    \item[$(0.3,0.7)$] $C=158089.611$, $E=22530.676$, $S=0.929620$,
    $a^\star=$\texttt{(3, 4, 3, 3, 3, 1, 2, 1, 3, 2, 1, 2, 2, 1, 1, 3, 2, 1, 2, 3, 3, 2, 1)}.

    \item[$(0.2,0.8)$] $C=160706.658$, $E=22303.080$, $S=0.950843$,
    $a^\star=$\texttt{(3, 4, 3, 3, 3, 1, 2, 1, 3, 2, 1, 2, 1, 1, 1, 3, 2, 1, 2, 3, 3, 2, 2)}.

    \item[$(0.1,0.9)$] $C=160758.214$, $E=22283.240$, $S=0.977517$,
    $a^\star=$\texttt{(3, 4, 3, 3, 3, 1, 2, 1, 3, 2, 1, 2, 2, 1, 1, 3, 2, 1, 2, 3, 3, 2, 2)}.

    \item[$(0,1)$] $C=160758.214$, $E=22283.240$, $S=1.000000$,
    $a^\star=$\texttt{(3, 4, 3, 3, 3, 1, 2, 1, 3, 2, 1, 2, 2, 1, 1, 3, 2, 1, 2, 3, 3, 2, 2)}.

\end{description}
}

\paragraph{Pareto-sequence solution.}
The successive one-group changes are:

\begin{center}
\small
\begin{tabular}{c l r r r}
\hline
Step & Change & $\Delta C$ & $\Delta E$ & $\rho=\Delta E/\Delta C$ \\
\hline
1 & $i=4:\ t_{4}=2\to 3$ & 32.053 & 1566.566 & 48.875 \\
2 & $i=3:\ t_{3}=1\to 3$ & 16.493 & 70.803 & 4.293 \\
3 & $i=19:\ t_{19}=3\to 2$ & 167.377 & 407.192 & 2.433 \\
4 & $i=7:\ t_{7}=3\to 2$ & 946.332 & 2141.708 & 2.263 \\
5 & $i=10:\ t_{10}=3\to 2$ & 385.830 & 221.121 & 0.573 \\
6 & $i=13:\ t_{13}=1\to 2$ & 51.555 & 19.840 & 0.385 \\
7 & $i=17:\ t_{17}=3\to 2$ & 662.286 & 250.178 & 0.378 \\
8 & $i=12:\ t_{12}=1\to 2$ & 1742.542 & 167.114 & 0.096 \\
9 & $i=23:\ t_{23}=1\to 2$ & 2668.603 & 247.436 & 0.093 \\
\hline
\end{tabular}
\end{center}

\paragraph{Complete list of Pareto-sequence solutions.}
The corresponding solution vectors are:

{\footnotesize
\begin{description}
    \item[Step 0] $C=154085.144$, $E=27375.199$,
    $a^{(0)}=$\texttt{(3, 4, 1, 2, 3, 1, 3, 1, 3, 3, 1, 1, 1, 1, 1, 3, 3, 1, 3, 3, 3, 2, 1)}.

    \item[Step 1] $C=154117.197$, $E=25808.633$,
    $a^{(1)}=$\texttt{(3, 4, 1, 3, 3, 1, 3, 1, 3, 3, 1, 1, 1, 1, 1, 3, 3, 1, 3, 3, 3, 2, 1)}.

    \item[Step 2] $C=154133.690$, $E=25737.829$,
    $a^{(2)}=$\texttt{(3, 4, 3, 3, 3, 1, 3, 1, 3, 3, 1, 1, 1, 1, 1, 3, 3, 1, 3, 3, 3, 2, 1)}.

    \item[Step 3] $C=154301.067$, $E=25330.637$,
    $a^{(3)}=$\texttt{(3, 4, 3, 3, 3, 1, 3, 1, 3, 3, 1, 1, 1, 1, 1, 3, 3, 1, 2, 3, 3, 2, 1)}.

    \item[Step 4] $C=155247.399$, $E=23188.930$,
    $a^{(4)}=$\texttt{(3, 4, 3, 3, 3, 1, 2, 1, 3, 3, 1, 1, 1, 1, 1, 3, 3, 1, 2, 3, 3, 2, 1)}.

    \item[Step 5] $C=155633.228$, $E=22967.808$,
    $a^{(5)}=$\texttt{(3, 4, 3, 3, 3, 1, 2, 1, 3, 2, 1, 1, 1, 1, 1, 3, 3, 1, 2, 3, 3, 2, 1)}.

    \item[Step 6] $C=155684.784$, $E=22947.968$,
    $a^{(6)}=$\texttt{(3, 4, 3, 3, 3, 1, 2, 1, 3, 2, 1, 1, 2, 1, 1, 3, 3, 1, 2, 3, 3, 2, 1)}.

    \item[Step 7] $C=156347.069$, $E=22697.790$,
    $a^{(7)}=$\texttt{(3, 4, 3, 3, 3, 1, 2, 1, 3, 2, 1, 1, 2, 1, 1, 3, 2, 1, 2, 3, 3, 2, 1)}.

    \item[Step 8] $C=158089.611$, $E=22530.676$,
    $a^{(8)}=$\texttt{(3, 4, 3, 3, 3, 1, 2, 1, 3, 2, 1, 2, 2, 1, 1, 3, 2, 1, 2, 3, 3, 2, 1)}.

    \item[Step 9] $C=160758.214$, $E=22283.240$,
    $a^{(9)}=$\texttt{(3, 4, 3, 3, 3, 1, 2, 1, 3, 2, 1, 2, 2, 1, 1, 3, 2, 1, 2, 3, 3, 2, 2)}.

\end{description}
}

\end{appendices}

\end{document}